# ESPACES SYMETRIQUES SYMPLECTIQUES

# Table des matières

# Introduction

Les espaces symétriques ont surtout été étudiés dans le cadre de la géométrie riemannienne ou pseudo-riemannienne, c'est à dire dans un cadre métrique. Il est dès lors tentant de définir une notion d'espace symétrique dans un cadre purement symplectique : *un espace symétrique symplectique est une variété symétrique munie d'une structure symplectique invariante par les symétries.*
Mon travail constitue une première étude de ces espaces. Les méthodes utilisées sont celles de la géométrie différentielle classique, de la théorie des groupes et algèbres de Lie, de la géométrie symplectique. La notion d'espace symétrique symplectique se place donc à un point de confluence de diverses théories; ce qui lui confère, à mon sens, un certain intérêt.

Le travail est divisé en six chapitres.
Dans le premier, on trouve une définition précise d'espace symétrique symplectique et les premiers outils permettant l'étude de ces espaces. Dans le second on trouve un théorème de décomposition "à la de Rham" des espaces symétriques symplectiques . Le troisième chapitre constitue une étude assez complète des espaces symétriques symplectiques dont le groupe des automorphismes est semisimple. Certains résultats concernant les espaces symétriques symplectiques dont le groupe des automorphismes n'est ni résoluble ni semisimple sont énoncés dans le chapitre 4. Des exemples d'espaces symétriques symplectiques et notement une classification en dimension 2 ou 4 sont exposés dans les cinquième et sixième chapitres.
Le plan général étant donné, je précise maintenant le contenu des différents chapitres.

**Chapitre 1. Généralités**

Dans son livre "*Symmetric spaces*" (1969) ([Lo]), O. Loos définit la notion d'espace symétrique de la manière suivante.
Un *espace symétrique* est un couple $(M, s)$ où $M$ est une variété différentiable connexe et $s : M \times M \to M$ est une application différentiable telle qu'en notant, pour tout $x$ dans $M$, $s_x = s(x, .)$ on ait

(i) $s_x$ est un difféomorphisme involutif de $M$ admettant $x$ comme point fixe isolé.

(ii) $s_x s_y s_x = s_{s_x(y)}$ pour tous $x, y$ dans $M$.

Il démontre l'existence et l'unicité d'une dérivée covariante affine $\nabla$ sur $M$ telle que les "symétries" ($\{s_x ; x \in M\}$) soient des transformations affines de $(M, \nabla)$; $(M, \nabla)$ est alors une *variété affine symétrique* ([Ko-No]). Dès lors, dans le cas simplement connexe, ces espaces peuvent être décrits de manière totalement algébrique ([Lo],[Ko-No]).

Dans un cadre de géométrie symplectique, il est naturel de définir la notion d'espace symétrique symplectique .
Un *espace symétrique symplectique* (cf. Définition 1.1.) est un triple $(M, \omega, s)$ où $(M, s)$ est un espace symétrique et où $\omega$ est une forme symplectique sur $M$ invariante par les symétries.

Dans ce cadre, l'existence et l'unicité de la dérivée covariante de Loos s'obtiennent très rapidement par des arguments élémentaires (cf. pages 5-7); on peut même écrire une formule explicite pour cette dérivée covariante (cf. Lemme 1.6.) qui est également valable dans un cadre pseudo-riemannien. On l'appelle *dérivée covariante canonique*; elle est symplectique (cf. Lemme 1.7.).

Deux espaces symétriques symplectiques $(M, \omega, s)$ et $(M', \omega', s')$ sont dits *isomorphes* s'il existe un difféomorphisme symplectique entre $(M, \omega)$ et $(M', \omega')$ qui "entrelace" $s$ et $s'$ (cf. Définition 1.2). Si $(M, \omega, s) = (M', \omega', s')$ on parle d'automorphisme. Le groupe $Aut(M, \omega, s)$ des automorphismes d'un espace symétrique symplectique $(M, \omega, s)$ est l'intersection du groupe des transformations affines de $(M, \nabla)$ et du groupe des difféomorphismes symplectiques de $(M, \omega)$ (cf. Proposition 1.8.); c'est un groupe de Lie de transformations de $M$ transitif sur $M$ (cf.page 8). Ceci permet une description des espaces symétriques symplectiques en termes d'espaces homogènes réductifs symplectiques (cf. Proposition 1.12.).

La notion de triple symétrique symplectique définit l'objet algébrique qui encode la structure géométrique d'un espace symétrique symplectique simplement connexe.

Un *triple symétrique symplectique* (cf. Définition 1.19.) est un triple $t = (\mathcal{G}, \sigma, \underline{\Omega})$ où

(i) $\mathcal{G}$ est une algèbre de Lie réelle de dimension finie.

(ii) $\sigma$ est un automorphisme involutif de $\mathcal{G}$ tel que si $\mathcal{G} = \mathcal{K} \oplus \mathcal{P}$ ($\sigma = id_\mathcal{K} \oplus (-id_\mathcal{P})$) on ait $[\mathcal{P}, \mathcal{P}] = \mathcal{K}$ et l'action de $\mathcal{K}$ sur $\mathcal{P}$ soit fidèle.

(iii) $\underline{\Omega}$ est un 2-cocycle de Chevalley pour la représentation triviale de $\mathcal{G}$ sur $\mathbb{R}$ tel que sa restriction à $\mathcal{P} \times \mathcal{P}$ soit de rang maximum et tel que pour tout $X$ dans $\mathcal{K}$, $i(X)\underline{\Omega} = 0$.

On définit de manière évidente la notion d'isomorphie de triples symétriques symplectiques (cf. définition 1.19.).

Il y a une correspondance bijective entre l'ensemble des classes d'isomorphie de triples symétriques symplectiques et l'ensemble des classes d'isomorphie d'espaces symétriques symplectiques simplement connexes (cf. Proposition 1.21.).

En bref, cette correspondance s'obtient de la manière suivante.

En fixant un point $o$ de $M$, la conjugaison par $s_o$ définit un automorphisme involutif $\tilde{\sigma}$ du groupe de Lie $Aut(M, \omega, s)$. Le groupe des transvections $G(M)$ — c'est-à-dire le sous groupe de $Aut(M, \omega, s)$ engendré par les produits d'un nombre pair de symétries — est un sous groupe de Lie connexe de $Aut(M, \omega, s)$; c'est le plus petit sous groupe de $Aut(M, \omega, s)$ transitif sur $M$ et stabilisé par $\tilde{\sigma}$, en particulier $M$ est un espace homogène pour $G(M)$ (cf. Proposition 1.17.).

L'algèbre $\mathcal{G}$ apparaissant dans la définition de triple symétrique symplectique n'est autre que l'algèbre de Lie de $G(M)$ et $\sigma$ en est l'automorphisme involutif induit par $\tilde{\sigma}$. La forme symplectique $\omega$ sur $M$ se relève en une forme invariante à gauche sur $G(M)$ dont la valeur au neutre définit $\underline{\Omega}$.

L'algèbre $\mathcal{K}$ est connue pour être isomorphe à l'algèbre de Lie du groupe d'holonomie en $o$ relativement à la dérivée covariante canonique $\nabla$ ([Ko-No]).

Un espace symétrique symplectique est un espace homogène symplectique relativement à l'action de son groupe des transvections, il est dès lors naturel de demander à quelle condition cette action est fortement hamiltonienne. Dans le cas d'un espace symétrique symplectique simplement connexe $(M, \omega, s)$ cette action est fortement hamiltonienne si et seulement si le 2-cocycle $\underline{\Omega}$ est exact (cf. Propostion 1.27.).

### Chapitre 2. Décompositions

En 1952, G. de Rham énonce le résultat suivant ([deR]).

*Soit $(M, g)$ une variété riemannienne (avec une métrique définie)*

*Soit $o$ un point de $M$ et $\Phi$ l'holonomie en $o$ relativement à la connexion de Levi-Civita.*

(i) *Si $\Phi$ préserve un sous espace propre $M_o^1$ de l'espace tangent $M_o$ et si $M$ est complète, connexe et simplement connexe, alors, $M$ est isométrique au produit direct des variétés intégrales maximales pour les distributions obtenues par transport parallèle de $M_o^1$ et $M_o^{1\,\perp}$ sur $M$.*

(ii) *$M_o$ se décompose de manière unique à l'ordre des facteurs près en une somme directe orthogonale de sous espaces $\Phi$-irréductibles et du noyau de $\Phi$*

Dix ans plus tard, H.Wu obtient un résultat analogue au point (i) du théorème de de Rham pour les variétés pseudo-riemanniennes ([Wu1]):
*Soit $(M,g)$ une variété pseudo-riemanienne.*
*Soit $o$ un point de $M$ et $\Phi$ l'holonomie en $o$.*
*Soit $M_o^1$ un sous espace propre $\Phi$-stable de l'espace tangent $M_o$ sur lequel la métrique se restreint de façon non singulière.*
*Si $M$ est complète, connexe et simplement connexe, alors, $M$ est isométrique au produit direct des variétés intégrales maximales pour les distributions obtenues par transport parallèle de $M_o^1$ et $M_o^{1\,\perp}$ sur $M$.*
(Si $M_o$ admet un tel sous espace on dit que $M_o$ est faiblement réductible; faiblement irréductible sinon.)
Il remarque une condition suffisante pour l'unicité à l'ordre près d'une décomposition de $M_o$ en une somme directe orthogonale de sous espaces faiblement irréductibles (relativement à l'action de $\Phi$) et d'un sous espace sur lequel $\Phi$ agit trivialement. Cette condition est :
*La métrique restreinte au noyau de $\Phi$ est non singulière* ([Wu1]). En 1967, il prouve que cette condition est cruciale ([Wu2]).
En 1980, M. Cahen et M. Parker obtiennent un résultat complet quant à l'unicité de la décompositon de de Rham-Wu dans le cadre des espaces pseudo-riemanniens symétriques ([Ca-Pa]) :
*Soit $(M,g)$ un espace pseudoriemannien symétrique connexe et simplement connexe.*
*Soit $o$ un point de $M$ et $\Phi$ l'holonomie en $o$.*
*Soit $I(M)$ le groupe des isométries de $M$.*
*Soient*
$$M_o = \bigoplus_{\alpha=0}^{p} M_o^\alpha = \bigoplus_{\beta=0}^{q} {M'}_o^\beta$$
*deux décompositons en somme directe orthogonale telles que :*

(i) *$\Phi$ agit trivialement sur $M_o^0$ et ${M'}_o^0$.*

(ii) *Pour tous $\alpha, \beta \geq 1$, $M_o^\alpha$ et ${M'}_o^\beta$ sont faiblement irréductibles relativement à l'action non triviale de $\Phi$.*

*Soient $\{M^\alpha\}$ et $\{{M'}^\beta\}$ les variétés intégrales maximales pour les distributions obtenues par transport parallèle de $\{M_o^\alpha\}$ et $\{{M'}_o^\beta\}$ sur $M$ respectivement (ce sont des sous-espaces symétriques).*
*Alors*

(a) *$p = q$*

(b) *Il existe un élément $h$ dans $I(M)$ et une permutation $\eta$ de $\{1,...,p\}$ tels que*
$$h(M^0) = {M'}^0 \quad et \quad h(M^\alpha) = {M'}^{\eta(\alpha)} \quad \forall \alpha \in \{1,...,p\}$$

On peut dès lors faire les remarques suivantes.

1. L'opération "somme directe" est un moyen de définir la notion d'"espace indécomposable".

2. Dans le cadre des variétés riemanniennes (avec métrique définie) et dans celui des espaces pseudoriemanniens symétriques l'opération "décomposition en somme directe d'epaces indécomposables" a un sens fonctoriel: deux espaces isomorphes admettent des décompositions isomorphes (cf. point (ii) du théorème de de Rham et théorème de décomposition de M. Cahen et M. Parker).

3. La notion de faible réductibilité est une très bonne notion dans le cadre des variétés pseudo-riemanniennes: elle ramène la décomposition d'un espace à la recherche des composantes faiblement irréductibles de l'holonomie.

J'ai obtenu deux théorèmes analogues à celui de M. Cahen et M. Parker. Ces deux résultats concernent respectivement le cas des variétés affines symétriques (non nécéssairement munies d'une structure métrique ou symplectique) et celui des espaces symétriques symplectiques.

(I) (cf. Théorème 2.1.3.)
Soit $(M,\nabla)$ une variété affine symétrique connexe et simplement connexe.
Soient
$$(M,\nabla) = \bigoplus_{\alpha=0}^{p}(M^\alpha,\nabla^\alpha) = \bigoplus_{\beta=0}^{q}(M'^\beta,\nabla'^\beta)$$
deux décompositons en somme directe de variétés affines symétriques telles que :

(i) $(M^0,\nabla^0)$ et $(M'^0,\nabla'^0)$ sont plates.

(ii) Pour tous $\alpha,\beta \geq 1$, $(M^\alpha,\nabla^\alpha)$ et $(M'^\beta,\nabla'^\beta)$ sont indécomposables et non plates.

Alors

(a) $p = q$

(b) Il existe une tranformation affine $h$ de $(M,\nabla)$ et une permutation $\eta$ de $\{1,...,p\}$ telles que
$$h(M^0,\nabla^0) = (M'^0,\nabla'^0) \quad et \quad h(M^\alpha,\nabla^\alpha) = (M'^{\eta(\alpha)},\nabla'^{\eta(\alpha)}) \quad \forall \alpha \in \{1,...,p\}$$

(II) (cf. Théorèmes 2.3.2. et 2.3.5. ) Soit $(M,\omega,s)$ un espace symétrique symplectique connexe et simplement connexe.
Soient
$$(M,\omega,s) = \bigoplus_{\alpha=0}^{p}(M^\alpha,\omega^\alpha,s^\alpha) = \bigoplus_{\beta=0}^{q}(M'^\beta,\omega'^\beta,s'^\beta)$$
deux décompositons en somme directe d'espaces symétriques symplectiques telles que :

(i) $(M^0,\omega^0,s^0)$ et $(M'^0,\omega'^0,s'^0)$ sont plats.

(ii) Pour tous $\alpha,\beta \geq 1$, $(M^\alpha,\omega^\alpha,s^\alpha)$ et $(M'^\beta,\omega'^\beta,s'^\beta)$ sont indécomposables et non plats.

Alors

(a) $p = q$

(b) Il existe un automorphisme $h$ de $(M,\omega,s)$ et une permutation $\eta$ de $\{1,...,p\}$ telles que
$$h(M^0,\omega^0,s^0) = (M'^0,\omega'^0,s'^0) \quad et \quad h(M^\alpha,\omega^\alpha,s^\alpha) = (M'^{\eta(\alpha)},\omega'^{\eta(\alpha)},s'^{\eta(\alpha)})$$
$$\forall \alpha \in \{1,...,p\}$$

(c) Dans le cas où l'action du groupe des transvections $G(M)$ est fortement hamiltonienne sur $(M,\omega)$ la décomposition est unique à l'ordre des facteurs près, le facteur plat n'y apparait pas et cette décomposition coïncide avec la décomposition de la variété affine $(M,\nabla)$ décrite dans (I) (où $\nabla$ désigne la dérivée covariante canonique); en particulier cette dernière décomposition est, elle aussi, unique à l'ordre des facteurs près.

La remarque 2. reste donc valable dans ces deux cas.

Quant à la notion de faible réductibilité elle semble ne pas être adaptée à une situation symplectique . Dans le but de montrer ceci, adoptons les notations suivantes. Soit $(M, \omega, s)$ un espace symétrique symplectique simplement connexe, $o$ un point de $M$ et $\Phi$ l'holonomie en $o$ pour la dérivée covariante canonique $\nabla$. Soit $M_o^1$ un sous espace propre $\Phi$-stable symplectique dans $M_o$ et $M^1$ la variété intégrale maximale pour la distribution obtenue par transport parallèle de $M_o^1$ sur $M$. $M^1$ est alors la variété sous-jacente à un sous espace symétrique symplectique $(M^1, \omega^1, s^1)$ de $(M, \omega, s)$ (cf. Lemme 1.20.). Soit $\Phi^1$ l'holonomie en $o$ relativement à la connexion canonique $\nabla^1$ sur $M^1$. On a alors les faits suivants.

a) Si $M^{1\perp}$ désigne la variété intégrale maximale pour la distribution obtenue par transport parallèle de $M_o^{1\perp}$ sur $M$, on n'a en général PAS

$$(M, \omega, s) \simeq (M^1, \omega^1, s^1) \times (M^{1\perp}, \omega^{1\perp}, s^{1\perp})$$

(voir contre-exemples en dimension 4, chapitre 6).Ceci rend les preuves des théorèmes (I) et (II) susmentionnés substanciellement différentes de celle du théorème de décomposition de M. Cahen et M. Parker pour les espaces pseudo-riemanniens symétriques.

b) En général un sous espace symplectique $\Phi^1$-stable dans $M_o^1$ n'est pas $\Phi$-stable dans $M_o$. J'exhibe, dans l'exemple suivant, une classe d'espaces symétriques symplectiques $(M, \omega, s)$ dont aucun sous espace $M_o^1$ $\Phi$-stable dans $M_o$ n'est faiblement irréductible relativement à $\Phi^1$.
$(M, \omega, s)$ est l'espace symétrique symplectique associé au triple symétrique symplectique $t = (\mathcal{G}, \sigma, \underline{\Omega})$ défini par
$$\mathcal{K} = \rangle u, v \langle \, ; \, \mathcal{P} = \rangle e_1, e_2, e_3, f_1, f_2, f_3 \langle$$
$$\underline{\Omega}(e_i, e_j) = \underline{\Omega}(f_i, f_j) = 0 \, ; \, \underline{\Omega}(e_i, f_j) = \delta_{ij}$$

et la table de $\mathcal{G}$ est

$$\begin{aligned}
[u, e_2] &= e_1 \\
[u, f_1] &= -f_2 \\
[u, e_3] &= f_3 \\
[v, f_1] &= e_1 \\
[f_1, e_3] &= \alpha u \\
[e_2, e_3] &= \alpha v \\
[f_1, e_2] &= \beta v \\
[f_1, f_2] &= \gamma v
\end{aligned}$$

où $\alpha \in \mathbb{R}_0$ ; $\beta, \gamma \in \mathbb{R}$

c) On peut ignorer la situation décrite dans a) et b) et "décomposer" à tout prix $M_o$ en sous espaces faiblement irréductibles sans se soucier de l'invariance par l'holonomie $\Phi$. Du point de vue du triple symétrique symplectique $t = (\mathcal{G}, \sigma, \underline{\Omega})$ associé à $(M, \omega, s)$, cela revient à écrire

$$\mathcal{P} = \bigoplus_{\alpha=1}^{p} \mathcal{P}_\alpha$$

avec $\underline{\Omega}(\mathcal{P}_\alpha, \mathcal{P}_\beta) = 0$  si $\alpha \neq \beta$ et $[[\mathcal{P}_\alpha, \mathcal{P}_\alpha], \mathcal{P}_\alpha] \subset \mathcal{P}_\alpha$ pour tout $\alpha$. Une telle décomposition n'a aucun caractère d'unicité: dans l'exemple suivant, j'exhibe un triple symétrique symplectique $t = (\mathcal{G}, \sigma, \underline{\Omega})$ dont $\mathcal{P}$ admet une "décomposition" en *trois* sous espaces faiblement irréductibles et une autre "décomposition" en *deux* tels sous espaces.

$t = (\mathcal{G}, \sigma, \underline{\Omega})$ est défini par

$$\mathcal{K} = \rangle u, v \langle \, ; \, \mathcal{P} = \rangle e_1, e_2, e_3, f_1, f_2, f_3 \langle$$

$$\underline{\Omega}(e_i, e_j) = \underline{\Omega}(f_i, f_j) = 0 \, ; \, \underline{\Omega}(e_i, e_j) = \delta_{ij}$$

et la table de $\mathcal{G}$ est

$$\begin{aligned} [u, e_2] &= e_1 \\ [u, f_1] &= -f_2 \\ [v, f_1] &= e_1 \\ [f_1, e_3] &= -u \\ [e_2, e_3] &= -v \\ [f_1, f_2] &= v \end{aligned}$$

En définissant $\mathcal{P}' = \rangle \{e_3, f_3\} \langle$ on a $\mathcal{P}'^\perp = \rangle \{e_1, f_1, e_2, f_2\} \langle$ et le triple symétrique symplectique induit par $[\mathcal{P}'^\perp, \mathcal{P}'^\perp] \oplus \mathcal{P}'^\perp$ est faiblement réductible.

Par contre, en définissant $\mathcal{P}' = \rangle \{e_3 + \alpha e_1 + \beta f_2, f_3 + \gamma e_1 + \delta f_2\} \langle$ on peut trouver $\alpha, \beta, \gamma, \delta \in \mathbb{R}$ tels que $[\mathcal{P}'^\perp, \mathcal{P}'^\perp] = \mathcal{K}$ et $\mathcal{K} \oplus \mathcal{P}^\perp$ induit un triple symétrique symplectique faiblement irréductible.

**Chapitre 3. Le cas réductif**

Elie Cartan a complètement élucidé la structure des espaces symétriques riemanniens simplement connexes ([Dieu]). Il a démontré le fait fondamental suivant : *classifier ces espaces est équivalent à classifier les algèbres de Lie simples réelles*.

L'idée de la démonstration est la suivante.

Soit $(M, g)$ un espace symétrique riemannien. Soit $(\mathcal{G}^\wedge, \sigma^\wedge, B^\wedge)$ le triple de Riemann associé. La métrique $g$ étant définie positive, l'holonomie est compacte; ceci impose à $\mathcal{G}^\wedge$ d'être réductive (cf. Lemme 3.1.2.). Le triple de Riemann s'écrit alors comme une somme directe

$$(\mathcal{G}^\wedge, \sigma^\wedge, B^\wedge) = (\mathcal{Z}, -id, B_\mathcal{Z}) \oplus (\mathcal{G}, \sigma, B)$$

où $(\mathcal{Z}, -id, B_\mathcal{Z})$ est un facteur euclidien et où $\mathcal{G}$ est semisimple. A nouveau comme l'holonomie est compacte, la sous algèbre $\mathcal{K}$ des point fixes de $\sigma$ dans $\mathcal{G}$ est compacte et dès lors $\sigma$ est "une involution de Cartan" de $\mathcal{G}$; elle est donc déterminée, à automorphisme intérieur de $\mathcal{G}$ près, par la structure d'algèbre de Lie sur $\mathcal{G}$.

Notons qu'en 1914, Cartan a, pour la première fois, donné une classification des algèbres de Lie simples réelles. Quinze ans plus tard, il donne une autre méthode de classification de ces algèbres :

*Etant donné une algèbre de Lie complexe simple $\mathcal{G}$, il y a une correspondance bijective entre l'ensemble des classes d'isomorphie de ses formes réelles et l'ensemble des classes de conjugaison des automorphismes involutifs de sa forme compacte $\mathcal{G}_u$.*

C'est cette correspondance qui induit la dualité entre les espaces hermitiens symétriques compacts et les espaces hermitiens symétriques non compacts.

Ces automorphismes ont été classifiés, au moyens de diverses méthodes, par Cartan (1929), Gantmacher (1939), Borel et de Siebenthal (1949), Araki (1962), Murakami (1965) et Wallach (1965) ([Bo-deSi], [Mu1]).

L'intersection des espaces symétriques symplectiques et des espaces symétriques riemanniens est constituée des espaces symétriques hermitiens. Ces espaces ont largement été étudiés par divers auteurs : Bergman, Chow, Koecher, Koranyi, Shapiro, Wolf, ...

Shapiro a également étudié les espaces pseudo-hermitiens symétriques ([Sha]).

Dans ce chapitre, je donne une description complète des espaces symétriques symplectiques simplement connexes $(M, \omega, s)$ dont l'holonomie $\Phi$ en un point $o$ de $M$ agit de manière complètement réductible sur $M_o$. Les triples symétriques symplectiques correspondant sont les triples symétriques symplectiques réductifs — c'est à dire les triples symétriques symplectiques $(\mathcal{G}, \sigma, \underline{\Omega})$ où $\mathcal{G}$ est une algèbre réductive (cf. Lemme 3.1.2.). Un triple symétrique symplectique réductif indécomposable et non plat est simple (cf. Proposition 3.5.4.). Dès lors, classifier les espaces symétriques symplectiques simplement connexes dont le groupe des transvections est réductif revient à donner un théorème analogue au théorème de classifcation des espaces riemanniens symétriques.

Par des méthodes analogues à celles utilisées par Borel, de Siebenthal et Murakami ([Bo-deSi], [Mu1]) et en utilisant certains résultats dûs à Sugiura, j'ai obtenu un tel théorème (cf. Théorème 3.7.1.). L'intérêt de ce résultat est qu'il décrit les triples symétriques symplectiques simples et leur isomorphie en termes des propriétés des racines de $\mathcal{G}$. Dans le cas symplectique, en général $\sigma$ n'est pas une involution de Cartan de $\mathcal{G}$; en particulier la notion de dualité "compact-non compact" disparait dans le cadre symplectique.

Dans le but d'obtenir une liste des triples symétriques symplectiques , on peut utiliser le critère suivant (cf. Théorème 3.1.8.) pour décider quels éléments dans la liste de Berger ([Be]) des espaces symétriques irréductibles livrent des espaces symétriques symplectiques .

Soit un couple symétrique $(\mathcal{G}, \sigma)$ (cf. Définition 2.1.1.) où $\mathcal{G}$ est simple et où $\sigma$ est automorphisme involutif de $\mathcal{G}$.
Alors il existe une bijection entre les 2-cocycles $\underline{\Omega} \in \mathcal{Z}^2(\mathcal{G})$ tels que $(\mathcal{G}, \sigma, \underline{\Omega})$ est un triple symétrique symplectique et les éléments non nuls du centre $\mathcal{Z}(\mathcal{K})$ de $\mathcal{K}$ ($\sigma = id_{\mathcal{K}} \oplus (-id_{\mathcal{P}})$).

Remarquons qu'un critère différent mais équivalent à celui-ci avait déjà été obtenu par M. Parker ([Pa]). L'algèbre de Lie $\mathcal{K}$ de l'holonomie d'un espace symétrique dont le groupe des transvections est semisimple est réductive dans $\mathcal{G}$ ([Che]); dans le cas des espaces symétriques symplectiques , elle a même rang que $\mathcal{G}$ (cf. Théorème 3.2.5.).

Un espace symétrique symplectique simplement connexe $(M, \omega, s)$ dont le groupe $G$ des transvections est simple est le revêtement universel symplectique d'une orbite coadjointe $\theta$ de $G$ dans $\mathcal{G}^\star$ — en effet les lemmes de Withehead nous disent que $\underline{\Omega}$ est exacte. Une telle orbite est symétrique.

La structure locale de l'espace homogène $G \to \theta$ étant identique à celle de l'espace homogène $G \to M$, caractériser et classifier les orbites symétriques de $G$ dans $\mathcal{G}^\star$ revient à décrire complètement les espaces symétriques symplectiques simplement connexes simples (et par là les réductifs). Je passe maintenant à l'énoncé du théorème, pour ce faire j'adopte les définitions et notations suivantes.

Soient $\mathcal{G}$ une algèbre de Lie simple complexe et $\boldsymbol{h}$ une sous-algèbre de Cartan de $\mathcal{G}$.

Soit $\phi$ le système de racines correspondant à $\boldsymbol{h}$ et $\mathcal{G} = \boldsymbol{h} \oplus \bigoplus_{\beta \in \phi} \mathcal{G}_\beta$ la décomposition radicielle de $\mathcal{G}$ par rapport à $\boldsymbol{h}$.

Un système admissible est un couple $\{\alpha, \Delta\}$ où $\Delta$ est une base de $\phi$ et $\alpha$ est un élément de $\Delta$ tels que si

$\mu \in \phi$ est la racine maximale relativement à $\Delta$, on ait :

$$\mu = \alpha + \sum_{\alpha' \in \Delta \setminus \{\alpha\}} n_{\alpha'} \cdot \alpha'$$

où $n_{\alpha'} \in \mathbb{N}$.

On note $h_\alpha^\Delta$ l'élément de la base duale (relativement à la forme de Killing) de $\Delta$ caractérisé par $\alpha(h_\alpha^\Delta) = 1$ et $\omega_\alpha^\Delta$ le poid fondemental associé à $\alpha$ relativement à $\Delta$ — on considère $\omega_\alpha^\Delta$ comme un élément de $\mathcal{G}^\star$ c'est à dire $\omega_\alpha^\Delta(\mathcal{G}_\beta) = 0 \quad \forall \beta \in \phi$.

Soient $\tau$ une conjugaison de $\mathcal{G}$ qui stabilise $\boldsymbol{h}$ et $\mathcal{G}_\tau$ la forme réelle de $\mathcal{G}$ associée; on note encore $\omega_\alpha^\Delta$ la resriction de $\omega_\alpha^\Delta$ à $\mathcal{G}_\tau$.

On a alors (cf. sections 3.3, 3.4, et 3.5)

(i) Une orbite coadjointe dans $\mathcal{G}^\star$ est symétrique si et seulement si elle contient un multiple (complexe) d'un poid fondemental $\omega_\alpha^\Delta$ où $\{\alpha, \Delta\}$ est un système admissible.
Le triple symétrique symplectique associé est

$$(\mathcal{G}, exp\ \pi\ i\ \mathrm{ad}(h_\alpha^\Delta), Re(\lambda.\delta\omega_\alpha^\Delta))$$

où $\lambda \in \mathbb{C}^\star$ ($\delta$ désigne le cobord de Chevalley).

(ii) Une orbite coadjointe dans $\mathcal{G}_\tau^\star$ est symétrique si et seulement si elle contient

a) un multiple imaginaire pur d'un poid fondemental resreint $\omega_\alpha^\Delta$ tel que $\tau^\star(\omega_\alpha^\Delta) = -\omega_\alpha^\Delta$ où $\{\alpha, \Delta\}$ est un système admissible. Dans ce cas, le triple symétrique symplectique associé est

$$(\mathcal{G}_\tau, exp\ \pi\ i\ \mathrm{ad}(h_\alpha^\Delta), r\,\delta(\,i\,\omega_\alpha^\Delta))$$

où $r \in \mathbb{R}_0$.

*ou bien*

b) un multiple réel d'un poid fondemental resreint $\omega_\alpha^\Delta$ tel que $\tau^\star(\omega_\alpha^\Delta) = \omega_\alpha^\Delta$ où $\{\alpha, \Delta\}$ est un système admissible. Dans ce cas, le triple symétrique symplectique associé est

$$(\mathcal{G}_\tau, exp\ \pi\ i\ \mathrm{ad}(h_\alpha^\Delta), r\delta\omega_\alpha^\Delta)$$

où $r \in \mathbb{R}_0$.

(iii) Les espaces symétriques symplectiques correspondant aux cas (i) et (ii) a) sont des espaces symétriques pseudo-kähleriens, la structure complexe au point $o$ est donnée par

$$\mathcal{J}_o = \pm \left. \frac{1}{\pi} \log \sigma \right|_\mathcal{P} \qquad (\mathcal{P} \sim M_o)$$

Les espaces symétriques symplectiques correspondant au cas (ii) b) ne sont jamais pseudo-kähleriens.

Remarquons que des multiples différents en valeur absolue donnent lieu à des triples symétriques symplectiques non isomorphes.

Soit $N(\boldsymbol{h})$ le groupe des automorphismes de $\mathcal{G}$ qui stabilisent $\boldsymbol{h}$.
Soit $\widetilde{N}_\tau$ le sous groupe de $N(\boldsymbol{h})$ constitué des éléments qui commutent avec $\tau$.
On note $Aut(\phi)$ le groupe des isométries de $\boldsymbol{h}_\mathbb{R}^\star$ qui conservent $\phi$.
On a l'épimorphisme canonique $\pi : N(\boldsymbol{h}) \to Aut(\phi)$; on note $N_\tau$ l'image de $\widetilde{N}_\tau$ par $\pi$.
On a alors (cf. section 3.6)

1) Les couples symétriques $(\mathcal{G}, exp\ \pi\ i\ \text{ad}(h_\alpha^\Delta))$ et $(\mathcal{G}, exp\ \pi\ i\ \text{ad}(h_{\alpha'}^{\Delta'}))$ sous jacents à deux triples symétriques symplectiques simples complexes (point (i) susmentionné) sont isomorphes si et seulement si les poids fondamentaux $\omega_\alpha^\Delta$ et $\omega_{\alpha'}^{\Delta'}$ sont dans la même orbite pour l'action de $Aut(\phi)$ sur $\boldsymbol{h}_\mathbb{R}^\star$.

2.a) Un choix de $\boldsymbol{h}$ tel que $\mathcal{G}_\tau \cap \boldsymbol{h}$ admet une partie compacte de dimension maximale livre un ensemble exhaustif des classes d'isomorphie des couples symétriques sous jacents aux triples symétriques symplectiques pseudo-kähleriens simples absolument simples. Dans cette situation, les couples symétriques $(\mathcal{G}_\tau, exp\ \pi\ i\ \text{ad}(h_\alpha^\Delta))$ et $(\mathcal{G}_\tau, exp\ \pi\ i\ \text{ad}(h_{\alpha'}^{\Delta'}))$ sous jacents à deux triples symétriques symplectiques pseudo-kähleriens simples absolument simples (point (ii) a) susmentionné) sont isomorphes si et seulement si les poids fondamentaux $\omega_\alpha^\Delta$ et $\omega_{\alpha'}^{\Delta'}$ sont dans la même orbite pour l'action de $N_\tau$ sur $\boldsymbol{h}_\mathbb{R}^\star$.

2.b) Un choix de $\boldsymbol{h}$ tel que $\mathcal{G}_\tau \cap \boldsymbol{h}$ admet une partie non compacte de dimension maximale livre un ensemble exhaustif des classes d'isomorphie des couples symétriques sous jacents aux triples symétriques symplectiques non pseudo-kähleriens simples absolument simples. Dans cette situation les couples symétriques $(\mathcal{G}_\tau, exp\ \pi\ i\ \text{ad}(h_\alpha^\Delta))$ et $(\mathcal{G}_\tau, exp\ \pi\ i\ \text{ad}(h_{\alpha'}^{\Delta'}))$ sous jacents à deux triples symétriques symplectiques non pseudo-kähleriens simples absolument simples (point (ii) b) susmentionné) sont isomorphes si et seulement si les poids fondamentaux $\omega_\alpha^\Delta$ et $\omega_{\alpha'}^{\Delta'}$ sont dans la même orbite pour l'action de $N_\tau$ sur $\boldsymbol{h}_\mathbb{R}^\star$.

Remarquons qu'un espace symétrique symplectique simplement connexe dont le groupe des transvections est simple est pseudo-kählerien si et seulement si le centre de l'holonomie contient un élément compact non trivial. Remarquons encore qu'en général, contrairement au cas hermitien symétrique, un espace symétrique symplectique n'est pas simplement connexe — l'hyperboloïde $SL(2,\mathbb{R})/SO(1,1)$ n'est pas simplement connexe.

En utilisant un résultat dû à Koh ([Koh]), on obtient la classification des espaces symétriques symplectiques compacts

1. Tout espace symétrique symplectique compact est kählerien.

2. Un espace symétrique symplectique compact indécomposable et non plat est simple; donc simplement connexe.

Les groupes de Lie sont des espaces symétriques, il est dès lors naturel de demander lesquels admettent une structure d'espace symétrique symplectique.
Soient $G$ un groupe de Lie d'algèbre $\mathcal{G}$ et $\nabla$ la dérivée covariante canonique sur $G$ déduite de sa structure d'espace symétrique. Le groupe $G^2 = G \times G$ est le groupe des transvections de $(G, \nabla)$, l'action étant

$$G^2 \times G \to G : ((g, h), x) \to gxh^{-1}$$

Dès lors, toute forme symplectique $G^2$-invariante $\omega$ sur $G$ est $G$-biinvariante. En posant $\Omega = \omega_e$ on a, pour tous $X, Y, Z$ dans $\mathcal{G}$ :

$$\Omega([X,Y],Z) = \Omega([X,Z],Y) = \Omega([Y,Z],X) = \Omega([Y,X],Z) = 0$$

Donc $[\mathcal{G}, \mathcal{G}] = 0$ par non singularité de $\Omega$; et on a :

Tout groupe de Lie symétrique symplectique est abélien.

**Chapitre 4. Résultats relatifs aux T.S.S. ni résolubles, ni semi-simples**

Soit $(M, \omega, s)$ un espace symétrique symplectique simplement connexe dont le groupe des transvections n'est ni résoluble, ni semi-simple. Soit $t = (\mathcal{G}, \sigma, \underline{\Omega})$ le triple symétrique symplectique associé.
les résultats principaux de ce chapitre sont les suivants.

(I) $t$ admet un unique facteur semisimple maximal — $t_{\mathcal{L}} = (\mathcal{L}, \sigma|_{\mathcal{L}}, \underline{\Omega}_{\mathcal{L}})$ est un facteur semisimple de $t$ s'il existe un triple symétrique symplectique $t_1$ tel qu'on ait la décomposition $t = t_{\mathcal{L}} \oplus t_1$ (cf. Proposition 4.1.2.).

(II) Soit $t_{\mathcal{L}}$ l'unique facteur semisimple maximal de $t$. Soient $\mathcal{R}$ le radical de $\mathcal{G}$ et $\mathcal{P}_{\mathcal{R}}$ les éléments antifixés par $\sigma$ dans $\mathcal{R}$. Alors,

  (i) $\mathcal{P}_{\mathcal{R}}$ est symplectique si et seulement si $\mathcal{L}$ est un facteur de Levi de $\mathcal{G}$ (cf. Proposition 4.1.6.).

  (ii) On a
  $$2 \leqslant \dim \mathcal{P}_{\mathcal{R}} \leqslant \dim M - 2$$
  (cf. Proposition 4.1.8.).

(III) Supposons $t$ sans facteur semisimple et $\mathcal{P}_{\mathcal{R}}$ isotrope. Alors,

  (i) $\mathcal{R}$ est abélien (cf. Lemme 4.2.1.).

  (ii) on a une constante $C$ qui ne dépend que du facteur de Levi de $t$ telle que
  $$\dim \mathcal{P}_{\mathcal{R}} \geq C$$
  — pour plus de précisions voir la Proposition 4.2.11..

(IV) (i) Si $t$ est sans facteur semisimple alors $\dim \mathcal{P}_{\mathcal{R}} = 2$ si et seulement si $\dim M = 4$ (cf. Théorème 4.2.13.).

  (ii) Ceci livre une classification des triples symétriques symplectiques indécomposables avec $\dim \mathcal{P}_{\mathcal{R}} = 2$. Les espaces symétriques symplectiques correspondants sont les fibrés cotangents à respectivement
   - La sphère.
   - Le disque.
   - L'hyperboloïde à une nappe.

  (cf. section 4.3)

(V) (cf. Théorème 4.4.1.) Soit $(M, \omega, s)$ un E.S.S. simplement connexe. Soient $o$ un point de $M$ et $K$ l'holonomie linéaire en $o$ relativement à la connexion affine canonique $\nabla$.
Soit $M_c$ une sous-variété de $M$ maximale pour les propriétés suivantes :

  (a) $M_c$ passe par $o$

  (b) $M_c$ est connexe

  (c) $M_c$ est totalement géodésique

  (d) $T_o(M_c)$ est invariant par $K$

  (e) $M_c$ est compacte

  Alors,

  (i) $M_c$ est simplement connexe et symplectique.

  (ii) $M_c$ est la sous-variété de $M$ sous-jacente à un sous-espace symétrique symplectique compact $(M_c, \omega_c, s_c)$ facteur direct de $(M, \omega, s)$.

(iii) $(M_c, \omega_c, s_c)$ admet un groupe des transvections semi-simple compact.

(iv) $M_c$ est l'unique sous-variété de $M$ maximale pour les propriétés (a)–(e)

**Chapitres 5 et 6. Exemples résolubles— Classification en dimension 2 ou 4**

Le but commun de ces deux chapitres est la liste des espaces symétriques symplectiques simplement connexes de dimension 2 ou 4 (cf. section 6.2.).

Le cas des espaces symétriques symplectiques dont le groupe des transvections est semisimple ou ni semisimple ni résoluble a été traité dans les chapitres 3 et 4; les chapitres 5 et 6 exposent donc le cas résoluble.

La raison pour laquelle j'ai scindé la preuve du théorème de classification en deux chapitres est la suivante.

Une grande proportion des structures symétriques locales sous jacentes aux espaces symétriques symplectiques résolubles de dimension 4 apparait dans deux classes, beaucoup plus générales, de structures symétriques résolubles (de toutes dimensions). Les algèbres y apparaissant sont des extensions d'algèbres abéliennes par des algèbres abéliennes.

Il est possible de décrire assez simplement ces classes ainsi que l'isomorphie des éléments qui les constituent; c'est le propos du chapitre 5. Le chapitre 6 traite de la classification proprement dite

# Chapitre 1

# Généralités

**Définition 1.1.** *Un espace symétrique symplectique (brièvement dénoté E.S.S.) est un triple $(M, \omega, s)$ où $M$ est une variété $C^\infty$ connexe, $\omega$ est une forme symplectique sur $M$, $s$ est une application $C^\infty$, $s : M \times M \to M$ telle que si on note $s(x,y) = s_x(y)$ on ait:*

(i) *Pour tout $x$ dans $M$, $s_x$ est un difféomorphisme symplectique involutif de $(M, \omega)$ appelé symétrie en $x$.*

(ii) *Pour tout $x$ dans $M$, $x$ est point fixe isolé de $s_x$.*

(iii) *Pour tous $x$ et $y$ dans $M$, on a $s_x s_y s_x = s_{s_x(y)}$.*

**Définition 1.2.** *Deux E.S.S. $(M, \omega, s)$ et $(M', \omega', s')$ sont isomorphes s'il existe un difféomorphisme symplectique $\varphi : (M, \omega) \to (M', \omega')$ tel que $\varphi s_x = s'_{\varphi(x)} \varphi$. Un tel $\varphi$ est appelé isomorphisme de $(M, \omega, s)$ sur $(M', \omega', s')$. Si $(M, \omega, s) = (M', \omega', s')$ on parle d'automorphisme de l'E.S.S. $(M, \omega, s)$; on note $Aut(M, \omega, s)$ le groupe des automorphismes de $(M, \omega, s)$.*

**Définition 1.3.** *Soit $(M, \omega, s)$ un E.S.S.. Une dérivée covariante affine $\nabla$ sur $M$ est dite admissible si*

(i) $\nabla \omega = 0$.

(ii) *Pour tout $x$ dans $M$, $s_x$ appartient à $Aff(M, \nabla)$, le groupe des transformations affines de l'espace affine $(M, \nabla)$.*

**Théorème 1.4.** *Un E.S.S. admet une unique dérivée covariante admissible appelée dérivée covariante canonique.*

L'existence va résulter des lemmes 1.5 et 1.6 ci-dessous.

**Lemme 1.5.**

(i) *Pour tout $x$ dans $M$ : $s_{x_{\star_x}} = -id|_{T_x(M)}$.*

(ii) *Définissons pour $y$ dans $M$ l'application*

$$\underline{s_y} : M \to M; \underline{s_y}(x) = s_x(y);$$

alors,

$$\underline{s_x}_{\star_x} = 2id|_{T_x(M)}.$$

PREUVE. Soit $\widetilde{g}$ une métrique riemannienne sur $M$. Alors $g = s_x^\star \widetilde{g} + \widetilde{g}$ est une métrique riemannienne sur $M$ admettant $s_x$ comme isométrie. Comme $s_{x_{\star_x}}^2 = id$, $s_{x_{\star_x}}$ n'admet que $\pm 1$ comme valeurs propres. Soit $X \neq 0$ un vecteur propre de valeur propre $+1$ et soit $\gamma(t)$ la géodésique maximale de $(M, g)$ telle que $\gamma(0) = x$ et $\dot{\gamma}(0) = X$. Cette géodésique est fixée point par point par $s_x$ et donc $x$ n'est pas un point fixe isolé de $s_x$. Comme $s_{x_{\star_x}}$ est une transformation orthogonale de $(T_x(M), g_x)$ on a bien $s_{x_{\star_x}} = -id|_{T_x(M)}$. Soit $x(t)$ une courbe $C^\infty$ dans $M$ telle que $x(0) = x$ on a alors;

$$s(x(t), x(t)) = x(t)$$

en dérivant en $t = 0$ on a

$$\underline{s_{x_{\star_x}}}(\dot{x}(0)) + s_{x_{\star_x}}(\dot{x}(0)) = \dot{x}(0)$$

donc $(ii)$. ∎

**Lemme 1.6.** *Soient $X, Y, Z$ trois champs de vecteurs $C^\infty$ sur $M$; alors la formule:*

$$\omega_x(\nabla_X Y, Z) = \frac{1}{2} X_x.\omega(Y + s_{x_\star} Y, Z)$$

*définit une dérivée covariante admissible sur $(M, \omega, s)$.*

PREUVE. Le membre de droite de la formule est $C^\infty(M)$-linéaire en $Z$ (car $Y_x + s_{x_{\star_x}} Y = 0$) et définit donc une 1-forme $C^\infty$ qui peut toujours s'écrire sous la forme $i(\nabla_X Y)\omega$. L'application $(X, Y) \to \nabla_X Y$ est $\mathbb{R}$-bilinéaire et $C^\infty(M)$-linéaire en $X$. Par ailleurs si $f \in C^\infty(M)$, on a

$$s_{x_\star}(fX) = (s_x^\star f)(s_{x_\star} X)$$

et donc:

$$\begin{aligned}
\frac{1}{2} X_x \omega(fY + s_{x_\star}(fY), Z) &= \frac{1}{2} X_x f \omega(Y, Z) \\
&\quad + \frac{1}{2} f(x) X_x \omega(Y, Z) + \frac{1}{2} X_x(s_x^\star f) \omega(s_{x_\star} Y, Z) \\
&\quad + \frac{1}{2} (s_x^\star f)(x) X_x \omega(s_{x_\star} Y, Z) \\
&= X_x f \omega(Y, Z) + \frac{1}{2} f(x) X_x \omega(Y + s_{x_\star} Y, Z)
\end{aligned}$$

ce qui montre que $\nabla_X fY = (Xf)Y + f \nabla_X Y$. Dès lors $\nabla$ est une dérivée covariante affine sur $M$. Soit $y \in M$, $\gamma$ une courbe $C^\infty$ telle que $\gamma(0) = y$ et posons $X_y = \dot{\gamma}(0)$ alors:

$$\begin{aligned}
X_y.\omega(s_{x_\star} Y, Z) &= \left.\frac{d}{dt}\right|_0 \omega_{\gamma(t)}(s_{x_\star} Y, Z) \\
&= \left.\frac{d}{dt}\right|_0 \omega_{s_x(\gamma(t))}(Y, s_{x_\star} Z) \\
&= s_{x_\star}(X_y).\omega(Y, s_{x_\star} Z)
\end{aligned}$$

En particulier

$$X_x.\omega(s_{x_\star} Y, Z) = -X_x.\omega(Y, s_{x_\star} Z),$$

donc

$$\begin{aligned}
\omega_x(\nabla_X Y, Z) &+ \omega_x(Y, \nabla_X Z) = \frac{1}{2} X_x \omega(Y, Z) \\
&+ \frac{1}{2} X_x \omega(s_{x_\star} Y, Z) + \frac{1}{2} X_x \omega(Y, Z) \\
&+ \frac{1}{2} X_x \omega(Y, s_{x_\star} Z) = X_x \omega(Y, Z)
\end{aligned}$$

c'est à dire

$$\nabla_X \omega = 0.$$

Soit $\varphi \in Aut(M, \omega, s)$ alors on a

$$\begin{aligned}
\omega_y(\varphi_\star Y + s_{x_\star} \varphi_\star Y, Z) &= \omega_{\varphi^{-1}(y)}(Y + s_{\varphi^{-1}(x)_\star} Y, \varphi_\star^{-1} Z) \\
&= \left(\varphi^{-1\star}(\omega(Y + s_{\varphi^{-1}(x)_\star} Y, \varphi_\star^{-1} Z))\right)(y)
\end{aligned}$$

c.à.d.

$$\omega(\varphi_\star Y + s_{x_\star} \varphi_\star Y, Z) = \varphi^{-1\star}(\omega(Y + s_{\varphi^{-1}(x)_\star} Y, \varphi_\star^{-1} Z))$$

donc on a

$$\begin{aligned}
\omega_x(\nabla_{\varphi_\star X} \varphi_\star Y, Z) &= \frac{1}{2}(\varphi_\star X)_x . \varphi^{-1\star}(\omega(Y + s_{\varphi^{-1}(x)_\star} Y, \varphi_\star^{-1} Z)) \\
&= \frac{1}{2} X_{\varphi^{-1}(x)} . \omega(Y + s_{\varphi^{-1}(x)_\star} Y, \varphi_\star^{-1} Z) \\
&= \omega_{\varphi^{-1}(x)}(\nabla_X Y, \varphi_\star^{-1} Z) \\
&= \omega_x(\varphi_\star \nabla_X Y, Z)
\end{aligned}$$

donc $Aut(M, \omega, s) \subset Aff(M, \nabla)$

■

PREUVE. (unicité)
Soit $\nabla$ et $\nabla'$ deux dérivées covariantes admissibles sur $(M, \omega, s)$; alors, si $X, Y, Z$ sont trois champs de vecteurs $C^\infty$, le champ $S$ défini par $S(X, Y, Z) = \omega(\nabla_X Y - \nabla'_X Y, Z)$ est un champ de tenseurs totalement symétrique. En effet, $\nabla \omega = 0$ implique

$$S(X, Y, Z) = S(X, Z, Y)$$

De plus, le tenseur de torsion $T$ de $\nabla$ est tel que

$$\begin{aligned}
s_{x_\star} T_x(Y, Z) &= -T_x(Y, Z) \\
&= s_{x_\star}(\nabla_Y Z - \nabla_Z Y - [Y, Z])_x \\
&= \left(\nabla_{s_{x_\star} Y} s_{x_\star} Z - \nabla_{s_{x_\star} Z} s_{x_\star} Y - [s_{x_\star} Y, s_{x_\star} Z]\right)_x \\
&= T_x(Y, Z)
\end{aligned}$$

et est donc identiquement nul. On en déduit

$$S(X, Y, Z) = S(Y, X, Z)$$

et donc la complète symétrie. D'autre part,

$$\begin{array}{rcl}(s_x^\star S)_x(X,Y,Z) & = & S_x(s_{x_\star}X, s_{x_\star}Y, s_{x_\star}Z) \\ & = & -S_x\ (X,Y,Z)\end{array}$$

et

$$\begin{array}{rcl}(s_x^\star S)_x(X,Y,Z) & = & \omega_x(\nabla_{s_{x_\star}X}\ s_{x_\star}Y - \nabla'_{s_{x_\star}X}s_{x_\star}Y, s_{x_\star}Z) \\ & = & \omega_{s_x(x)}(s_{x_\star}(\nabla_X Y - \nabla'_X Y), s_{x_\star}Z) \\ & = & \omega_x(\nabla_X Y - \nabla'_X Y, Z) \\ & = & S_x(X,Y,Z)\end{array}$$

Donc $S \equiv 0$ et $\nabla$ est unique. ∎

Une connexion affine $\nabla$ sur $(M, \omega)$ est dite symplectique si $\nabla \omega = 0$ et si $\nabla$ est sans torsion. Dès lors

**Lemme 1.7.** *Sur un E.S.S. $(M, \omega, s)$ la connexion canonique $\nabla$ est symplectique et on a*

$$Aut(M, \omega, s) \subset Aff(M, \nabla)$$

**Proposition 1.8.**

(i) *La connexion canonique est complète.*

(ii) *Si $Symp(M, \omega)$ désigne le groupe des difféomorphismes symplectiques de $(M, \omega)$, on a :*

$$Aut(M, \omega, s) = Aff(M, \nabla) \cap Symp(M, \omega)$$

*En particulier, $Aut(M, \omega, s)$ a une structure de groupe de Lie de transformations de $M$.*

PREUVE.

(i) Comme $s_x$ est une transformation affine, on a, partout où cela a un sens

$$s_x\ Exp_x\ v = Exp_x\ -v \qquad\qquad v \in T_x(M)$$

Il existe donc un voisinage $U_x$ de $x$, stabilisé par $s_x$ et tel que $s_x|_{U_x}$ coïncide avec la symétrie géodésique au point $x$. Il en résulte immédiatement que toute géodésique maximale de $(M, \nabla)$ est définie sur toute la droite réelle; la connexion est donc complète.

(ii) On a, par le théorème 1, l'inclusion

$$Aut(M, \omega, s) \subset Aff(M, \nabla) \cap Symp(M, \omega)$$

De plus, soit $\psi \in Aff(M, \nabla) \cap Symp(M, \omega)$; alors $\psi s_x \psi^{-1}$ est une affinité stabilisant le point $\psi(x)$ et dont la différentielle au point $\psi(x)$ vaut $-id|_{T_{\psi(x)}(M)}$. Par connexité de $M$, on a donc $\psi s_x \psi^{-1} = s_{\psi(x)}$ et $\psi \in Aut(M, \omega, s)$.
$Aut(M, \omega, s)$ est donc un sous-groupe fermé de $Aff(M, \nabla)$, et par conséquent, un groupe de Lie de transformations de $(M, \omega, s)$. ∎

**Définition 1.9.** *On note $\mathcal{D}er(M)$ l'algèbre de Lie du groupe $Aut(M, \omega, s)$.*

**Définition 1.10.** *Le groupe $G = G(M)$ engendré par $\{s_x \circ s_y;\ x,y \in M\}$ est appelé groupe des transvections de $(M,\omega,s)$.*

Dans un espace affine connexe, deux points peuvent toujours être reliés entre eux par une géodésique brisée; donc, si $x,y \in M$, il existe $z_1,\ldots,z_r \in M$ tels que

$$\begin{aligned} y &= s_{z_r} \circ \ldots \circ s_{z_1}(x) \\ &= s_{z_r} \circ \ldots \circ s_{z_1} \circ s_x(x) \end{aligned}$$

ceci montre que $G$ agit transitivement sur $M$; il en est donc de même pour $Aut(M,\omega,s)$ et pour sa composante connexe au neutre $Aut_0(M,\omega,s)$. En particulier $M$ admet une structure de variété analytique.

**Définition 1.11.** *Fixons $o \in M$; le quadruple $(M,\omega,s,o)$ est alors appelé E.S.S. pointé.*
La conjugaison involutive $\widetilde{\sigma}(g) = s_o g s_o$ est appelée automorphisme canonique de
$Aut(M,\omega,s)$. Sa différentielle au neutre $e$ de $Aut(M,\omega,s)$ est notée $\sigma$; elle induit une décomposition de $\mathcal{D}er(M)$ en sous-espaces propres pour les valeurs propres $+1$ et $-1$; on note respectivement $\mathcal{D}er^+(M)$ et $\mathcal{P}$ ces sous-espaces.

**Proposition 1.12.** *Soit $\widetilde{H}$ le stabilisateur de $o$ dans $Aut(M,\omega,s)$ et soit $\mathcal{H}$ son algèbre de Lie. Alors*

$$\mathcal{H} = \mathcal{D}er^+(M)$$

*Si $H = \widetilde{H} \cap Aut_o(M,\omega,s)$, l'espace homogène*

$$Aut_o(M,\omega,s)/_H$$

*est réductif.*

PREUVE. Si $X \in \mathcal{D}er(M)$ et si $\varphi_t$ est le sous-groupe à un paramètre de $Aut(M,\omega,s)$ associé, on a (en identifiant un élément de $\mathcal{D}er(M)$ au champ de vecteurs associé sur $M$)

$$\sigma(X)_x = \left.\frac{d}{dt} s_o \varphi_t s_o(x)\right|_o = s_{o_{\star_{s_o(x)}}}(X)$$

Si $X \in \mathcal{H}$, on a $\varphi_t s_o \varphi_{-t} = s_{\varphi_t(o)} = s_o$ donc

$$\left.\frac{d}{dt}\varphi_t s_o \varphi_{-t}\right|_o = 0$$

c'est-à-dire

$$X_{s_o(x)} = s_{o_{\star_x}}(X)$$

Remplaçant $x$ par $s_0(x)$, on trouve

$$X_x = \sigma(X)_x$$

c'est-à-dire

$$X \in \mathcal{D}er^+(M)$$

Inversément, si $X \in \mathcal{D}er^+(M)$, on a

$$X_o = s_{o_{\star_{s_o(o)}}}(X) = s_{o_{\star_o}}(X) = -X_o$$

c'est-à-dire $X_o = 0$; donc

$$\begin{aligned} \left.\frac{d}{dt}\varphi_t(o)\right|_{t=s} &= \left.\frac{d}{du}\varphi_{u+s}(o)\right|_0 \\ &= \varphi_{s_{\star_o}}(X_o) = 0 \end{aligned}$$

c'est-à-dire $\varphi_t(o) = \varphi_0(o) = o$ et $X \in \mathcal{H}$.

On a, $\forall d \in Aut(M, \omega, s)$ et $X \in \mathcal{D}er(M)$,

$$\sigma(\mathrm{Ad}(d)X) = \mathrm{Ad}(\widetilde{\sigma}(d))\ \sigma(X)$$

De plus, $\forall h \in \widetilde{H}$,
$$h \cdot o = o = s_o\ h\ s_o(o)$$

et
$$h_{\star_o} = s_{o_{\star_o}}\ h_{\star_o}\ s_{o_{\star_o}}$$

$M$ étant connexe, on a $h = \widetilde{\sigma}(h)$.

Maintenant, $\forall p \in \mathcal{P}$ :
$$\begin{aligned}\sigma(\mathrm{Ad}(h)\ p) &= \mathrm{Ad}(\widetilde{\sigma}(h))\ \sigma(p) \\ &= \mathrm{Ad}(h)\ \sigma(p) \\ &= -\mathrm{Ad}(h)\ p\end{aligned}$$

dès lors $\mathrm{Ad}(\widetilde{H})\ \mathcal{P} \subset \mathcal{P}$ et la réductivité. ∎

Observons que, comme $\sigma$ est un automorphisme de $\mathcal{D}er(M)$, on a

$$\begin{aligned}[\mathcal{H}\ \mathcal{H}] &\subset \mathcal{H} \\ [\mathcal{H}\ \mathcal{P}] &\subset \mathcal{P} \\ [\mathcal{P}\ \mathcal{P}] &\subset \mathcal{H}\end{aligned}$$

Soit $M = D/H$ un espace homogène réductif, c'est-à-dire où $D$ est un groupe de Lie connexe agissant effectivement sur $M$ et où l'algèbre de Lie $\mathcal{D}$ de $D$ se décompose en une somme directe vectorielle $\mathcal{D} = \mathcal{H} \oplus \mathcal{P}$ où $\mathcal{H}$ est l'algèbre de l'isotropie $H$ et $\mathcal{P}$ est un sous-vectoriel $H$-invariant de $\mathcal{D}$. Notons $\pi : D \to M$ la projection canonique. En choisissant $o = \pi(e)$ dans $M$, on peut écrire pour tout $d \in D$, $\pi(d) = d \cdot o$. La différentielle de $\pi$ au neutre de $D$ induit un isomorphisme linéaire $\pi_{\star_e|\mathcal{P}} : \mathcal{P} \to T_o(M)$. L'action de $D$ sur $M$ étant effective, l'homomorphisme $\widetilde{\lambda} : H \to GL(T_o M)$, $\widetilde{\lambda}(h) = h_{\star_o}$, est injectif. Notons $m = \dim M$ et soit $\widetilde{u}_o \in Iso(\mathbb{R}^m, \mathcal{P})$. En $o$, on a : $u_o = \pi_{\star_e|\mathcal{P}} \circ \widetilde{u}_o$ d'où l'isomorphisme de fibrés principaux :

$$D \times_\lambda GL(m, \mathbb{R}) \to B(M)$$
$$[d, A] \to d_{\star_o} \circ u_o \circ A$$

où $\lambda : H \to GL(m, \mathbb{R})$, $\lambda(h) = u_o^{-1} \circ \widetilde{\lambda}(h) \circ u_o$ et où $B(M)$ désigne le fibré des repères linéaires au-dessus de $M$ (par "repère" en un point $x$ de $M$, on entend "élément de $Iso(\mathbb{R}^m, T_x(M))$").

Dès lors, le choix de $u_o$ permet d'identifier $D$ à un $H$-sous-fibré principal de $B(M)$; explicitement :
$i : D \to B(M)$

$$i(d) = d_{\star_o} \circ u_o$$

Le triple $(D, i, \lambda)$ est alors une restriction de $B(M)$ à $H$, c'est-à-dire une $H$-structure sur $M$ (pour plus de précisions voir Dieudonné, tome IV, pages 257-261).

Cette structure est clairement $D$-invariante. Sur le fibré $D$, considérons la 1-forme de connexion canonique $\alpha$

$$\alpha_d = pr_\mathcal{H} \circ L_{d^{-1}_\star}$$

($pr_\mathcal{H} : \mathcal{D} \to \mathcal{H}$ est la projection parallèlement à $\mathcal{P}$).

En notant $\mu$ l'action adjointe $H \times \mathcal{P} \to \mathcal{P}$, le fibré tangent à $M = D/H$ s'identifie au fibré associé $D \times_{(\mu, H)} \mathcal{P}$; explicitement l'isomorphisme est

$$D \times_{(\mu, H)} \mathcal{P} \xrightarrow{j} T(M)$$

$$[d, X] \longrightarrow d_{\star_o}(\pi_{\star_e}(X))$$

Notons $X_x^\star = \frac{d}{dt}\big|_o expt\, X \cdot x$ le champ fondamental sur $M$ associé à $X \in \mathcal{D}$. On a alors

$$X_{\pi(d)}^\star = j\left[d, pr_\mathcal{P}(\mathrm{Ad}(d^{-1})\, X)\right]$$

Le relevé horizontal est donné par

$$\overline{X^\star}\big|_d = \left(\ X\ - \widetilde{\left(pr_\mathcal{H}(\mathrm{Ad}(d^{-1})\, X)\right)}\ \right)_d$$

(où $X$ (resp. $\widetilde{X}$) désigne le champ invariant à droite (resp. à gauche) sur $D$ associé à $X \in \mathcal{D}$); et la fonction $H$-équivariante sur $D$ à valeurs dans $\mathcal{P}$ associée au champ $X_{\pi(d)}^\star$ est donnée par

$$\widehat{X}(d) = pr_\mathcal{P}(\mathrm{Ad}(d^{-1})\, X)$$

Dès lors, en notant $\nabla^{(\alpha)}$ la dérivée covariante sur $T(M)$ associée à $\alpha$, on a

$$\left(\nabla_{X^\star}^{(\alpha)} Y^\star\right)_{\pi(d)} = \left[Y, \mathrm{Ad}(d)\left(\mathrm{Ad}(d^{-1})\, X\right)^\mathcal{P}\right]_{\pi(d)}^\star \tag{1}$$

Comme $d_{1\star_{\pi(d)}} X^\star = (\mathrm{Ad}(d_1)\, X)_{\pi(d_1\, d)}^\star\ \forall d_1, d \in D$, on a $d_{1\star} X^\star = (\mathrm{Ad}(d_1)\, X)^\star$ et dès lors

$$\begin{aligned}\left(\nabla_{d_{1\star} X^\star}^{(\alpha)} d_{1\star} Y^\star\right)_{\pi(d)} &= \left[\mathrm{Ad}(d_1)Y, \mathrm{Ad}(d)\left(\mathrm{Ad}(d^{-1})\, \mathrm{Ad}(d_1)\, X\right)^\mathcal{P}\right]_{\pi(d)}^\star \\ &= \left(\mathrm{Ad}(d_1)\left[Y, \mathrm{Ad}(d_1^{-1}\, d)\left(\mathrm{Ad}(d_1^{-1}\, d)^{-1}\, X\right)^\mathcal{P}\right]\right)_{\pi(d)}^\star \\ &= \left(d_{1\star}\, \nabla_{X^\star}^{(\alpha)} Y^\star\right)_{\pi(d)}\end{aligned}$$

En d'autres termes, $D$ agit par affinités sur $(M, \nabla^{(\alpha)})$. De plus, la formule (1) implique le parallélisme de tout champ de tenseurs $D$-invariant sur $M$.

**Proposition 1.13.** *Soit $(M, \omega, s, o)$ un E.S.S. pointé. Alors, en posant $D = Aut_o(M, \omega, s)$, on a que la connexion canonique $\nabla$ sur l'E.S.S. $(M, \omega, s)$ coïncide avec la connexion canonique $\nabla^{(\alpha)}$ déduite de la structure d'espace homogène réductif $M = D/H$.*

PREUVE. On a

$$(s_{o\star} X^\star)_{\pi(d)} = (\sigma(X))_{\pi(d)}^\star$$
$$\sigma(\mathrm{Ad}(d)\, X) = \mathrm{Ad}(\widetilde{\sigma}d)\, \sigma(X)$$

et

$$pr_\mathcal{P}(\sigma(Y)) = -pr_\mathcal{P}(Y)$$

donc

$$\begin{aligned}\left(\nabla_{s_{o\star} X^\star}^{(\alpha)} s_{o\star} Y^\star\right)_{\pi(d)} &= \left(\sigma\left[Y, \mathrm{Ad}(\widetilde{\sigma}d)\left(\mathrm{Ad}(\widetilde{\sigma}d^{-1})X\right)^\mathcal{P}\right]\right)_{\pi(d)}^\star \\ &= s_{o\star_{s_o\ \pi(d)}}\left[Y, \mathrm{Ad}(\widetilde{\sigma}\, d)\left(\mathrm{Ad}(\widetilde{\sigma}\, d)^{-1}\, X\right)^\mathcal{P}\right]^\star \\ &= s_{o\star_{\widetilde{\pi(\sigma\ d)}}}\left[Y, \mathrm{Ad}(\widetilde{\sigma}\, d)\left(\mathrm{Ad}(\widetilde{\sigma}d)^{-1}\, X\right)^\mathcal{P}\right]^\star \\ &= \left(s_{o\star}\, \nabla_{X^\star}^{(\alpha)} Y^\star\right)_{\pi(d)}\end{aligned}$$

Donc, les symétries sont des affinités de $\nabla^{(\alpha)}$. ∎

**Courbure**. La formule pour la dérivée covariante des champs de vecteurs fondamentaux fournit une expression de l'endomorphisme de courbure

$$(R(X^\star, Y^\star)Z^\star)_{\pi(d)} = -\left[\operatorname{Ad} d\left[\left(\operatorname{Ad} d^{-1}X\right)^{\mathcal{P}}, \left(\operatorname{Ad} d^{-1}Y\right)^{\mathcal{P}}\right], Z\right]^\star_{\pi(d)}$$

En particulier, au point $o = \pi(e)$ de $M$, si $A, B, C \in T_o(M)$ et si $\overline{A}, \overline{B}, \overline{C}$ désignent les éléments correspondants dans $\mathcal{P}$, on a

$$\overline{R(A,B)\,C} = -\,[[\overline{A}\ \overline{B}]\ \overline{C}]$$

On notera donc, assez naturellement, pour trois éléments de $\mathcal{P}$ :

$$R(X,Y)\,Z = -\,[[X\ Y]\ Z]$$

**Définition 1.14.** *[Ko-No, vol 2, th1, pg 222 et 223]*
Une variété affine $(M, \nabla)$ est dite localement symétrique si la torsion et la dérivée covariante de la courbure sont nulles : $T^\nabla = 0$ et $\nabla R^\nabla = 0$.

On sait [Ko-No, vol 2, th 1, pg 222] que $(M, \nabla)$ est localement symétrique si et seulement si, pour tout point $x$ de $M$, il existe un voisinage $U_x$ de ce point stable par la symétrie géodésique $s_x$ au point $x$ et tel que $s_x|_{U_x}$ soit une affinité de $(U_x, \nabla)$.
Une variété affine $(M, \nabla)$ est dite affine symétrique, si la symétrie géodésique locale $s_x$ en $x$ s'étend en une transformation affine de $(M, \nabla)$.

**Proposition 1.15.** *Soit $(M, \omega, s)$ un E.S.S. et soit $\nabla$ sa connexion canonique. Alors $(M, \nabla)$ est une variété affine symétrique.*

Si on note $exp$ l'application exponentielle

$$exp : \mathcal{D}er(M) \to Aut_o(M, \omega, s)$$

on a le

**Corollaire 1.16.** *En identifiant $\mathcal{P}$ à $T_o(M)$, on a, pour tout $X$ dans $\mathcal{P}$ :*

$$exp\,tX = s_{Exp_o(\frac{t}{2}\,X)}\,s_o$$

*En particulier, $s_{Exp_o(\frac{t}{2}\,X)}\,s_o$ induit le transport parallèle le long de la géodésique*
$\pi(\exp\,tX) = Exp_o(tX)$.

PREUVE. Si $Y$ est parallèle le long de $Exp_o(tX)$, $(s_{Exp_o(\frac{t}{2}\,X)}s_o)_\star Y$ l'est aussi, mais $s_{o\star}Y = -Y$ le long de $Exp_o(tX)$ et $s_{Exp_o(\frac{t}{2}\,X)\star}(-Y) = Y$ donc $s_{Exp_o(\frac{t}{2}X)}s_o$ induit le transport parallèle le long de $Exp_o(tX)$. Mais $exp\,tX$ induit aussi le transport parallèle le long de $Exp_o(tX)$ ([Ko-No, vol. 2, page 192]) et

$$\pi(\exp\,tX) = exp\,tX\,o = Exp_o(tX) = s_{Exp_o(\frac{t}{2}\,X)}s_o(o)$$

Donc les deux affinités $exp\,X$ et $s_{Exp_o(\frac{t}{2}\,X)}s_o$ ont le même jet d'ordre 1 en $o$, elles coïncident. ∎

**Proposition 1.17.** *Soit $\mathcal{G} = [\mathcal{P},\,\mathcal{P}] \oplus \mathcal{P} \stackrel{\text{not}}{=} \mathcal{K} \oplus \mathcal{P}$. C'est une sous-algèbre de $\mathcal{D}er(M)$. Le sous-groupe de Lie connexe de $Aut(M, \omega, s)$ d'algèbre $\mathcal{G}$ est le groupe des transvections $G(M)$. C'est le plus petit sous-groupe de $Aut(M, \omega, s)$ transitif sur $M$ et stable par $\widetilde{\sigma}$.*

PREUVE. Le corollaire 1.16 montre que le sous-groupe de Lie connexe $G_1$ de $Aut(M, \omega, s)$ d'algèbre $\mathcal{G}$ est contenu dans $G(M)$. Comme $G(M)$ est engendré par les produits $s_x\,s_y$ et que $s_x\,s_y = s_x\,s_o\,(s_y\,s_o)^{-1}$, pour prouver que $G(M)$ est contenu dans $G_1$, il suffit de prouver que $s_x\,s_o \in G_1, \forall x \in M$.

Le point $x$ pouvant être joint à $o$ par une géodésique brisée, il résulte du corollaire 1.16. qu'il existe $g$ dans $G_1$ tel que $x = g \cdot o$. On a donc, $s_x \, s_o = g \, s_o \, g^{-1} \, s_o = g \, \widetilde{\sigma}(g^{-1})$. $\mathcal{G}$ étant stable par $\sigma$, $G_1$ l'est par $\widetilde{\sigma}$; d'où $G_1 = G(M)$.

Si $G'$ est un sous-groupe de $Aut(M, \omega, s)$ transitif sur $M$ et stable par $\widetilde{\sigma}$, l'argument ci-dessus prouve que $G' \supset G(M)$. ∎

**Remarque 1.18**. Notons $K$ le stabilisateur de $o$ dans $G(M)$; c'est un sous-groupe de Lie d'algèbre $\mathcal{K}$. Comme $M$ est connexe, $K$ est contenu dans le sous-groupe $\widetilde{G^\sigma}$ des points fixés par $\widetilde{\sigma}$ dans $G(M)$. D'autre part, comme $o$ est un point fixe isolé de $s_o$, la composante connexe au neutre $(\widetilde{G^\sigma})_o$ de $\widetilde{G^\sigma}$ est contenue dans $K$.

**Définition 1.19.** *Un sous-E.S.S. d'un E.S.S. $(M, \omega, s)$ est un quadruple $(M_1, \omega_1, s_1, j_1)$ où*

(i) *$(M_1, \omega_1, s_1)$ est un E.S.S.*

(ii) *$j_1 : M_1 \to M$ est une immersion symplectique injective telle que pour tous $x$ et $y$ dans $M_1$ on ait $s_{j_1(x)}(j_1(y)) = j_1(s_x(y))$.*

**Lemme 1.20.** *Soit $(M_1, \omega_1, s_1, j_1)$ un sous-E.S.S. d'un E.S.S. $(M, \omega, s)$. Alors $j_1$ est une affinité relativement aux dérivées covariantes $\nabla$ et $\nabla^1$ canoniquement déduites des structures d'E.S.S. sur respectivement $M$ et $M_1$. En particulier, $(M_1, \nabla^1)$ est une sous variété totalement géodésique complète de $(M, \nabla)$.*

PREUVE. Soit $\gamma : ]-r, r[ \to M_1$ une courbe différentiable dans $M_1$. Notons $\xi_0 = \gamma(0)$ et soient $X$ un champ de vecteurs sur $M_1$ et $\epsilon > 0$ tels que $X_{\gamma(t)} = \dot{\gamma}(t)$ pour tout $t : |t| < \epsilon$. Soit $v$ un champ de vecteurs sur $\gamma$, parallèle le long de $\gamma$ (relativement à la dérivée covariante $\nabla^1$). Soit $Y$ un champ de vecteurs sur $M_1$ tel que $Y_{\gamma(t)} = v_{\gamma(t)}$ pour tout $t : |t| < \epsilon$. Soient $U_0$ un voisinage de $\xi_0$ et $\widetilde{X}$ et $\widetilde{Y}$ deux champs de vecteurs sur $M$ tels que $\widetilde{X}_{j_1(\xi)} = j_{1_{\star_\xi}}(X), \widetilde{Y}_{j_1(\xi)} = j_{1_{\star_\xi}}(Y) \; \forall \xi \in U_0$. Soit $\delta > 0$ avec $\gamma(t) \in U_0 \; \forall t : |t| < \delta$ et soit $\mu :]-\delta, \delta[ \to M$ défini par $\mu = j_1 \circ \gamma$. Soit $\tau \in ]-\delta, \delta[$ tel qu'il existe un voisinage symétrique $U_1$ de $\gamma(\tau)$ avec $\xi_0 \in U_1 \supset U_0$. Soit $Z$ un champ de vecteurs sur $M$. La fonction définie par

$$f_{\tau, Z}(x) = \omega_x(\widetilde{Y} + s_{\mu(\tau)_\star}\widetilde{Y}, Z)$$

est différentiable et on a, pour tout $\xi$ dans $U_1$ :

$$(j_1^\star f_{\tau, Z})(\xi) = \omega_{j_1(\xi)}(j_{1_{\star_\xi}}(Y + s_{1_{\gamma(\tau)_\star}}Y), \pi_\xi(Z_{j_1(\xi)}))$$

où $\pi_\xi$ désigne la projection $\omega_{j_1(\xi)}$−orthogonale de $T_{j_1(\xi)}(M)$ sur $j_{1_{\star_\xi}} T_\xi(M_1)$.

En écrivant localement l'immersion $j_1$ comme la restriction à un ouvert de $\mathbb{R}^{2n}$ ($2n = \dim M_1$) d'une application linéaire injective de $\mathbb{R}^{2n}$ dans $\mathbb{R}^{2m}$ ($2m = \dim M$), on construit, à l'aide d'une section locale du fibré des bases symplectiques au dessus de $M_1$, $2n$ champs de vecteurs locaux, $\{e_1, \ldots, e_n, f_1, \ldots, f_n\}$, au voisinage de $j_1(\xi_0)$ tels que

(i) $\rangle \{e_1, \ldots, e_n, f_1, \ldots, f_n\}_{j_1(\xi)} \langle = j_{1_{\star_\xi}} T_\xi(M_1)$ pour tout $\xi$ dans un voisinage $U_2$ de $\xi_0$.

(ii) $\omega_{j_1(\xi)}(e_i, f_j) = \delta_{ij}, \omega_{j_1(\xi)}(e_i, e_j) = \omega_{j_1(\xi)}(f_i, f_j) = 0 \; \forall \xi \in U_2$.

$\widetilde{Z} = \omega(Z, e_i) f_i + \omega(f_j, Z) e_j$ est un champ de vecteurs local au voisinage de $j_1(\xi_0)$ et, en remarquant que $\widetilde{Z}_{j_1(\xi)} = \pi_\xi(Z) \; \forall \xi \in U_2$, on voit qu'il existe un champ $\overline{Z}$ sur $M_1$ tel qu'au voisinage de $\xi_0$ on ait

$$j_1^\star f_{\tau, Z} = \omega_1(Y + s_{1_{\gamma(\tau)_\star}} Y, \overline{Z})$$

Dès lors, en utilisant la formule du lemme 1.6., on a

$$\begin{array}{rcl} \omega_{\mu(\tau)}(\nabla_{\widetilde{X}} \widetilde{Y}, Z) &=& \frac{1}{2} \widetilde{X}_{\mu(\tau)} \cdot f_{\tau, Z} \\ &=& \frac{1}{2} X_{\gamma(\tau)} \cdot j_1^\star f_{\tau, Z} \\ &=& \omega_{1_{\gamma(\tau)}}(\nabla^1_X Y, \overline{Z}) \\ &=& 0 \end{array}$$

On conclut en utilisant la stabilité de $M_1$ par les symétries centrées en ses points. ∎

**Définition 1.21.** *Un triple symplectique symétrique (T.S.S.) est un triple $t = (\mathcal{G}, \sigma, \Omega)$ où*

(i) *$\mathcal{G}$ est une algèbre de Lie réelle de dimension finie.*

(ii) *$\sigma$ est un automorphisme involutif de $\mathcal{G}$ tel que si $\mathcal{G} = \mathcal{K} \oplus \mathcal{P}$ ($\sigma = id_{\mathcal{K}} \oplus (-id_{\mathcal{P}})$) on ait*

$$[\mathcal{P}\ \mathcal{P}] = \mathcal{K}$$

*et l'action de $\mathcal{K}$ sur $\mathcal{P}$ soit fidèle.*

(iii) *$\Omega$ est une 2-forme symplectique $\mathcal{K}$-invariante sur l'espace vectoriel $\mathcal{P}$.*

*La décomposition $\mathcal{G} = \mathcal{K} \oplus \mathcal{P}$ est appelée décomposition canonique et la dimension de $\mathcal{P}$ est appelée la dimension du triple $t$.*

**Définition 1.22.** *Deux tels triples $t_i = (\mathcal{G}_i, \sigma_i, \Omega_i)$ $i = 1, 2$ seront dits isomorphes si il existe un isomorphisme d'algèbres $\varphi : \mathcal{G}_1 \to \mathcal{G}_2$ tel que*

(i) $\varphi \circ \sigma_1 = \sigma_2 \circ \varphi$

(ii) $\varphi^\star \underline{\Omega}_2 = \underline{\Omega}_1$

*où $\underline{\Omega}_i$ désigne l'extension de $\Omega_i$ à $\mathcal{G}_i \times \mathcal{G}_i$ caractérisée par $i(k_i)\underline{\Omega}_i = 0$ pour tout $k_i$ dans $\mathcal{K}_i$.*

Notons $C^p(\mathcal{G})$ les $p$-cochaines de $\mathcal{G}$ à valeurs dans $\mathbb{R}$ et $\delta$ l'opérateur de cobord de Chevalley associé à la représentation triviale de $\mathcal{G}$ sur $\mathbb{R}$. La $\mathcal{K}$-invariance de $\underline{\Omega}$ implique que $\underline{\Omega}$ est un 2-cocycle de Chevalley.

Etant donné un E.S.S. pointé $(M, \omega, s, o)$, on lui associe un T.S.S. $(\mathcal{G}, \sigma, \Omega)$ où $\mathcal{G}$ est l'algèbre de Lie du groupe des transvections $G$ de $M$; $\sigma$ est la différentielle au neutre $e$ de l'automorphisme $\widetilde{\sigma}$ de $G$ définie par $\widetilde{\sigma}(g) = s_o g s_o$; et où $\underline{\Omega} = (\pi^\star \omega_o)_e$; si $\pi : G \to M$ est la projection canonique.

**Proposition 1.23.** *Cette correspondance induit une bijection entre l'ensemble des classes d'isomorphie d'espaces symétriques symplectiques simplement connexes et l'ensemble des classes d'isomorphie de triples symétriques symplectiques.*

PREUVE

(i) Soient $(M_i, \omega_i, s_i)$ ($i = 1, 2$) deux E.S.S. et soit $\varphi : M_1 \to M_2$ un isomorphisme d'E.S.S.. Associons à $(M_i, \omega_i, s_i)$ un T.S.S. $t_i$ en choisissant un point base $o_i$ dans $M_i$. En modifiant, si nécessaire, $\varphi$ par une transvection de $M_2$, on peut supposer $\varphi(o_1) = o_2$. Cet isomorphisme induit un isomorphisme du groupe des transvections $G(M_1)$ sur le groupe des transvections $G(M_2)$ par

$$\phi(s_x\ s_y) = s_{\varphi(x)}\ s_{\varphi(y)}$$

Le Lemme 1.20. assure que $\phi$ est de Lie. On vérifie que $\phi$ entrelace les automorphismes fondamentaux

$$\phi \circ \widetilde{\sigma}_1 = \widetilde{\sigma}_2 \circ \phi$$

et donc que la différentielle de $\phi$ à l'identité définit un isomorphisme de $t_1$ sur $t_2$.

(ii) Soit $t = (\mathcal{G}, \sigma, \Omega)$ un T.S.S. Soit $\widehat{G}$ le groupe de Lie connexe simplement connexe d'algèbre $\mathcal{G}$, soit $\widehat{\sigma}$ l'automorphisme involutif de $\widehat{G}$ dont la différentielle au neutre est $\sigma$; soit $\underline{\omega}$ la 2-forme invariante à gauche sur $\widehat{G}$ dont la valeur au neutre est $\underline{\Omega}$; soit $\widehat{K}$ le groupe des points fixés par $\widehat{\sigma}$ dans $\widehat{G}$. C'est un sous-groupe fermé et connexe [Koh].

La forme $\underline{\omega}$ est $\widehat{K}$-invariante à droite. Notons $\widehat{\pi} : \widehat{G} \to \widehat{G}/\widehat{K} = M$ la projection canonique; l'espace $M$ est simplement connexe. Si $X$ est un vecteur tangent aux fibres de cette projection, on a $i(X) \underline{\omega} = 0$. La forme $\underline{\omega}$ se projette donc sur la 2-forme $\omega$ sur $M$ (c-à-d $\underline{\omega} = \widehat{\pi}^\star \omega$). Comme $\underline{\Omega}$ est un 2-cocycle, $\underline{\omega}$ est fermée et donc $\omega$ est symplectique. De plus $\omega$ est $\widehat{G}$-invariante.

Définissons la symétrie au point $\widehat{\pi}(g)$ par

$$s_{\widehat{\pi}(g)}(\widehat{\pi}(g')) = \widehat{\pi}(g\ \widehat{\sigma}(g^{-1}g'))$$

On vérifie que ces symétries sont des difféomorphismes symplectiques et ont toutes les propriétés requises pour définir sur $M$ une structure d'E.S.S.

Soit $N = \{g \in \widehat{G} \mid g \cdot x = x\ \forall x \in M\}$ le noyau de l'action $\widehat{G} \times M \to M$.

En remarquant que $N$ coïncide avec l'intersection de tous les stabilisateurs dans $\widehat{G}$ des points de $M$, on voit que $N$ est un sous-groupe fermé, discret et distingué de $\widehat{G}$. $N$ est dès lors central et on a le revêtement universel

$$\widehat{G} \xrightarrow{\rho} G = \widehat{G}/N$$

$G$ agit par automorphismes sur $(M, \omega, s)$, il s'identifie donc à un sous-groupe de Lie connexe de $Aut(M, \omega, s)$. Si $\sum$ désigne l'automorphisme involutif de $Aut(M, \omega, s)$ induit par la symétrie $s_o$ au point $o = \widehat{\pi}(e)$, on a

$$\Sigma \rho(g) = \rho \widehat{\sigma}(g) \qquad\qquad \forall g \in \widehat{G}$$

Dès lors $\Sigma$ stabilise $G$ et comme $\mathcal{G} = [\mathcal{P}\ \mathcal{P}] \oplus \mathcal{P}$, la proposition 1.17 nous dit que $G$ est le groupe des transvections de $(M, \omega, s)$. On voit donc que le T.S.S. associé à $(M, \omega, s)$ est isomorphe à $t$.

(iii) Soient $(\mathcal{G}_i, \sigma_i, \omega_i)$ $(i = 1, 2)$ deux T.S.S. isomorphes et soit $\varphi : \mathcal{G}_1 \to \mathcal{G}_2$ l'isomor- phisme. Il résulte immédiatement de (ii) que $\varphi$ induit un isomorphisme des E.S.S. simplement connexes associés. ∎

**Définition 1.24.** *Soient $(M, \omega)$ une variété symplectique connexe simplement connexe, $G$ un groupe de Lie, $\tau : G \times M \to M$ une action symplectique de $G$ sur $M$. On définit*

(i) *l'isomorphisme $^b : T(M) \to T^\star(M)$, $X^b = i(X)\omega$ et $^\# : T^\star(M) \to T(M)$ son inverse.*

(ii) *$\mathcal{X}_f$ le champ hamiltonien associé à la fonction $f \in C^\infty(M)$ ($-df^\# = \mathcal{X}_f$). On note $\mathcal{H}am(M, \omega)$ l'espace des champs hamiltoniens.*

(iii) *$^\star X$, le champ fondamental associé à un élément $X$ de l'algèbre de Lie $\mathcal{G}$ de $G$ : $^\star X_x = \frac{d}{dt}\big|_o \tau(exp - tX, x)$, $x \in M$.*

(iv) *le crochet de Poisson de deux fonctions $C^\infty$ $f$ et $g$ :*

$$\{f, g\}(x) = \omega_x(\mathcal{X}_f, \mathcal{X}_g)$$

Les faits suivants sont classiques :

(a) $C^\infty(M), \{,\}$ est une algèbre de Lie.

(b) $\mathcal{H}am(M, \omega), [,]$ est une algèbre de Lie et on a l'extension centrale

$$0 \to \mathbb{R} \to C^\infty(M) \xrightarrow{\mathcal{X}} \mathcal{H}am(M, \omega) \to 0$$

(c) L'application $\mathcal{G} \to \Gamma^\infty(T(M)) : X \to {}^\star X$ est un homomorphisme; de plus, si $M$ est simplement connexe, les champs de vecteurs fondamentaux sont hamiltoniens.

**Définition 1.25.** *Avec les notations qui précèdent, on dira que $\tau$ est une action fortement hamiltonienne si*

(i) *les champs de vecteurs fondamentaux sont hamiltoniens,*

(ii) il existe un homomorphisme d'algèbres de Lie $\lambda$ tel que le diagramme suivant commute :

$$0 \longrightarrow \mathbb{R} \longrightarrow C^\infty(M) \longrightarrow \mathcal{H}am(M,\omega) \longrightarrow 0$$

avec $\lambda : \mathcal{G} \to C^\infty(M)$ et $\mathcal{G} \to \mathcal{H}am(M,\omega)$.

**Définition 1.26.** On note pour tout $x$ dans $M$ et $g$ dans $G$, $\tau_x : G \to M$ et $\tau(g) : M \to M$ les applications données par $\tau_x(g') = \tau(g', x)$; $\tau(g)(y) = \tau(g, y)$ ($g' \in G$, $y \in M$) et on définit la fonction $\psi : M \to C^2(\mathcal{G})$ :

$$\psi(x) = (\tau_x^\star \omega)_e$$

La forme $\tau_x^\star \omega$ est fermée et invariante à gauche sur $G$ donc $\psi(x)$ est un 2-cocycle de Chevalley.

**Proposition 1.27.** *Soit $G$ un groupe de Lie connexe agissant transitivement et symplectiquement sur $(M, \omega)$. Si l'action est fortement hamiltonienne, $\psi(x)$ est un 2-cobord pour tout $x$ dans $M$.*
*Réciproquement, si $M$ est simplement connexe et $\psi(o)$ est un 2-cobord pour un $o$ dans $M$, l'action est fortement hamiltonienne.*

PREUVE. La première partie résulte du fait que $M$ est un revêtement symplectique $G$-équivariant d'une orbite coadjointe de $G$ et de la forme particulière de la forme symplectique sur une orbite coadjointe. Pour la seconde partie, notons $H$ le stabilisateur de $o$ dans $G$ et observons que la connexité de $G$ et la simple connexité de $M$ assurent la connexité de $H$. De plus, si $\pi : G \to G/H = M$ est la submersion canonique, on a, avec les notations précédentes, $\psi(o) = (\pi^\star \omega)_e = \delta\xi$ pour un certain $\xi$ dans $\mathcal{G}^\star$. Dès lors, pour tout $X$ dans $\mathcal{H}$ (l'algèbre de Lie de $H$) et $Y$ dans $\mathcal{G}$, on a

$$\xi([X\ Y]) = 0$$

Par connexité de $H$, $\xi$ est $H$-invariant.
L'application $P : M \to \mathcal{G}^\star$ définie par

$$P(\pi(g)) = \text{Ad}^\star(g) \cdot \xi$$

a donc un sens et on vérifie que c'est l'application moment associée à l'action de $G$ sur $M$. ∎

**Définition 1.28.** *Un triple exact (T.E.) est un triple $h = (\mathcal{H}, \sigma, \Omega)$ où*

(i) *$\mathcal{H}$ est une algèbre de Lie réelle de dimension finie*

(ii) *$\sigma$ est un automorphisme involutif de $\mathcal{H}$ tel que si on note $\mathcal{H} = \mathcal{L} \oplus \mathcal{P}$ la décomposition induite par $\sigma$ (i.e. $\sigma = id_{\mathcal{L}} \oplus (-id)_{\mathcal{P}}$), on ait*

$$[\mathcal{P}\ \mathcal{P}] = \mathcal{L}$$

(iii) *$\Omega$ est un 2-forme symplectique sur $\mathcal{P}$, $\mathcal{L}$-invariante et telle qu'il existe $\xi \in \mathcal{G}^\star$ avec $\delta\xi = \underline{\Omega}$ où $\underline{\Omega}$ est l'extension par $o$ sur $\mathcal{L}$ de $\Omega$ à $\mathcal{H}$.*

**Remarques**

(i) $\xi([\mathcal{K}\ \mathcal{P}]) = \xi([\mathcal{K}\ \mathcal{K}]) = 0$

(ii) $\xi$ peut être choisie telle que $\xi(\mathcal{P}) = 0$.

Comme corollaire de la proposition 1.27, on a la

**Proposition 1.29.** *Soit $(M, \omega, s, o)$ un E.S.S. pointé simplement connexe. Soient $G$ son groupe des transvections et $t = (\mathcal{G}, \sigma, \Omega)$ son T.S.S. associé. Alors*

(i) *L'action de $G$ sur $M$ est fortement hamiltonienne si et seulement si $t$ est exact.*

(ii) *Dans ce cas, si $\theta$ désigne l'orbite coadjointe par $\xi$, l'application moment relative à l'action de $G$ : $M \xrightarrow{P} \theta$ est un revêtement analytique symplectique $G$ équivariant.*

**Lemme 1.30.** *Supposons le T.S.S. $t = (\mathcal{G}, \sigma, \Omega)$ non exact. Alors le triple $(\mathcal{H}(\mathcal{G}), \sigma_{\mathcal{H}(\mathcal{G})}, \Omega_{\mathcal{H}(\mathcal{G})})$ défini par*

- $\mathcal{H}(\mathcal{G}) = \mathbb{R} \cdot E \oplus \mathcal{G}$ avec
$$[(t, X), (s, Y)]_{\mathcal{H}(\mathcal{G})} = \underline{\Omega}(X, Y)E \oplus [X, Y]_{\mathcal{G}}$$

- $\sigma_{\mathcal{H}(\mathcal{G})} = id_{\mathbb{R} \cdot E} \oplus \sigma$

- $\Omega_{\mathcal{H}(\mathcal{G})}(0 \oplus p,\ 0 \oplus p') = \Omega(p, p') \quad \forall p, p' \in \mathcal{P}$

*est exact. On note $\mathcal{L}(\mathcal{G}) = \mathbb{R} \cdot E \oplus \mathcal{K}$ et $\mathcal{H}(\mathcal{G}) = \mathcal{L}(\mathcal{G}) \oplus \mathcal{P}$.*

PREUVE. Par $\mathcal{K}$-invariance de $\Omega$ et $\delta$-fermeture de $\underline{\Omega}$, on vérifie que

(a) $\mathcal{H}(\mathcal{G})$ est une algèbre de Lie;

(b) $\sigma_{\mathcal{H}(\mathcal{G})}$ en est un automorphisme involutif;

(c) $\Omega_{\mathcal{H}(\mathcal{G})}$ est une forme symplectique $\mathcal{L}(\mathcal{G})$ invariante sur $\mathcal{P}$;

(d) $\underline{\Omega}_{\mathcal{H}(\mathcal{G})} = -\delta E_\star$ où $E_\star$ est donné par $E_\star(tE \oplus X) = t\ \forall t \in \mathbb{R}, X \in \mathcal{G}$.

Reste à voir $[\mathcal{P}\ \mathcal{P}]_{\mathcal{H}(\mathcal{G})} = \mathcal{L}(\mathcal{G})$.

Posons $\mathcal{H} = [\mathcal{P}\ \mathcal{P}]_{\mathcal{H}(\mathcal{G})} \oplus \mathcal{P}$ et considérons la suite exacte

$$\mathbb{R} \xrightarrow{i} \mathcal{H}(\mathcal{G}) \xrightarrow{\widetilde{\pi}} \mathcal{G}$$

Si $E \notin \mathcal{H}$, $\widetilde{\pi}$ induit un isomorphisme d'algèbres de Lie :

$$\mathcal{H} \xrightarrow{\widetilde{\pi}} \mathcal{G}$$

Comme

$$\Omega_{\mathcal{H}(\mathcal{G})}(p, p') = (\delta E_\star)[p, p'] = \Omega(p, p') \qquad \forall p, p' \in \mathcal{P}$$

$E_\star$ se restreint non-trivialement à $\mathcal{H}$; ce qui contredit l'hypothèse que $t$ n'est pas exact.

On a donc $E \in \mathcal{H}$.

De plus, $[\mathcal{P},\ \mathcal{P}]_{\mathcal{H}} \supset \mathcal{K}$ car si

$$k = \sum_{\alpha, \beta}[p_\alpha, p'_\beta]_{\mathcal{G}} \in \mathcal{K}; \quad p_\alpha,\quad p'_\beta \in \mathcal{P}$$

on a

$$\sum_{\alpha, \beta}[p_\alpha, p'_\beta]_{\mathcal{H}} = k \oplus \left(\sum_{\alpha, \beta} \Omega(p_\alpha, p'_\beta)\right) E$$

donc $k \in \mathcal{H}$ car $E \in \mathcal{H}$. Dès lors

$$[\mathcal{P}\ \mathcal{P}]_{\mathcal{H}(\mathcal{G})} = \mathcal{L}(\mathcal{G})$$

∎

**Remarque 1.31**.
L'application injective
$$\mathcal{K} \longrightarrow \mathcal{H}(\mathcal{G}) : k \longrightarrow 0 \oplus k$$
est un homomorphisme.

De manière analogue au cas exact, on a la

**Proposition 1.32.** *Soit $(M, \omega, s, o)$ un E.S.S. pointé simplement connexe. Soit $G$ son groupe des transvections et $t = (\mathcal{G}, \sigma, \Omega)$ son T.S.S. associé. Soit $H(\mathcal{G})$ le groupe de Lie connexe simplement connexe d'algèbre $\mathcal{H}(\mathcal{G})$. Alors $H(\mathcal{G})$ agit sur $(M, \omega, s)$ par automorphismes et l'application moment relative à cette action est un revêtement symplectique de $M$ sur l'orbite $\theta$ par $E_\star$ dans $\mathcal{H}(\mathcal{G})^\star$.*

<u>Preuve</u>. On a la suite exacte
$$\mathbb{R} \longrightarrow \mathcal{H}(\mathcal{G}) \xrightarrow{\widetilde{\pi}} \mathcal{G}$$

$H(\mathcal{G})$ étant connexe simplement connexe et $G$ connexe, il existe un unique homomorphisme surjectif $C^\infty$

$$m : H(\mathcal{G}) \longrightarrow G$$

tel que $m_{\star_e} = \widetilde{\pi}$.
D'où l'action annoncée.
Cette action est fortement hamiltonienne. ∎

# Chapitre 2

# Décompositions

## 2.1 Couples symétriques

**Définition 2.1.1.** *Un couple symétrique est un couple $c = (\mathcal{G}, \sigma)$ où*

(i) *$\mathcal{G}$ est une algèbre de Lie réelle de dimension finie*

(ii) *$\sigma$ est un automorphisme involutif de $\mathcal{G}$ tel que si $\mathcal{G} = \mathcal{K} \oplus \mathcal{P}$ est la décomposition associée, on a*

  (a) *$\mathcal{K}$ agit fidèlement sur $\mathcal{P}$*

  (b) *$[\mathcal{P}\ \mathcal{P}] = \mathcal{K}$.*

*Un tel couple est dit plat si $\mathcal{K} = \{0\}$.*

**Définition 2.1.2.** *On définit de manière évidente la somme directe de deux couples symétriques. Un couple symétrique est dit décomposable s'il s'écrit comme une somme directe de deux couples symétriques non triviaux. Un couple symétrique est indécomposable si il n'est pas décomposable.*

**Théorème 2.1.3.** *Soit c un couple symétrique. Soient*

$$c = \bigoplus_{\alpha=0}^{p} c_\alpha = \bigoplus_{\beta=0}^{q} c'_\beta$$

*deux décompositions de c telles que*

- *$c_0$ et $c'_0$ sont plats.*
- *Pour tous $\alpha, \beta \geqslant 1$, $c_\alpha$ et $c'_\beta$ sont indécomposables et non plats.*

*Alors*

(i) *$p = q$*

(ii) *Il existe un automorphisme $\varphi$ de c et une permutation $\eta$ de $\{1, \ldots, p = q\}$ tels que*
  *(ii.1) $\varphi\, c_0 = c'_0$. (ii.2) $\varphi\, c_\alpha = c'_{\eta(\alpha)}$ pour tout $\alpha$ dans $\{1, \ldots, p\}$.*

Ceci va résulter d'une succession de lemmes.
Le théorème étant clair si $c$ est plat, on supposera $c = (\mathcal{G}, \sigma)$ non plat. On note

$$\begin{aligned}
\mathcal{K}_\alpha &= \mathcal{K} \cap \mathcal{G}_\alpha \\
\mathcal{K}'_\beta &= \mathcal{K} \cap \mathcal{G}'_\beta \\
\mathcal{P}_\alpha &= \mathcal{P} \cap \mathcal{G}_\alpha \quad \forall \alpha, \beta \geqslant 0 \\
\mathcal{P}'_\beta &= \mathcal{P} \cap \mathcal{G}'_\beta
\end{aligned}$$

$$\overline{c_\alpha} = \bigoplus_{\gamma \neq \alpha} c_\gamma \qquad \overline{c'_\beta} = \bigoplus_{\delta \neq \beta} c'_\delta$$

$$\begin{aligned}
\overline{c_\alpha} &= (\overline{\mathcal{G}_\alpha}, \overline{\sigma_\alpha}) & \overline{c'_\beta} &= (\overline{\mathcal{G}'_\beta}, \overline{\sigma'_\beta}) \\
\overline{\mathcal{K}_\alpha} &= \mathcal{K} \cap \overline{\mathcal{G}_\alpha} \\
\overline{\mathcal{K}'_\beta} &= \mathcal{K} \cap \overline{\mathcal{G}'_\beta} \\
\overline{\mathcal{P}'_\beta} &= \mathcal{P} \cap \overline{\mathcal{G}'_\beta} \\
\overline{\mathcal{P}_\alpha} &= \mathcal{P} \cap \overline{\mathcal{G}_\alpha}
\end{aligned}$$

$\pi_\alpha$ (resp. $\pi'_\beta$) : $\mathcal{G} \to \mathcal{G}_\alpha$ (resp. $\mathcal{G}'_\beta$) la projection parallèlement à $\overline{\mathcal{G}_\alpha}$ (resp. $\overline{\mathcal{G}'_\beta}$),
$\overline{\pi_\alpha}$ (resp. $\overline{\pi'_\beta}$) : $\mathcal{G} \to \overline{\mathcal{G}_\alpha}$ (resp. $\overline{\mathcal{G}'_\beta}$) la projection parallèlement à $\mathcal{G}_\alpha$ (resp. $\mathcal{G}'_\beta$). Remarquons que $\pi_\alpha, \pi'_\beta, \overline{\pi_\alpha}, \overline{\pi'_\beta}$ sont des homomorphismes d'algèbres.

**Lemme 2.1.4.** *Pour tout $\alpha \geqslant 1$, il existe $\beta \geqslant 1$ avec*

$$\mathcal{P}_\alpha \cap \mathcal{P}'_\beta \neq \{0\}$$

PREUVE. Ceci résulte des deux faits suivants :

(a) $[\mathcal{K}_\alpha, \mathcal{P}'_\beta] \subset \mathcal{P}_\alpha \cap \mathcal{P}'_\beta \quad \forall \alpha, \beta \geqslant 0$.

(b) $[\mathcal{K}_\alpha, \mathcal{P}] = \sum_{\beta \geqslant 1} [\mathcal{K}_\alpha, \mathcal{P}'_\beta]$.

∎

Dès maintenant, on fixe $\alpha$ et $\beta$ tels que

$$\alpha, \beta \geqslant 1 \quad \text{et} \quad U \stackrel{def.}{=} \mathcal{P}_\alpha \cap \mathcal{P}'_\beta \neq \{0\}$$

**Lemme 2.1.5.**

(i) $[\pi_\alpha(\overline{\mathcal{P}'_\beta}), \pi_\alpha(\mathcal{P}'_\beta)] = 0 = [\overline{\pi_\alpha}(\mathcal{P}'_\beta), \overline{\pi_\alpha}(\overline{\mathcal{P}'_\beta})]$

(ii) $[\mathcal{K}, \pi_\alpha(\overline{\mathcal{P}'_\beta})] \subset \mathcal{P}_\alpha \cap \overline{\mathcal{P}'_\beta}$

(iii) $[\mathcal{K}, \pi_\alpha(\mathcal{P}'_\beta)] \subset U$

(iv) $U \cap \mathcal{P}_\alpha \cap \overline{\mathcal{P}'_\beta} = 0$

PREUVE.

(i) $[\mathcal{P}'_\beta, \overline{\mathcal{P}'_\beta}] = 0$ donc $\pi_\alpha[\mathcal{P}'_\beta, \overline{\mathcal{P}}'_\beta] = [\pi_\alpha(\mathcal{P}'_\beta), \pi_\alpha(\overline{\mathcal{P}'_\beta})] = 0$; l'autre égalité s'obtient de manière analogue.
(ii) $[\mathcal{K}, \pi_\alpha(\overline{\mathcal{P}'_\beta})] = [\mathcal{K}_\alpha, \pi_\alpha(\overline{\mathcal{P}'_\beta})] = [\mathcal{K}_\alpha, \pi_\alpha(\overline{\mathcal{P}'_\beta}) \oplus \overline{\pi_\alpha}(\overline{\mathcal{P}'_\beta})] = [\mathcal{K}_\alpha, \overline{\mathcal{P}'_\beta}] \subset \overline{\mathcal{P}'_\beta} \cap \mathcal{P}_\alpha$; et de même pour (iii)
(iv) est clair. ∎

**Lemme 2.1.6.** $\pi_\alpha(\mathcal{P}'_\beta) = \mathcal{P}_\alpha$ et $\mathcal{P}_\alpha \cap \overline{\mathcal{P}'_\beta} = 0$.

PREUVE. Par le lemme 2 $(iv)$, il existe $D_1 \subset \pi_\alpha(\overline{\mathcal{P}'_\beta})$ et $C_1 \subset \pi_\alpha(\mathcal{P}'_\beta)$ tels que

$$\begin{aligned} D_1 &\oplus C_1 = \mathcal{P}_\alpha \\ \mathcal{P}_\alpha \cap \overline{\mathcal{P}'_\beta} &\subset D_1 \\ U &\subset C_1 \end{aligned}$$

Mais,
$$[\mathcal{K}\ C_1] \subset [\mathcal{K}\ \pi_\alpha(\mathcal{P}'_\beta)] \subset U \subset C_1$$

et
$$[\mathcal{K}\ D_1] \subset [\mathcal{K}\ \pi_\alpha(\overline{\mathcal{P}'_\beta})] \subset \mathcal{P}_\alpha \cap \overline{\mathcal{P}'_\beta} \subset D_1$$

de plus, $[C_1, D_1] \subset [\pi_\alpha(\mathcal{P}'_\beta), \pi_\alpha(\overline{\mathcal{P}'_\beta})] = 0$.
Donc $\mathcal{G}_\alpha = ([C_1\ C_1] \oplus C_1) \oplus ([D_1,\ D_1] \oplus D_1)$. Par indécomposabilité, on a $C_1$ ou $D_1$ nul; or $C_1 \supset U \neq \{0\}$
donc
$$C_1 = \pi_\alpha(\mathcal{P}'_\beta) = \mathcal{P}_\alpha \qquad \text{et} \qquad \mathcal{P}_\alpha \cap \overline{\mathcal{P}'_\beta} \subset D_1 = \{0\}$$

∎

**Lemme 2.1.7.**

(i) $\dim \mathcal{P}'_\beta \geqslant \dim \mathcal{P}_\alpha$

(ii) Pour tout $\alpha \geqslant 1$, il existe un et un seul $\beta$ tel que $\mathcal{P}'_\beta \cap \mathcal{P}_\alpha \neq \{0\}$

PREUVE.

(i) car $\pi_\alpha(\mathcal{P}'_\beta) = \mathcal{P}_\alpha$.

(ii) car $\mathcal{P}_\alpha \cap \overline{\mathcal{P}'_\beta} = 0$.

∎

Donc, par symétrie, on a le

**Lemme 2.1.8.**

(i) $\dim \mathcal{P}'_\beta = \dim \mathcal{P}_\alpha$

(ii) $p = q$ et la correspondance $\alpha \leftrightarrow \beta \Leftrightarrow \mathcal{P}_\alpha \cap \mathcal{P}'_\beta \neq \{0\}$ induit une bijection de $\{\alpha\}$ sur $\{\beta\}$.

**Lemme 2.1.9.** $\overline{\pi_\alpha}(\mathcal{P}'_\beta)$ est central.

PREUVE.

- $\begin{aligned}[t] [\overline{\pi_\alpha}(\mathcal{P}'_\beta), \overline{\mathcal{P}_\alpha}] &= [\pi_\alpha(\mathcal{P}'_\beta) \oplus \overline{\pi_\alpha}(\mathcal{P}'_\beta), \overline{\mathcal{P}_\alpha}] \\ &= [\mathcal{P}'_\beta, \overline{\mathcal{P}_\alpha}] \\ &\subset \mathcal{K}'_\beta \cap \overline{\mathcal{K}_\alpha} \end{aligned}$

- Soit $k \in \mathcal{K}_\gamma$ ($\gamma \neq \alpha$) avec $[k, \mathcal{P}'_\beta] \neq 0$ alors $\mathcal{P}_\gamma \cap \mathcal{P}'_\beta \supset [\mathcal{K}_\gamma, \mathcal{P}'_\beta] \neq \{0\}$ ce qui contredit le lemme 2.1.8.(ii).

$$\Rightarrow \quad \forall \gamma \neq \alpha \quad [\mathcal{K}_\gamma, \mathcal{P}'_\beta] = 0$$

$$\Rightarrow \quad [\overline{\mathcal{K}_\alpha}, \mathcal{P}'_\beta] = 0$$

$$\Rightarrow \quad [\mathcal{K}'_\beta \cap \overline{\mathcal{K}_\alpha}, \mathcal{P}'_\beta] = 0$$

$$\Rightarrow \quad \mathcal{K}'_\beta \cap \overline{\mathcal{K}_\alpha} \subset Ker\, \mathrm{ad}(\mathcal{K}'_\beta)\Big|_{\mathcal{P}'_\beta} = 0$$

$$\Rightarrow \quad [\overline{\pi_\alpha}(\mathcal{P}'_\beta), \overline{\mathcal{P}_\alpha} \oplus \mathcal{P}_\alpha] = [\overline{\pi_\alpha}(\mathcal{P}'_\beta), \mathcal{P}] = 0$$

■

PREUVE DU THÉORÈME. Soit $\eta$ la permutation de $\{1, \ldots, p = q\}$ définie par

$$\mathcal{P}_\alpha \cap \mathcal{P}'_{\eta(\alpha)} \neq \{0\}$$

Comme, pour tout $\alpha \geqslant 1$ on a $\dim \mathcal{P}_\alpha = \dim \mathcal{P}'_{\eta(\alpha)}$, on a

$$\dim \mathcal{P}_0 = \dim \mathcal{P}'_0$$

Soit $\xi : \mathcal{P}'_0 \to \mathcal{P}_0$ un isomorphisme linéaire et définissons $\phi \in \mathrm{End}(\mathcal{G})$ par

$$\phi = id_\mathcal{K} \oplus \xi \oplus \bigoplus_{\alpha=1}^{p} \pi_\alpha\Big|_{\mathcal{P}'_{\eta(\alpha)}}$$

Comme $[\mathcal{K}, \mathcal{P}_\alpha] \subset \mathcal{P}_\alpha \cap \mathcal{P}'_{\eta(\alpha)}$ on a

$$\phi\Big|_{[\mathcal{G},\,\mathcal{G}]} = id\Big|_{[\mathcal{G},\,\mathcal{G}]}$$

et dès lors par centralité de $\overline{\pi_\alpha}(\mathcal{P}'_{\eta(\alpha)})$, $\phi$ est un automorphisme de $c$. ■

## 2.2 Triples symétriques

**Définition 2.2.1.** *On définit de manière évidente la somme directe de deux T.S.S. Un T.S.S est dit décomposable s'il s'écrit comme une somme directe de deux T.S.S. non triviaux. Un T.S.S. est indécomposable s'il n'est pas décomposable.*

**Théorème 2.2.2.** *Soit $t = (\mathcal{G}, \sigma, \Omega)$ un T.S.S. Soient*

$$t = \bigoplus_{\alpha=0}^{p} t_\alpha = \bigoplus_{\beta=0}^{q} t'_\beta$$

*deux décompositions de t telles que*

- *$t_0$ et $t'_0$ sont plats (Un T.S.S. est plat si son couple symétrique sous-jacent l'est).*
- *$\forall \alpha, \beta \geqslant 1$, $t_\alpha$ et $t'_\beta$ sont indécomposables et non plats.*

*Alors*

(i) *$p = q$*

(ii) *Il existe une permutation $\eta$ de $\{1, \ldots, p = q\}$ et un automorphisme $\varphi$ de t tels que*
  (ii.1) *$\varphi t_0 = t'_0$*
  (ii.2) *$\varphi t_\alpha = t'_{\eta(\alpha)}$ pour tout $\alpha$ dans $\{1, \ldots, p\}$*

La preuve consiste en quelques lemmes et une proposition.

**Lemme 2.2.3.** *Posons $\mathcal{Z}_0 = rad_\Omega(\mathcal{Z}(\mathcal{G}))$.*
*Soient $\mathcal{P}_0$ et $\mathcal{P}'_0$ deux supplémenaires de $\mathcal{Z}_0$ dans $\mathcal{Z}(\mathcal{G})$ :*

$$\mathcal{Z}(\mathcal{G}) = \mathcal{Z}_0 \oplus \mathcal{P}_0 = \mathcal{Z}_0 \oplus \mathcal{P}'_0$$

*Posons*

$$\begin{aligned}
t_0 &= (\mathcal{P}_0, -id, \Omega|_{\mathcal{P}_0 \times \mathcal{P}_0}) \\
t'_0 &= (\mathcal{P}'_0, -id, \Omega|_{\mathcal{P}'_0 \times \mathcal{P}'_0}) \\
\mathcal{P}_1 &= \mathcal{P}_0^\perp \\
\mathcal{P}'_1 &= \mathcal{P}'_0{}^\perp \\
\mathcal{G}_1 &= [\mathcal{P}_1, \mathcal{P}_1] \oplus \mathcal{P}_1 \\
\mathcal{G}'_1 &= [\mathcal{P}'_1, \mathcal{P}'_1] \oplus \mathcal{P}'_1 \\
t_1 &= (\mathcal{G}_1, \sigma|_{\mathcal{G}_1}, \Omega|_{\mathcal{P}_1 \times \mathcal{P}_1}) \\
t'_1 &= (\mathcal{G}'_1, \sigma|_{\mathcal{G}'_1}, \Omega|_{\mathcal{P}'_1 \times \mathcal{P}'_1})
\end{aligned}$$

*Alors*

(i) *$t = t_0 \oplus t_1 = t'_0 \oplus t'_1$ et $\mathcal{Z}(\mathcal{G}_1)$ (resp. $\mathcal{Z}(\mathcal{G}'_1)$) est isotrope dans $\mathcal{P}_1$ (resp. $\mathcal{P}'_1$).*

(ii) *Il existe un automorphisme $\varphi$ de t tel que $\varphi(t_0) = t'_0$ et $\varphi(t_1) = t'_1$.*

PREUVE. Par centralité de $\mathcal{P}_0$ et $\mathcal{P}'_0$, $(i)$ est immédiat.
$(ii)$. Soit $\pi : \mathcal{G} \to \mathcal{P}_0$ la projection parallèlement à $\mathcal{G}_1$; $\pi$ est un homomorphisme d'algèbres.
Soit $\overline{\pi} : \mathcal{G} \to \mathcal{G}_1$ la projection parallèlement à $\mathcal{P}_0$; $\overline{\pi}$ est un homomorphisme d'algèbres.
La restriction de $\pi$ à $\mathcal{P}'_0$ est un isomorphisme de $\mathcal{P}'_0$ sur $\mathcal{P}_0$; en effet, par l'absurde, supposons qu'il existe $z \neq 0$ dans $\mathcal{P}_1 \cap \mathcal{P}'_0$; on a alors $y$ dans $\mathcal{P}'_0$ avec $\Omega(z, y) = 1$ donc $\mathbb{R}z \oplus \mathbb{R}\overline{\pi}(y)$ est un sous-espace symplectique de $\mathcal{P}_1 \cap \mathcal{Z}(\mathcal{G}) = \mathcal{Z}(\mathcal{G}_1)$, ce qui est impossible.
Soit $\overline{Z} = \overline{\pi}(\mathcal{P}'_0)$; comme $\overline{Z} \subset \mathcal{Z}(\mathcal{G}_1)$, $\overline{Z}$ est isotrope.

$\overline{Z}$ et $\overline{Z}^\perp \cap \mathcal{P}_1$ sont contenus dans $\mathcal{P}'_1$.
En effet, on a $\Omega(\overline{Z}, \overline{Z} \oplus \mathcal{P}_0) = 0$ et comme $\mathcal{P}'_0 \subset \overline{Z} \oplus \mathcal{P}_0$ on a $\overline{Z} \subset \mathcal{P}'_1$; de plus $\Omega(\overline{Z}^\perp \cap \mathcal{P}_1, \overline{Z} \oplus \mathcal{P}_0) = 0$ donc $\overline{Z}^\perp \cap \mathcal{P}_1 \subset \mathcal{P}'_1$.
Remarquons que, comme $[\mathcal{K}, \mathcal{P}] = [\mathcal{K}_1, \mathcal{P}'_1]$ et $\Omega([\mathcal{K}, \mathcal{P}], \overline{Z}) = 0$, on a $[\mathcal{K}\,\mathcal{P}] \subset \overline{Z}^\perp \cap \mathcal{P}_1 \cap \mathcal{P}'_1 = \overline{Z}^\perp \cap \mathcal{P}_1$.

Soit $\overline{V}$ un sous-espace isotrope de $\mathcal{P}'_1$ en dualité avec $\overline{Z}$; par un argument dimensionnel on a $\overline{V} \oplus (\overline{Z}^\perp \cap \mathcal{P}'_1) = \mathcal{P}'_1$. Soit $V = \overline{\pi}(\overline{V})$; on vérifie que $V$ est un sous-espace de $\mathcal{P}_1$ en dualité avec $\overline{Z}$. Donc

$$V \oplus (\overline{Z}^\perp \cap \mathcal{P}_1) = \mathcal{P}_1$$

Comme $\overline{Z}^\perp \cap \mathcal{P}_1 \subset \mathcal{P}'_1$ on a par un argument dimensionnel :

$$\overline{Z}^\perp \cap \mathcal{P}_1 = \overline{Z}^\perp \cap \mathcal{P}'_1$$

et

$$\dim V = \dim \overline{V}$$

En particulier, $\overline{\pi}$ établit un isomorphisme linéaire entre $\overline{V}$ et $V$.
Considérons le sous-espace symplectique $V \oplus \overline{Z}$ de $\mathcal{P}_1$. Soit $V \oplus \overline{Z} \xrightarrow{r} V$ la projection parallèlement à $\overline{Z}$.
Soit $\ell : V \to V \oplus \overline{Z}$ une application linéaire telle que $L = \ell(V)$ est lagrangien dans $V \oplus \overline{Z}$ et

$$r \circ \ell = id\big|_V$$

$L$ est alors en dualité avec $\overline{Z}$, on a

$$\mathcal{P}_1 = L \oplus \left(\overline{Z}^\perp \cap \mathcal{P}_1\right)$$

et l'application $\Psi : \overline{V} \to L$ définie par

$$\Psi = \ell \circ \overline{\pi}\big|_{\overline{V}}$$

est un isomorphisme linéaire; on observe de plus que $\forall \overline{v} \in \overline{V}$, $\Psi(\overline{v}) - \overline{v}$ est un élément central.
On vérifie alors que l'application $\varphi : \mathcal{G} \to \mathcal{G}$ définie par

$$\begin{array}{ll} \varphi\big|_{\mathcal{K}} = id\big|_{\mathcal{K}} & \varphi\big|_{\overline{Z}^\perp \cap \mathcal{P}_1} = id\big|_{\overline{Z}^\perp \cap \mathcal{P}_1} \\ \varphi\big|_{\overline{V}} = \Psi & \varphi\big|_{\mathcal{P}'_0} = \pi\big|_{\mathcal{P}'_0} \end{array}$$

livre un automorphisme du triple $t$ comme annoncé. ∎

Nous sommes donc ramenés à prouver le théorème pour un triple $t = (\mathcal{G}, \sigma, \Omega)$ où $\mathcal{Z}(\mathcal{G})$ est totalement isotrope dans $\mathcal{P}$. Soit donc un tel triple et soient

$$t = \bigoplus_{\alpha=1}^{p} t_\alpha = \bigoplus_{\beta=1}^{q} t'_\beta$$

deux décompositions où pour tous $\alpha$ et $\beta$, $t_\alpha$ et $t'_\beta$ sont des triples indécomposables et non plats.
Nous adopterons les notations suivantes, pour tous $\alpha$ et $\beta$ :

- $t_\alpha = (\mathcal{G}_\alpha, \sigma_\alpha, \Omega_\alpha);\ t'_\beta = (\mathcal{G}'_\beta, \sigma'_\beta, \Omega'_\beta)$
- $\mathcal{K}_\alpha = \mathcal{G}_\alpha \cap \mathcal{K};\ \mathcal{K}'_\beta = \mathcal{G}'_\beta \cap \mathcal{K}$
- $\mathcal{P}_\alpha = \mathcal{G}_\alpha \cap \mathcal{P};\ \mathcal{P}'_\beta = \mathcal{G}'_\beta \cap \mathcal{P}$

- $\overline{t_\alpha} = (\overline{\mathcal{G}_\alpha}, \overline{\sigma_\alpha}, \overline{\Omega_\alpha}) = \bigoplus_{\gamma \neq \alpha} t_\gamma$

- $\overline{t'_\beta} = (\overline{\mathcal{G}'_\beta}, \overline{\sigma'_\beta}, \overline{\Omega'_\beta}) = \bigoplus_{\delta \neq \beta} t'_\delta$

- $\overline{\mathcal{K}_\alpha} = \overline{\mathcal{G}_\alpha} \cap \mathcal{K}; \ \overline{\mathcal{K}'_\beta} = \overline{\mathcal{G}'_\beta} \cap \mathcal{K}$

- $\overline{\mathcal{P}_\alpha} = \mathcal{P}_\alpha^\perp = \overline{\mathcal{G}_\alpha} \cap \mathcal{P}; \ \overline{\mathcal{P}'_\beta} = {\mathcal{P}'_\beta}^\perp = \overline{\mathcal{G}'_\beta} \cap \mathcal{P}$

- $\pi_\alpha : \mathcal{G} \to \mathcal{G}_\alpha$ est la projection parallèlement à $\overline{\mathcal{G}_\alpha}$

- $\overline{\pi_\alpha} : \mathcal{G} \to \overline{\mathcal{G}_\alpha}$ est la projection parallèlement à $\mathcal{G}_\alpha$

- $\pi'_\beta : \mathcal{G} \to \mathcal{G}'_\beta$ est la projection parallèlement à $\overline{\mathcal{G}'_\beta}$

- $\overline{\pi'_\beta} : \mathcal{G} \to \overline{\mathcal{G}'_\beta}$ est la projection parallèlement à $\mathcal{G}'_\beta$

De manière identique aux couples symétriques, on a :

**Lemme 2.2.4.** *Pour tout $\alpha$, il existe $\beta$ tel que $\mathcal{P}_\alpha \cap \mathcal{P}'_\beta \neq \{0\}$.*

On fixe alors $\alpha$ et $\beta$ tels que
$$U \stackrel{def.}{=} \mathcal{P}_\alpha \cap \mathcal{P}'_\beta \neq \{0\}$$

Comme pour les couples symétriques, on a :

**Lemme 2.2.5.**

(i) $[\pi_\alpha(\overline{\mathcal{P}'_\beta}), \pi_\alpha(\mathcal{P}'_\beta)] = [\overline{\pi_\alpha}(\mathcal{P}'_\beta) \overline{\pi_\alpha}(\overline{\mathcal{P}'_\beta})] = 0$

(ii) $[\mathcal{K}, \pi_\alpha(\overline{\mathcal{P}'_\beta})] \subset \mathcal{P}_\alpha \cap \overline{\mathcal{P}'_\beta}$

(iii) $[\mathcal{K}, \pi_\alpha(\mathcal{P}'_\beta)] \subset U$

(iv) $[\mathcal{K}, \overline{\pi_\alpha}(\mathcal{P}'_\beta)] \subset \mathcal{P}_\alpha^\perp \cap \mathcal{P}'_\beta$

(v) $[\mathcal{K}, \overline{\pi_\alpha}(\overline{\mathcal{P}'_\beta})] \subset \mathcal{P}_\alpha^\perp \cap {\mathcal{P}'_\beta}^\perp$

(vi) $U \cap \mathcal{P}_\alpha \cap \overline{\mathcal{P}'_\beta} = \{0\}$

**Lemme 2.2.6.** $\overline{\pi_\alpha}(\mathcal{P}'_\beta) \cap \overline{\pi_\alpha}(\overline{\mathcal{P}'_\beta}) \subset \mathcal{Z}(\mathcal{G}) \supset \pi_\alpha(\mathcal{P}'_\beta) \cap \pi_\alpha(\overline{\mathcal{P}'_\beta})$.

PREUVE.
$$\begin{aligned}[\pi_\alpha(\mathcal{P}'_\beta) \cap \pi_\alpha(\overline{\mathcal{P}'_\beta}), \ \pi_\alpha(\mathcal{P}'_\beta) + \pi_\alpha(\overline{\mathcal{P}'_\beta})] = 0 &= [\pi_\alpha(\mathcal{P}'_\beta) \cap \pi_\alpha(\overline{\mathcal{P}'_\beta}), \mathcal{P}_\alpha] \\ &= [\pi_\alpha(\mathcal{P}'_\beta) \cap \pi_\alpha(\overline{\mathcal{P}'_\beta}), \mathcal{P}]\end{aligned}$$

L'autre inclusion s'obtient de la même manière. ■

**Lemme 2.2.7.** $rad_\Omega(\pi_\alpha(\mathcal{P}'_\beta)) \subset rad_\Omega(\mathcal{P}_\alpha \cap \overline{\mathcal{P}'_\beta})$.

PREUVE.
Soit $x \in rad_\Omega(\pi_\alpha(\mathcal{P}'_\beta))$, $y \in \mathcal{P}'_\beta$; alors
$$\Omega(x, y) = \Omega(x, \pi_\alpha(y) + \overline{\pi}_\alpha(y)) = 0$$

et donc $rad_\Omega(\pi_\alpha(\mathcal{P}'_\beta)) \subset \mathcal{P}_\alpha \cap \overline{\mathcal{P}'_\beta} \subset \pi_\alpha(\overline{\mathcal{P}'_\beta})$; en particulier, par le lemme 2.2.6., $rad_\Omega(\pi_\alpha(\mathcal{P}'_\beta)) \subset \mathcal{Z}(\mathcal{G})$.
De plus si $z \in \mathcal{P}_\alpha \cap \overline{\mathcal{P}'_\beta}$ et $t \in \mathcal{P}'_\beta$, $\Omega(z, \pi_\alpha(t)) = \Omega(z, t) = 0$ et donc $rad_\Omega(\pi_\alpha(\mathcal{P}'_\beta)) \subset rad_\Omega(\mathcal{P}_\alpha \cap \overline{\mathcal{P}'_\beta})$ ■

A ce stade ci, deux cas se présentent :
- soit $rad_\Omega(\pi_\alpha(\mathcal{P}'_\beta))$ est nul
- soit il ne l'est pas

**Lemme 2.2.8.** *Supposons $rad_\Omega(\pi_\alpha(\mathcal{P}'_\beta)) = \{0\}$. Alors on a*

(i) $\pi_\alpha(\mathcal{P}'_\beta) = \mathcal{P}_\alpha$

(ii) $\mathcal{P}_\alpha \cap \overline{\mathcal{P}'_\beta} = \{0\}$

(iii) $\pi_\alpha(\overline{\mathcal{P}'_\beta})$ est central

(iv) $\beta$ est l'unique indice tel que $\mathcal{P}_\alpha \cap \mathcal{P}'_\beta \neq \{0\}$

(v) $\dim \mathcal{P}'_\beta \geqslant \dim \mathcal{P}_\alpha$

PREUVE.

(i) $\Omega\left((\pi_\alpha(\mathcal{P}'_\beta))^\perp \cap \mathcal{P}_\alpha, \pi_\alpha(\mathcal{P}'_\beta) \oplus \overline{\pi_\alpha}(\mathcal{P}'_\beta)\right) = 0$ donc, comme $\mathcal{P}'_\beta \subset \pi_\alpha(\mathcal{P}'_\beta) \oplus \overline{\pi_\alpha}(\mathcal{P}'_\beta)$, on a $(\pi_\alpha(\mathcal{P}'_\beta))^\perp \cap \mathcal{P}_\alpha \subset {\mathcal{P}'_\beta}^\perp$ ; dès lors

$$[\pi_\alpha(\mathcal{P}'_\beta),\ (\pi_\alpha(\mathcal{P}'_\beta))^\perp \cap \mathcal{P}_\alpha] = \pi_\alpha[\mathcal{P}'_\beta,\ (\pi_\alpha(\mathcal{P}'_\beta))^\perp \cap \mathcal{P}_\alpha] = 0$$

Par hypothèse et par indécomposabilité de $\mathcal{P}_\alpha$, on a $\pi_\alpha(\mathcal{P}'_\beta) = 0$ ou $(\pi_\alpha(\mathcal{P}'_\beta))^\perp \cap \mathcal{P}_\alpha = 0$. Comme $\pi_\alpha(\mathcal{P}'_\beta) \supset U \neq 0$, on a $\pi_\alpha(\mathcal{P}'_\beta) = \mathcal{P}_\alpha$.

(ii) Comme $\pi_\alpha(\mathcal{P}'_\beta) = \mathcal{P}_\alpha$, $\mathcal{P}_\alpha \cap \overline{\mathcal{P}'_\beta} \cap \pi_\alpha(\mathcal{P}'_\beta) = \mathcal{P}_\alpha \cap \overline{\mathcal{P}'_\beta}$. La preuve du lemme 2.2.7. montre que $\Omega(\pi_\alpha(\mathcal{P}'_\beta),\ \mathcal{P}_\alpha \cap \overline{\mathcal{P}'_\beta})$. Donc

$$\mathcal{P}_\alpha \cap \overline{\mathcal{P}'_\beta} \cap \pi_\alpha(\mathcal{P}'_\beta) \subset rad_\Omega(\pi_\alpha(\mathcal{P}'_\beta)) = \{0\}$$

(iii) $[\pi_\alpha(\overline{\mathcal{P}'_\beta}),\ \mathcal{P}_\alpha] = [\pi_\alpha(\overline{\mathcal{P}'_\beta}), \pi_\alpha(\mathcal{P}'_\beta)] = 0$

(iv) et (v) découlent de (ii) et (i).

∎

**Lemme 2.2.9.** *Supposons $rad_\Omega(\pi_\alpha(\mathcal{P}'_\beta)) = \{0\}$. Si, de plus, $\dim \mathcal{P}_\alpha = \dim \mathcal{P}'_\beta$, alors il existe un automorphisme $\varphi$ de $t$ tel que $\varphi(t'_\beta) = t_\alpha$ et $\varphi(\overline{t'_\beta}) = \overline{t_\alpha}$.*

PREUVE.

(a) $\overline{\pi_\alpha}(\mathcal{P}'_\beta)$ est central.
   En effet,
   $$\begin{aligned}
   [\overline{\pi_\alpha}(\mathcal{P}'_\beta), \mathcal{P}_\alpha^\perp] &= [\pi_\alpha(\mathcal{P}'_\beta) \oplus \overline{\pi_\alpha}(\mathcal{P}'_\beta),\ \mathcal{P}_\alpha^\perp] \\
   &= [\mathcal{P}'_\beta, \mathcal{P}_\alpha^\perp] \\
   &\subset \mathcal{K}'_\beta \cap \overline{\mathcal{K}_\alpha}
   \end{aligned}$$

Le lemme 2.2.8.(ii) affirme que $\mathcal{P}_\alpha \cap {\mathcal{P}'}^\perp_\beta = 0$ ; donc $\mathcal{P}_\alpha^\perp + \mathcal{P}'_\beta = \mathcal{P}$ et par argument dimensionnel la somme est directe, donc $\mathcal{P}_\alpha^\perp \cap \mathcal{P}'_\beta = 0$. Or

$$[\overline{\mathcal{K}_\alpha}, \mathcal{P}'_\beta] \subset \mathcal{P}_\alpha^\perp \cap \mathcal{P}'_\beta = 0$$

Par effectivité, $\overline{\mathcal{K}_\alpha} \cap \mathcal{K}'_\beta = 0$ et donc $[\overline{\pi_\alpha}(\mathcal{P}'_\beta), \mathcal{P}_\alpha^\perp] = [\overline{\pi_\alpha}(\mathcal{P}'_\beta), \mathcal{P}] = 0$.

(b) Comme $\mathcal{P}_\alpha \cap \overline{\mathcal{P}'_\beta} = 0$, $\overline{\pi_\alpha}|_{\overline{\mathcal{P}'_\beta}}$ est injective et par argument dimensionnel un isomorphisme linéaire. Donc $\overline{\pi_\alpha}(\overline{\mathcal{P}'_\beta}) = \mathcal{P}_\alpha^\perp$

(c) On vérifie que $\varphi = id_\mathcal{K} \oplus \pi_\alpha|_{\mathcal{P}'_\beta} \oplus \overline{\pi_\alpha}|_{\overline{\mathcal{P}'_\beta}}$ est comme annoncé.

∎

**Lemme 2.2.10.** *Supposons* $rad_\Omega(\pi_\alpha(\mathcal{P}'_\beta)) = \{0\}$.
*Si pour tous* $\gamma, \delta$ *on a*

$$\mathcal{P}_\gamma \cap \mathcal{P}'_\delta \neq \{0\} \Longrightarrow rad_\Omega\ \pi_\gamma(\mathcal{P}'_\delta) = rad_\Omega\ \pi'_\delta(\mathcal{P}_\gamma) = \{0\}$$

*alors*

$$\dim \mathcal{P}_\alpha = \dim \mathcal{P}'_\beta$$

*En particulier, le théorème 2.2.2. est vrai dans ce cas.*

PREUVE. Ceci est immédiat par le lemme 2.2.8. $(v)$. ∎

Le second cas $(rad_\Omega(\pi_\alpha(\mathcal{P}'_\beta)) \neq 0)$ ne se présente pas :

**Proposition 2.2.11.** *Soit* $t = (\mathcal{G}, \sigma, \Omega)$ *un T.S.S. avec* $\mathcal{Z}(\mathcal{G})$ *isotrope dans* $\mathcal{P}$. *Soient*

$$t = \bigoplus_{\alpha=1}^{p} t_\alpha = \bigoplus_{\beta=1}^{q} t'_\beta$$

*deux décompositions de t en triples indécomposables non plats.*
*Soient* $\alpha$ *et* $\beta$ *tels que* $\mathcal{P}_\alpha \cap \mathcal{P}'_\beta \neq \{0\}$.
*Alors*

$$rad_\Omega(\pi_\alpha(\mathcal{P}'_\beta)) = rad_\Omega(\pi'_\beta(\mathcal{P}_\alpha)) = \{0\}$$

PREUVE. Nous prouvons la proposition 2.2.11. par l'absurde et ceci en plusieurs étapes. Soient donc $\alpha$ et $\beta$ tels que $U = \mathcal{P}_\alpha \cap \mathcal{P}'_\beta \neq \{0\}$ et $rad_\Omega(\pi_\alpha(\mathcal{P}'_\beta)) \neq \{0\}$.
Soit $\overline{c}'_\beta$ le couple symétrique sous-jacent au T.S.S. déduit de $\mathcal{P}'_\beta{}^\perp$.
Soit

$$\overline{c}'_\beta = \bigoplus_{s=0}^{\widetilde{q}} (\overline{c}'_\beta)_s \qquad\qquad (\mathcal{P}'_\beta{}^\perp = \bigoplus_s (\mathcal{P}'_\beta{}^\perp)_s)$$

une décomposition en couples indécomposables non plats et un couple plat $(s = 0)$.
Soit $r \in rad_\Omega(\pi_\alpha(\mathcal{P}'_\beta))$, $r \neq 0$; alors il existe $s \in \{1, \ldots, \widetilde{q}\}$ et $p \in (\mathcal{P}'_\beta{}^\perp)_s$ tels que $\Omega(r, p) = 1$ (Sinon pour tout $s$, on a $\Omega(r, (\mathcal{P}'_\beta{}^\perp)_s) = 0$ donc $\Omega(r, \mathcal{P}'_\beta{}^\perp) = 0$ et $r = 0$ car $r \in \mathcal{P}_\alpha \cap \overline{\mathcal{P}'_\beta}$. En outre si $s = 0$, $p \in \mathcal{Z}(\mathcal{G})$; comme $r$ est aussi central, ceci contredit l'isotropie du centre). Comme $r \in \pi_\alpha(\overline{\mathcal{P}'_\beta})$ (voir la preuve du lemme 2.2.7.), $\pi_\alpha(p) \neq 0$. Dans la suite, nous fixons $r$ et $p$ comme plus haut.

**Première étape** $\overline{\pi}_\alpha(p) \in \mathcal{Z}(\mathcal{G}) \setminus \{0\}$.

En effet, si $\overline{\pi}_\alpha(p) = 0$, $p \in \mathcal{P}_\alpha \cap \overline{\mathcal{P}'_\beta}$ et $\Omega(r, p) = 1$ ce qui contredit $r \in rad_\Omega(\pi_\alpha(\mathcal{P}'_\beta)) \subset rad_\Omega(\mathcal{P}_\alpha \cap \overline{\mathcal{P}'_\beta})$. Comme $\mathcal{Z}(\mathcal{G})$ est isotrope et que $r \in rad_\Omega(\pi_\alpha(\mathcal{P}'_\beta)) \subset \mathcal{Z}(\mathcal{G})$, $\pi_\alpha(p) \notin \mathcal{Z}(\mathcal{G})$.
Soit $c_\gamma$ le couple symétrique sous-jacent au triple $t_\gamma$ ($\gamma \in \{1, \ldots, p\}$).
Soit

$$c_\gamma = \bigoplus_{i_\gamma} (c_\gamma)_{i_\gamma}$$

une décomposition en couples indécomposables non plats et un couple plat ($i_\gamma = 0$). ($\mathcal{P}_\gamma = \bigoplus_{i_\gamma}(\mathcal{P}_\gamma)_{i_\gamma}$). Par le lemme 2.1.7.($ii$), il existe un unique $j_\mu \neq 0$ ($\mu \in \{1, \ldots, p\}$) tel que

$$(\overline{\mathcal{P}}'_\beta)_s \cap (\mathcal{P}_\mu)_{j_\mu} \neq \{0\}$$

et on a, avec les notations du paragraphe précédent

$$p = \pi_{j_\mu}(p) + \overline{\pi}_{j_\mu}(p)$$

où $\overline{\pi}_{j_\mu}(p) \in \mathcal{Z}(\mathcal{G})$ (cf. lemme 2.1.9.).
Si $\mu \neq \alpha$, on a $\pi_\alpha(p) = \pi_\alpha \overline{\pi}_{j_\mu}(p) \in \mathcal{Z}(\mathcal{G})$ ce qui est impossible, donc $\alpha = \mu$
On a alors $\overline{\pi}_\alpha(\overline{\pi}_{j_\mu}(p)) = \overline{\pi}_\alpha(p) \in \mathcal{Z}(\mathcal{G})$.
Ceci achève la preuve de la première étape.
Soit $\xi$ l'élément de $\pi_\alpha(\mathcal{P}'_\beta)^\star$ défini par $\xi(x) = \Omega(\pi_\alpha(p), x)$ pour tout $x$ dans $\pi_\alpha(\mathcal{P}'_\beta)$.
Comme $\xi(r) = -1$, on a $\pi_\alpha(\mathcal{P}'_\beta) = Ker(\xi) \oplus \mathbb{R} \cdot r$.
Posons $\Psi \stackrel{def.}{=} \pi_\alpha^{-1}(Ker(\xi)) \cap \mathcal{P}'_\beta$. On a $\mathcal{P}_\alpha^\perp \cap \mathcal{P}'_\beta \subset \Psi$. Fixons $y \in \mathcal{P}'_\beta$ tel que $\pi_\alpha(y) = r$. Posons $D \stackrel{def.}{=} (\mathbb{R} \cdot r \oplus \mathbb{R}\pi_\alpha(p))^\perp \cap \pi_\alpha(\overline{\mathcal{P}'_\beta})$ et $\Phi \stackrel{def.}{=} \pi_\alpha^{-1}(D) \cap \overline{\mathcal{P}}'_\beta$. On a $\Phi \subset \mathcal{P}_\alpha^\perp \cap \mathcal{P}'_\beta{}^\perp$

**Deuxième étape**

($i$) $\Phi \cap (\mathbb{R} \cdot r \oplus \mathbb{R} \cdot p) = \{0\}$

($ii$) $\Phi \oplus (\mathbb{R} \cdot r \oplus \mathbb{R} \cdot p) = \overline{\mathcal{P}}'_\beta$

($iii$) $\Psi \cap \mathbb{R} \cdot y = \{0\}$

($iv$) $\Psi \oplus \mathbb{R} \cdot y = \mathcal{P}'_\beta$

En effet,
Si ($i$) est faux cela implique qu'il existe $a, b \in \mathbb{R}$ tels que $a\, r + b\, \pi_\alpha(p) \in D$; ceci contredit la définition de $D$. ($ii$) résulte de ($i$) et du fait que $D$ est de codimension au maximum 2 dans $\pi_\alpha(\overline{\mathcal{P}'_\beta})$. ($iii$) résulte de la définition de $\Psi$ et ($iv$) résulte de ($iii$) comme ($ii$) de ($i$). Ceci achève la preuve de la deuxième étape.
On peut donc définir $A$ dans $End(\mathcal{G})$ par

$$A\big|_{\mathcal{K} \oplus \Phi \oplus \Psi} = id\big|_{\mathcal{K} \oplus \Phi \oplus \Psi}$$

$$\begin{array}{rcl} A(p) & = & \pi_\alpha(p) \\ A(y) & = & \overline{\pi}_\alpha(y) \\ A(r) & = & r \end{array}$$

**Troisième étape**

($i$) $A$ est un automorphisme du couple symétrique $c = (\mathcal{G}, \sigma)$.

($ii$) $A(\mathcal{P}'_\beta)$ est un sous-espace symplectique de $(\mathcal{P}, \Omega)$.

($iii$) La restriction $A\big|_{\mathcal{P}'_\beta} : \mathcal{P}'_\beta \to A(\mathcal{P}'_\beta)$ est une application symplectique.

($iv$) $A(\mathcal{P}'_\beta{}^\perp) = (A(\mathcal{P}'_\beta))^\perp$.

En effet, commençons par remarquer que les restrictions de $A$ à $\mathcal{P}'_\beta$ et $\overline{\mathcal{P}}'_\beta$ sont sans noyau.
Supposons que $A(\psi + y) = \psi + \overline{\pi}_\alpha(y) = 0$, pour un certain $\psi \in \Psi$; on a donc $\overline{\pi}_\alpha(y) \in \Psi \subset \mathcal{P}'_\beta$ et donc $r = \pi_\alpha(y) \in \mathcal{P}'_\beta$ ce qui contredit $\Omega(r,p) = 1$.
Supposons que $A(\varphi + a\,r + b\,p) = \varphi + a\,r + b\,\pi_\alpha(p) = 0$, pour un certain $\varphi \in \Psi$. On peut supposer $b \neq 0$; comme $r \in {\mathcal{P}'_\beta}^\perp$, on a $\pi_\alpha(p) \in {\mathcal{P}'_\beta}^\perp$ ce qui contredit $\Omega(r,p) = 1$.

Maintenant, soient $u, v \in \mathcal{P}'_\beta$, on a $u = \psi + a\,y$, $v = \psi' + a'\,y$ où $\psi,\ \psi' \in \Psi$ et $a, a' \in \mathbb{R}$;

$$\begin{aligned} \Omega(u,v) &= \Omega(\psi + a\,A(y) + a\,r, \psi' + a'\,A(y) + a'\,r) \\ &= \Omega(A(u), A(v)) \qquad (r \in rad_\Omega(\pi_\alpha(\mathcal{P}'_\beta)) \subset \overline{\mathcal{P}}'_\beta) \end{aligned}$$

donc $(ii)$ et $(iii)$.
Comme $\Omega(\pi_\alpha(p), Ker(\xi)) = 0$, on a

$$\Omega(A(\mathcal{P}'_\beta), A(\overline{\mathcal{P}}'_\beta)) = 0$$

Donc $A(\overline{\mathcal{P}}'_\beta) \subset A(\mathcal{P}'_\beta)^\perp$; comme $A(\overline{\mathcal{P}}'_\beta)$ a la même dimension que $\overline{\mathcal{P}}'_\beta$ on a l'égalité et comme $A(\mathcal{P}'_\beta)$ est symplectique on a la somme directe orthogonale $A(\mathcal{P}'_\beta) \oplus A(\overline{\mathcal{P}}'_\beta)$; en particulier $A$ est un isomorphisme linéaire.
On a $\Omega(U, \pi_\alpha(p)) = \Omega(U, p) = 0$ donc $U \subset \Psi$. soit $\eta \in (\mathcal{P}_\alpha \cap \overline{\mathcal{P}'_\beta})^\star$ défini par $\Omega(\pi_\alpha(p)\,, .)$; comme $r \in \mathcal{P}_\alpha \cap \overline{\mathcal{P}'_\beta}$ et que $\Omega(r,p) = 1$, on a $\eta \neq 0$ d'où $\mathcal{P}_\alpha \cap \overline{\mathcal{P}'_\beta} = Ker(\eta) \oplus \mathbb{R}.r$. Comme $r \in rad_\Omega(\mathcal{P}_\alpha \cap \overline{\mathcal{P}'_\beta})$, $Ker(\eta) \subset D$ et donc $\mathcal{P}_\alpha \cap \overline{\mathcal{P}'_\beta} \subset \Phi \oplus \mathbb{R}.r$. On a enfin $\overline{\mathcal{P}}_\alpha \cap \overline{\mathcal{P}'_\beta} \subset \Phi$.
Dès lors, on a

$$A\big|_{(\mathcal{P}_\alpha \cap \overline{\mathcal{P}'_\beta}) \oplus U \oplus (\mathcal{P}_\alpha^\perp \cap {\mathcal{P}'_\beta}^\perp) \oplus (\mathcal{P}_\alpha^\perp \cap \mathcal{P}'_\beta)} = id\big|_{(\mathcal{P}_\alpha \cap \overline{\mathcal{P}'_\beta}) \oplus U \oplus (\mathcal{P}_\alpha^\perp \cap {\mathcal{P}'_\beta}^\perp) \oplus (\mathcal{P}_\alpha^\perp \cap \mathcal{P}'_\beta)}$$

Par le lemme 2.2.5., on a

$$[\mathcal{G},\ \mathcal{G}] \subset \mathcal{K} \oplus (\mathcal{P}_\alpha \cap \overline{\mathcal{P}'_\beta}) \oplus U \oplus (\mathcal{P}_\alpha^\perp \cap {\mathcal{P}'_\beta}^\perp) \oplus (\mathcal{P}_\alpha^\perp \cap \mathcal{P}'_\beta)$$

donc $A\big|_{[\mathcal{G},\ \mathcal{G}]} = id\big|_{[\mathcal{G},\ \mathcal{G}]}$ et dès lors, par centralité de $\overline{\pi}_\alpha(p)$ et de $r$, on a $(i)$. Ceci achève la preuve de la troisième étape et nous livre une nouvelle décomposition de $t$, en effet en posant

$$\begin{aligned} \mathcal{G}_1^{(1)} &= \mathcal{K}'_\beta \oplus A(\mathcal{P}'_\beta) \\ \mathcal{P}_1^{(1)} &= A(\mathcal{P}'_\beta) \\ \overline{\mathcal{G}}_1^{(1)} &= \overline{\mathcal{K}}'_\beta \oplus A(\overline{\mathcal{P}}'_\beta) \\ t_1^{(1)} &= (\mathcal{G}_1^{(1)}, \sigma\big|_{\mathcal{G}_1^{(1)}}, \Omega\big|_{A(\mathcal{P}'_\beta) \times A(\mathcal{P}'_\beta)}) \\ \overline{t}_1^{(1)} &= (\overline{\mathcal{G}}_1^{(1)}, \sigma\big|_{\overline{\mathcal{G}}_1^{(1)}}, \Omega\big|_{A(\overline{\mathcal{P}}'_\beta) \times A(\overline{\mathcal{P}}'_\beta)}) \end{aligned}$$

on a

$$t = t_1^{(1)} \oplus \overline{t}_1^{(1)}$$

Le triple $t_1^{(1)}$ est indécomposable et isomorphe au triple $t'_\beta$.

**Quatrième étape**

(i) $\mathcal{P}_1^{(1)} \cap \mathcal{P}_\alpha \supseteq \mathcal{P}'_\beta \cap \mathcal{P}_\alpha = U$

(ii) $\pi_\alpha(\mathcal{P}_1^{(1)}) \subset \pi_\alpha(\mathcal{P}'_\beta)$

(iii) $\dim(rad_\Omega \, \pi_\alpha(\mathcal{P}_1^{(1)})) < \dim(rad_\Omega(\pi_\alpha(\mathcal{P}_\beta')))$

En effet,

(i) On a vu que $A\big|_U = id\big|_U$ donc $U \subset \mathcal{P}_1^{(1)} \cap \mathcal{P}_\alpha$

(ii) $\pi_\alpha(\mathcal{P}_1^{(1)}) = \pi_\alpha(\Psi \oplus \mathbb{R} \cdot \overline{\pi}_\alpha(y)) = Ker(\xi) \subset \pi_\alpha(\mathcal{P}_\beta')$

(iii) On a $rad_\Omega(Ker(\xi)) \oplus \mathbb{R} \cdot r \subseteq rad_\Omega(\pi_\alpha(\mathcal{P}_\beta'))$.

Ceci achève la preuve de la quatrième étape.
Par récurrence, on construit un facteur direct

$$t_1^{(n)} = \left(\mathcal{G}_1^{(n)} = \mathcal{K}_\beta' \oplus \mathcal{P}_1^{(n)}, \; \sigma\big|_{\mathcal{G}_1^{(n)}}, \Omega_{\mathcal{P}_1^{(n)} \times \mathcal{P}_1^{(n)}}\right)$$

de $t$, isomorphe à $t_\beta'$ tel que

$$\mathcal{P}_1^{(n)} \cap \mathcal{P}_\alpha \supseteq U$$

$$\dim rad_\Omega(\pi_\alpha(\mathcal{P}_1^{(n)})) = 0$$

et

$$\pi_\alpha(\mathcal{P}_1^{(n)}) \subset \pi_\alpha(\mathcal{P}_\beta')$$

Mais le lemme 2.2.8.(i) nous dit que $\pi_\alpha(\mathcal{P}_1^{(n)}) = \mathcal{P}_\alpha$ donc $\pi_\alpha(\mathcal{P}_\beta') = \mathcal{P}_\alpha$ et $rad_\Omega(\pi_\alpha(\mathcal{P}_\beta')) = 0$; une contradiction.

∎

**Remarque 2.2.12.** *Soit $t = (\mathcal{G}, \sigma, \Omega)$ un triple symétrique symplectique. Supposons qu'il existe une paire $(V, W)$ de sous espaces de $\mathcal{P}$ telle que*

(i) $\mathcal{P} = V \oplus W$

(ii) $V$ et $W$ sont $\mathcal{K}$-invariants

(iii) $[V, W] = 0$

(iv) $V$ est symplectique.

*Alors, en notant $\pi : \mathcal{P} \to V^\perp$ la projection $\Omega$-orthogonale parallèlement à $V$ et en posant*

$$\widetilde{\Omega} = (\pi\big|_W)^\star(\Omega\big|_{V^\perp \times V^\perp}) \oplus (\Omega\big|_{V \times V})$$

*on a :*

(a) *$\tilde{t} = (\mathcal{G}, \sigma, \tilde{\Omega})$ est un triple symétrique symplectique*

(b) $\tilde{\Omega}(V, W) = 0$

(c) $\tilde{\Omega}\big|_{V \times V} = \Omega\big|_{V \times V}$

Cette remarque motive la définition suivante :

**Définition 2.2.13.** *Soit $t = (\mathcal{G}, \sigma, \Omega)$ un triple symétrique symplectique. On dit que $t$ est presque affinement indécomposable si pour tous sous espaces $V$ et $W$ de $\mathcal{P}$ tels que (i),...,(iv) on a $V$ ou $W$ nul.*

**Théorème 2.2.14.** *Soit $t = (\mathcal{G}, \sigma, \Omega)$ un triple symétrique symplectique . Alors il existe une forme symplectique $\Omega'$ $\mathcal{K}$-invariante sur $\mathcal{P}$ telle que $t' = (\mathcal{G}, \sigma, \Omega')$ admet une décomposition*

$$t' = \bigoplus_{\alpha=0}^{p} t'_\alpha$$

*où*

1. *$t'_0$ est plat*
2. *$t'_\alpha$ est presque affinement indécomposable et non plat pour tout $\alpha \geq 1$.*

*Une telle décomposition de $t'$ est unique à automorphisme de $t'$ près.*

**Remarque 2.2.15.** *Soit $t = (\mathcal{G}, \sigma, \Omega)$ un triple symétrique symplectique . Si $\mathcal{Z}(\mathcal{G}) = \{0\}$ alors*

(i) *la décomposition de $c = (\mathcal{G}, \sigma)$ dans le théorème 2.1.3. est unique à l'ordre des facteurs près.*

(ii) *la décomposition de $t$ dans le théorème 2.2.2. est unique à l'ordre des facteurs près.*

*(ceci est immédiat par les lemmes 2.1.6., 2.1.9., 2.2.8.)*

**Théorème 2.2.16.** *$t = (\mathcal{G}, \sigma, \Omega)$ un triple symétrique symplectique exact. Alors,*

(i) *$\mathcal{Z}(\mathcal{G}) = \{0\}$.*

(ii) *Les décompositions affine (théorème 2.1.3.) et symplectique (théorème 2.2.2.) coïncident.*

PREUVE. Soit $\alpha \in \mathcal{G}^\star$ tel que $\underline{\Omega} = \delta\alpha$. Comme $\Omega$ est non dégénérée on a (i).
(ii) est immédiat par le fait suivant :
si $V$ et $W$ sont tels que

1. $\mathcal{P} = V \oplus W$
2. $[V, W] = 0$

alors $V$ et $W$ sont symplectiques, $\Omega$-orthogonaux et $\mathcal{K}$-invariants.
En effet, comme $\Omega(V, W) = \alpha[V, W] = 0$, seule la $\mathcal{K}$-invariance reste à voir. On a

$$\begin{aligned}
\Omega(V, [\mathcal{K}\,W]) &= \Omega\left(V, [[V\,V] + [W\,W], W]\right) \\
&= \Omega(V, [[W\,W]W]) \quad \text{(Jacobi)} \\
&= \Omega([[W\,W]V], W) \quad (\mathcal{K}\text{-invariance de } \Omega) \\
&= 0 \quad \text{(Jacobi)}
\end{aligned}$$

donc $[\mathcal{K}\,W] \subset V^\perp = W$; idem pour $V$. ∎

Des contre-exemples à la réciproque du théorème 2.2.16. existent en dimension 4 (voir chapitres 5 et 6).

## 2.3 Décompositions des E.S.S.

**Définition 2.3.1.** *On définit de manière évidente le produit de deux E.S.S..Un E.S.S. est dit décomposable s'il est isomorphe au produit de deux E.S.S. non triviaux.Un E.S.S. est indécomposable s'il n'est pas décomposable.*

Les lemme 1.20.,proposition 1.23. et théorème 2.2.2. nous livrent le

**Théorème 2.3.2.** *Tout E.S.S. $(M, \omega, s)$ simplement connexe est le produit d'E.S.S. indécomposables non plats et d'un E.S.S. plat. Une telle décomposition est unique à automorphisme de $(M, \omega, s)$ près.*

**Définition 2.3.3.** *[Wu1 page 307].*
*Deux variétés affines sont CP-difféomorphes s'il existe un difféomorphisme affin entre elles.*

**Définition 2.3.4.** *Soit $(M, \nabla)$ une variété affine symétrique. Une variété affine symétrique $(M_1, \nabla^1)$ est un facteur direct affin symétrique de $(M, \nabla)$ s'il existe une variété affine symétrique $(M_2, \nabla^2)$ telle que $(M_1, \nabla^1) \times (M_2, \nabla^2)$ est CP-difféomorphe à $(M, \nabla)$.*

Le théorème 2.2.16. livre alors le

**Théorème 2.3.5.** *Soit $(M, \omega, s)$ un espace symétrique symplectique simplement connexe tel que l'action du groupe des transvections est fortement hamiltonienne. Alors*

(i) *La décomposition de théorème 2.3.2. est unique à l'ordre des facteurs près et le facteur plat n'y apparait pas.*

(ii) *Tout facteur direct affin symétrique de $(M, \nabla)$ (où $\nabla$ est la connexion canonique sur $(M, \omega, s)$) est un produit de facteurs symplectiques apparaissant dans la décomposition du théorème 2.3.2..*

# Chapitre 3

# Le cas réductif

## 3.1 Espaces symétriques à holonomie complètement réductible – Couples symétriques semi-simples – Corollaires du théorème de Chevalley

Soit $q = (\mathcal{G}, \sigma)$ un couple symétrique.

Nous dirons que $q$ est un couple symétrique réductif (resp. semi-simple, resp. simple, ...) si $\mathcal{G}$ est une algèbre de Lie réductive (resp. semi-simple, resp. simple, ...)

**Théorème 3.1.1.** *[Ko-No, vol II, pg. 327; Ca-Pa pg. 22]*

(i) *Si $A$ est un groupe compact d'automorphismes d'une algèbre de Lie $\mathcal{G}$ de dimension finie, il existe un facteur de Levi $\mathcal{S}$ dans $\mathcal{G}$ stable par $A$*

(ii) *Si $\sigma$ est un automorphisme involutif de $\mathcal{G}$ et si $\mathcal{S}_1$ et $\mathcal{S}_2$ sont deux facteurs de Levi $\sigma$-stables, il existe un automorphisme $\varphi$ du couple $q = (\mathcal{G}, \sigma)$ tel que $\varphi \mathcal{S}_1 = \mathcal{S}_2$.*

Si $q = (\mathcal{G}, \sigma)$ est un C.S. et si $\mathcal{S}$ est un facteur de Levi $\sigma$-stable, on notera $\mathcal{K}_{\mathcal{S}} = \mathcal{K} \cap \mathcal{S}$, $\mathcal{P}_{\mathcal{S}} = \mathcal{P} \cap \mathcal{S}$, $\mathcal{K}_{\mathcal{R}} = \mathcal{K} \cap \mathcal{R}$, $\mathcal{P}_{\mathcal{R}} = \mathcal{P} \cap \mathcal{R}$, de sorte que $\mathcal{K} = \mathcal{K}_{\mathcal{S}} \oplus \mathcal{K}_{\mathcal{R}}$, $\mathcal{P} = \mathcal{P}_{\mathcal{S}} \oplus \mathcal{P}_{\mathcal{R}}$.

**Lemme 3.1.2.** *Soit $q = (\mathcal{G}, \sigma)$ un C.S.. Alors $\mathcal{K}$ agit complètement réductiblement sur $\mathcal{P}$ si et seulement si $\mathcal{G}$ est réductive.*

PREUVE. Rappelons que l'action de $\mathcal{K}_{\mathcal{S}}$ sur $\mathcal{P}$ est complètement réductible [Ca-Pa, pg. 23]; celle de $\mathcal{K}_{\mathcal{R}}$ étant nilpotente ($\mathcal{K}_{\mathcal{R}} \subset [\mathcal{G} \, \mathcal{G}] \cap \mathcal{R}$ le radical nilpotent de $\mathcal{G}$).

Supposons $ad(\mathcal{K})|_{\mathcal{P}}$ complètement réductible et soit $\mathcal{P} = \bigoplus_{i=1}^{r} V_i$ une décomposition en sous-espaces $\mathcal{K}$-irréductibles.

Soit $W_i = \{p \in V_i \mid [\mathcal{K}_{\mathcal{R}}, p] = 0\}$, alors, par Jacobi, $[\mathcal{K}_{\mathcal{R}} [\mathcal{K}_s \ W_i]] = 0$ c.à.d., par définition de $W_i$, $[\mathcal{K}_s \ W_i] \subset W_i$. Dès lors, on a $W_i = V_i$ ou bien $W_i = 0$; mais $W_i \neq 0$ car $ad(\mathcal{K}_{\mathcal{R}})|_{\mathcal{P}}$ est nilpotente dès lors, $[\mathcal{K}_{\mathcal{R}}, \mathcal{P}] = 0$ et $\mathcal{K}_{\mathcal{R}} = 0$ par effectivité.

Ceci livre, par Jacobi, $[\mathcal{P}_{\mathcal{R}}, \mathcal{P}_s] = [\mathcal{K}_s, \mathcal{P}_{\mathcal{R}}] = [\mathcal{P}_{\mathcal{R}} \ \mathcal{P}_{\mathcal{R}}] = 0$ ($[\mathcal{P}_s \ \mathcal{P}_s] = \mathcal{K}_s$ car $[\mathcal{S} \ \mathcal{S}] = \mathcal{S}$) c.à.d. $\mathcal{P}_{\mathcal{R}} = \mathcal{Z}$ est central et $\mathcal{G} = \mathcal{S} \oplus \mathcal{Z}$.

Inversement, si $\mathcal{G} = \mathcal{S} \oplus \mathcal{Z}$ où $\mathcal{S}$ est $\sigma$-stable et $\mathcal{Z} = \mathcal{R}$ est le centre de $\mathcal{G}$, on a $\mathcal{Z} \subset \mathcal{P}$ par effectivité. On conclut par le fait que $\mathcal{K}_{\mathcal{S}} = \mathcal{K}$ est réductive dans $\mathcal{S}$, donc dans $\mathcal{G}$ (voir théorème 3.1.3.). ∎

On dira d'un C.S. $q = (\mathcal{G}, \sigma)$ ou d'un T.S.S. $t = (\mathcal{G}, \sigma, \Omega)$ qu'il est réductif (resp. semi-simple, resp. résoluble, ...) si $\mathcal{G}$ est réductive (resp. semi-simple, resp. résoluble, ...).
Si $\widetilde{q} = (\widetilde{\mathcal{G}}, \widetilde{\sigma})$ est un C.S. réductif de décomposition canonique $\widetilde{\mathcal{G}} = \widetilde{\mathcal{K}} \oplus \widetilde{\mathcal{P}}$, on a la décomposition en somme directe de C.S. :
$$\widetilde{q} = q_0 \oplus q$$
avec
$$q_0 = (\mathcal{Z}, -id) \qquad q = (\mathcal{G} = [\widetilde{\mathcal{G}} \, \widetilde{\mathcal{G}}], \, \sigma = \widetilde{\sigma}|_{\mathcal{G}}$$
où $\mathcal{Z}$ est le centre de $\widetilde{\mathcal{G}}$, $q_0$ est plat et $q$ est semi-simple.

L'étude des C.S. réductifs se ramène donc à celle des semi-simples. Dès maintenant, et tout au long de ce paragraphe, $q = (\mathcal{G}, \sigma)$ désigne un C.S. semi-simple de décomposition canonique $\mathcal{G} = \mathcal{K} \oplus \mathcal{P}$, $\beta$ est la forme de Killing de $\mathcal{G}$, $B$ en est sa restriction à $\mathcal{P} \times \mathcal{P}$ et $\mathcal{Z}(\mathcal{K})$ désigne le centre de $\mathcal{K}$.

**Théorème 3.1.3.** *[Che, pg. 292, prop. 8] $\mathcal{K}$ est réductive dans $\mathcal{G}$.*

**Lemme 3.1.4.** *$\beta(\mathcal{K}, \mathcal{P}) = 0$; en particulier, $\beta|_{\mathcal{K} \times \mathcal{K}}$ et $B$ sont non-singulières.*

PREUVE. En tant qu'automorphisme de $\mathcal{G}$, $\sigma$ conserve $\beta$ et
$$\beta(\mathcal{K}, \mathcal{P}) = \beta(\sigma\mathcal{K}, \sigma\mathcal{P}) = -\beta(\mathcal{K}, \mathcal{P}) = 0$$

■

**Proposition 3.1.5.** *Si $q$ est simple et si $V$ est un sous-espace propre non-trivial $\mathcal{K}$-stable de $\mathcal{P}$, alors*

*(i) $V$ est $B$-isotrope maximal,*

*(ii) $[V \, V] = 0$,*

*(iii) $\mathcal{K}$ agit irréductiblement sur $V$.*

PREUVE. Soit $rad(V)$ le $B$-radical de $V$. Par le théorème 3.1.3., il existe un sous-espace $W$ $\mathcal{K}$-invariant dans $V$ avec
$$rad(V) \oplus W = V$$
($rad(V)$ est $\mathcal{K}$-invariant).
$B|_{W \times W}$ est non-singulière; dès lors,
$$\mathcal{P} = W \oplus W^{\perp_B}$$
Maintenant, en posant $\boldsymbol{a} = [W \, W^{\perp_B}]$, on a :
$$\begin{aligned} B([\boldsymbol{a}, W], W) &= \beta(\boldsymbol{a}, [W \, W]) \\ &= B([[W \, W] \, W^{\perp_B}], \, W) \\ &= 0 \end{aligned}$$
car $W^{\perp_B}$ est $\mathcal{K}$-stable.
De même, $B([\boldsymbol{a} \, W^{\perp_\beta}], \, W^{\perp_\beta}) = 0$ donc, par effectivité, $\boldsymbol{a} = 0$ et $W \oplus [W \, W]$ est un idéal de $\mathcal{G}$. Comme $\mathcal{G}$ est simple,
$$W = 0$$
Soit $\overline{V}$ un sous-espace $\mathcal{K}$-invariant tel que
$$V \oplus \overline{V} = \mathcal{P}$$

Par un raisonnement analogue à celui utilisé plus haut, $\overline{V}$ est $B$-isotrope. Dès lors, $V$ et $\overline{V}$ sont $B$-duaux, isotropes maximaux et donc $\mathcal{K}$ irréductibles.
Enfin, $\forall v_1, v_2, v_3 \in V, \overline{v} \in \overline{V}$, on a

$$\begin{aligned} B([[v_1, v_2]\, \overline{v}],\, v_3) &= B(v_3, [[\overline{v}, v_1]\, v_2] + [[v_2, \overline{v}]\, v_1]) \\ &= 0 \end{aligned}$$

par $\mathcal{K}$-invariance et isotropie de $V$.
Ceci livre, par effectivité : $[V\, V] = 0$. ∎

**Proposition 3.1.6.** *Si $q$ est simple, l'action sur $\mathcal{P}$ de tout élément non nul de $\mathcal{Z}(\mathcal{K})$ est sans noyau.*

PREUVE. Soit $u \in \mathcal{Z}(\mathcal{K})$.
Soit $\mathcal{Z}_u = \{p \in \mathcal{P} \mid [u\, p] = 0\}$. Supposons $\mathcal{Z}_u \neq 0$. Par Jacobi, $\mathcal{Z}_u$ est $\mathcal{K}$-invariant; dès lors, par la proposition 3.1.5., $\mathcal{Z}_u$ est $B$-isotrope maximal et il existe $V$ $\mathcal{K}$-invariant et dual de $\mathcal{Z}_u$ tel que

$$\mathcal{P} = V \oplus \mathcal{Z}_u$$

Pour tous $v, v_1 \in V$ et $p \in \mathcal{Z}_u$, on a :

$$\begin{aligned} B([u\, v], v_1 + p) &= B([u\, v],\, p) \\ &= B(v,\, [u\, p]) \\ &= 0 \end{aligned}$$

donc $[u\, \mathcal{P}] = 0$ et $u = 0$ par effectivité. ∎

**Définition 3.1.7.** *Si $X \in \mathcal{G}$, on note $^\flat X$ l'élément de $\mathcal{G}^\star$ tel que $^\flat X(Y) = \beta(X, Y)\, \forall Y \in \mathcal{G}$ et $\underline{X}$ l'élément de $\Lambda^2(\mathcal{G})$ tel que*

$$\underline{X}(Y, Z) = -^\flat X([Y\, Z]) \qquad\qquad \forall Y, Z \in \mathcal{G}$$

Remarque : $\delta \underline{X} = 0$.

**Théorème 3.1.8.** *Supposons $q$ simple.*
*L'application $\mathcal{Z}(\mathcal{K}) \setminus \{0\} \to \mathcal{Z}^2(\mathcal{G}) : z \mapsto \underline{z}$ définit une bijection entre $\mathcal{Z}(\mathcal{K}) \setminus \{0\}$ et l'ensemble des formes symplectiques $\mathcal{K}$-invariantes sur $\mathcal{P}$.*

PREUVE. Comme $\beta(\mathcal{Z}(\mathcal{K}), [\mathcal{K}\, \mathcal{K}]) = 0$, $\beta|_{\mathcal{Z}(\mathcal{K}) \times \mathcal{Z}(\mathcal{K})}$ est non-singulière; or si $\Omega$ est une forme symplectique $\mathcal{K}$-invariante sur $\mathcal{P}$ on a, comme $H^2(\mathcal{G}) = 0, \underline{\Omega} = \delta \xi$. Comme $H^1(\mathcal{G}) = 0$, un tel élément $\xi \in \mathcal{G}^\star$ est unique et par la remarque page 27 :

$$\xi(\mathcal{P}) = \xi[\mathcal{K}\, \mathcal{K}] = 0$$

c.à.d. $\xi = {}^\flat z$ avec $z \in \mathcal{Z}(\mathcal{K}) \setminus \{0\}$,
c.à.d. $\underline{\Omega} = \underline{z}$.
Soit maintenant $z \in \mathcal{Z}(\mathcal{K}) \setminus \{0\}$ ; alors

$$\Omega = \underline{z}|_{\mathcal{P} \times \mathcal{P}}$$

est symplectique, en effet, si $p \in \mathcal{P}$ est tel que $\underline{z}(p, \mathcal{P}) = 0$, on a $\beta([z\, p],\, \mathcal{P}) = 0$ c.à.d. $p \in \ker(ad(z)|_{\mathcal{P}})$.
Dès lors, $p = 0$ par la proposition 3.1.6.
$\underline{z}$ est visiblement $\mathcal{K}$-invariante; et si $z_1, z_2 \in \mathcal{Z}(\mathcal{K}) \setminus \{0\}$ sont tels que $\underline{z_1} = \underline{z_2}$, on a

$$\beta([z_1 - z_2, \mathcal{P}],\, \mathcal{P}) = 0$$

c.à.d. $z_1 = z_2$ par effectivité; l'application $z \to \underline{z}$ est donc injective. ∎

## 3.2 T.S.S. semi-simples – Indécomposabilité – Rangs

Dans ce paragraphe, $t = (\mathcal{G}, \sigma, \Omega)$ désigne un T.S.S. semi-simple. Le C.S. sous-jacent $q = (\mathcal{G}, \sigma)$ est semi-simple et on adopte les mêmes notations que dans les paragraphes précédents.

Comme $\mathcal{G}$ est semi-simple, $H^1(\mathcal{G}) = H^2(\mathcal{G}) = 0$; dès lors $t$ est exact et il existe un unique $\xi \in \mathcal{G}^\star$ tel que $\delta \xi = \underline{\Omega}$. Dans ce cas, on a vu (cf. théorème 3.1.8) que $\mathcal{Z}(\mathcal{K})$ est un sous-espace non trivial de $\mathcal{K}$ sur lequel la restriction de $\beta$ est non-singulière. En particulier, on a le

**Théorème 3.2.1.** *Un T.S.S. semi-simple est indécomposable si et seulement s'il est simple.*

La proposition 3.1.6. nous livre le

**Lemme 3.2.2.** *Le centralisateur dans $\mathcal{G}$ de $\mathcal{Z}(\mathcal{K})$ est égal à $\mathcal{K}$.*

**Proposition 3.2.3.** *Si $\mathcal{G}$ admet une structure complexe $\mathcal{J}$ (c-à-d $\mathcal{J} \in \mathrm{End}(\mathcal{G}) : \mathcal{J}^2 = -id$ et $\mathcal{J}[X,Y] = [\mathcal{J}X, Y] = [X, \mathcal{J}Y] \; \forall X, Y \in \mathcal{G}$), $\sigma$ commute avec $\mathcal{J}$.*

PREUVE. Par $\beta$-orthogonalité de $\mathcal{P}$ et $\mathcal{K}$, il suffit de prouver que $\mathcal{JK}$ est contenu dans $\mathcal{K}$, ce qui est immédiat par le lemme précédent. ∎

**Théorème 3.2.4.** *Supposons $t$ simple. Alors*

(i) *Si $\mathcal{G}$ admet une structure d'algèbre de Lie complexe (c-à-d $\mathcal{G} \otimes \mathbb{C} = \mathcal{G}^\mathbb{C}$ n'est pas simple), on a :*
$$\dim_\mathbb{C}(\mathcal{Z}(\mathcal{K})) = 1$$

(ii) *Si $\mathcal{G}$ est absolument simple ($\mathcal{G}^\mathbb{C}$ est simple) alors*
$$\dim(\mathcal{Z}(\mathcal{K})) = 1$$

(iii) *$\mathcal{Z}(\mathcal{K})$ est formé d'éléments ad-semi-simples.*

PREUVE.

(i) Par la proposition 3.2.3., $\mathcal{Z}(\mathcal{K})$ est une sous-algèbre complexe. Soit $z \in \mathcal{Z}(\mathcal{K}) \setminus \{0\}$; en tant qu'endomorphisme complexe de $\mathcal{P}$, $z$ possède un vecteur propre. La valeur propre associée est non nulle par la proposition 3.1.6.. Soit alors
$$V_{z,\lambda} = \{p \in \mathcal{P} \mid [z\;p] = \lambda\;p\} \qquad \lambda \in \mathbb{C}_0$$
En utilisant Jacobi, on voit que $V_{z,\lambda}$ est $\mathcal{K}$-invariant. Soit alors $z' \in \mathcal{Z}(\mathcal{K}) \setminus \{0\}$, en tant qu'endomorphisme complexe de $V_{z,\lambda}$, $z'$ possède un vecteur propre $p'$ de valeur propre $\lambda' \in \mathbb{C}_0$. On a $[z - \frac{\lambda}{\lambda'} z', p'] = \lambda\;p' - \lambda p' = 0$ et donc par la proposition 3.1.6., $z = \frac{\lambda}{\lambda'}\;z'$.

(ii) En notant $\sigma^\mathbb{C}$, l'extension $\mathbb{C}$-linéaire de $\sigma$ à $\mathcal{G}^\mathbb{C}$ et $\Omega^\mathbb{C}$, l'extension $\mathbb{C}$-bilinéaire de $\Omega$ à $\mathcal{G}^\mathbb{C}$, on voit que $t^\mathbb{C} = (\mathcal{G}^\mathbb{C}, \sigma^\mathbb{C}, \mathrm{Re}\;\Omega^\mathbb{C})$ définit un T.S.S. simple non-absolument simple; (i) nous dit alors que $\dim_\mathbb{C}(\mathcal{Z}(\mathcal{K}^\mathbb{C})) = 1$. Mais $\mathcal{Z}(\mathcal{K}^\mathbb{C}) = \mathcal{Z}(\mathcal{K})^\mathbb{C}$, donc
$$\dim_\mathbb{C}(\mathcal{Z}(\mathcal{K}^\mathbb{C})) = \dim_\mathbb{R}(\mathcal{Z}(\mathcal{K})) = 1$$

(iii) Par la proposition 3.1.5., $V_{z,\lambda}$ est $B$-isotrope maximal et est un $\mathcal{K}$-module irréductible. Soit $\overline{V}_{z,\lambda}$ un sous-espace $\mathcal{K}$-invariant dans $\mathcal{P}$ tel que $\mathcal{P} = V_{z,\lambda} \oplus \overline{V}_{z,\lambda}$. On a vu qu'alors $\overline{V}_{z,\lambda}$ est en $\beta$-dualité avec $V_{z,\lambda}$; ceci livre $\overline{V}_{z,\lambda} = V_{z,-\lambda}$. ∎

Dès lors, $\mathcal{Z}(\mathcal{K})$ est contenu dans une sous-algèbre de Cartan de $\mathcal{G}$; le Lemme 3.2.2. nous livre alors le

**Théorème 3.2.5.** *Supposons $t$ simple. Toute sous-algèbre de Cartan contenant $\mathcal{Z}(\mathcal{K})$ est contenue dans $\mathcal{K}$; en particulier les rangs de $\mathcal{K}$ et de $\mathcal{G}$ sont égaux.*

## 3.3 T.S.S. simples complexes – Racines admissibles

**Définition 3.3.1.** *Un T.S.S. complexe est un triple $t_c = (\mathcal{G}, \sigma, \Omega_c)$ où*

(i) *$\mathcal{G}$ est une algèbre de Lie complexe de dimension finie et $\sigma$ en est un automorphisme involutif complexe tels que, $\mathcal{G}$ étant vue sur les réels, $(\mathcal{G}, \sigma)$ est un couple symétrique.*

(ii) *En notant $\mathcal{G} = \mathcal{K} \oplus \mathcal{P}$ la décomposition relativement à $\sigma$, $\Omega_c$ est une 2-forme antisymétrique complexe sur $\mathcal{P}$ qui est $\mathbb{C}$-bilinéaire, non singulière et $\mathcal{K}$-invariante.*

Nous avons vu que dans la cas où $t = (\mathcal{G}, \sigma, \Omega)$ est un T.S.S. simple non absolument simple, $\sigma$ préserve la structure complexe sur $\mathcal{G}$. De plus, si on note $\Lambda^2(\mathcal{P})_{\mathbb{C}}$ l'espace des 2-formes antisymétriques complexes sur l'espace $\mathcal{P}$ et $\Lambda^2(\mathcal{P})_{\mathbb{R}}$ l'espace des 2-formes antisymétriques sur l'espace réel $\mathcal{P}$ et si l'on remarque que l'application Re : $\Lambda^2(\mathcal{P})_{\mathbb{C}} \to \Lambda^2(\mathcal{P})_{\mathbb{R}}$, définie par $(\mathrm{Re}\,\Omega)(X,Y) = \mathrm{Re}(\Omega(X,Y))$ pour tous $X$ et $Y$ dans $\mathcal{P}$, est un isomorphisme d'espaces vectoriels réels qui envoie toute forme $\mathcal{K}$-invariante sur une forme $\mathcal{K}$ invariante, on voit que

$$t_c = (\mathcal{G}, \sigma, \Omega_c = \mathrm{Re}^{-1}(\Omega))$$

est un T.S.S. simple complexe.
En outre, il est clair que si $t = (\mathcal{G}_0, \sigma_0, \Omega_0)$ est un T.S.S. simple absolument simple, le triple $t_c = t^{\mathbb{C}} = (\mathcal{G} = \mathcal{G}_0^{\mathbb{C}}, \sigma = \sigma_0^{\mathbb{C}}, \Omega_c = \Omega_0^{\mathbb{C}})$ est un T.S.S. simple complexe.

Ces deux points motivent l'étude des T.S.S. simples complexes. Dans toute la suite de ce paragraphe, on adoptera les notations suivantes :

- $t_c = (\mathcal{G}, \sigma, \Omega_c)$ est un T.S.S. simple complexe,

- $\boldsymbol{h}$ est une sous-algèbre de Cartan de $\mathcal{G}$ contenant $\mathcal{Z}(\mathcal{K})$,

- $\phi$ est le système de racines correspondant à $\boldsymbol{h}$

- $\mathcal{G} = \boldsymbol{h} \oplus \bigoplus_{\alpha \in \phi} \mathcal{G}_\alpha$ est la décomposition radicielle de $\mathcal{G}$ par rapport à $\boldsymbol{h}$.

On a alors immédiatement le

**Lemme 3.3.2.**
$$\mathcal{K} = \boldsymbol{h} \oplus \bigoplus_{\alpha \in \phi_\mathcal{K}} \mathcal{G}_\alpha$$

*et*

$$\mathcal{P} = \bigoplus_{\alpha \in \phi_\mathcal{P}} \mathcal{G}_\alpha$$

*où $\phi_\mathcal{K} = \{\alpha \in \phi \mid \alpha(\mathcal{Z}(\mathcal{K})) = 0\}$, $\phi_\mathcal{P} = \phi \setminus \phi_\mathcal{K}$.*

Remarquons que si on note, comme plus haut, $\underline{\Omega_c} = \delta\xi = \underline{z}$ où $\xi \in \mathcal{G}^\star$ et $z \in \mathcal{Z}(\mathcal{K})$, on a $\mathcal{Z}(\mathcal{K}) = \mathbb{C} \cdot z$ et $\alpha \in \phi_\mathcal{K} \iff \alpha(z) = 0$.

**Définition 3.3.3.** *Une base $\Delta$ de $\phi$ est dite compatible si pour tous $\alpha$ dans $\phi_\mathcal{P}$ et $\beta$ dans $\phi_\mathcal{K}$,*

$$\beta \prec \alpha$$

*où $\prec$ désigne l'ordre sur $\phi$ induit par le choix de $\Delta$.*

**Lemme 3.3.4.** *Il existe une base compatible, $\Delta$, de $\phi$.*

PREUVE. On a $\boldsymbol{h} = (\boldsymbol{h}_{\mathbb{R}})^{\mathbb{C}}$ où

$$\boldsymbol{h}_{\mathbb{R}} = \sum_{\alpha \in \phi} \mathbb{R} \cdot H_\alpha$$

avec $\beta(H_\alpha, H) = \alpha(H)\ \forall H \in \boldsymbol{h}$.

Donc, il existe $z_1$ et $z_2$ dans $\boldsymbol{h}_{\mathbb{R}}$ tels que $z = z_1 + i\, z_2$. On voit alors facilement que $z_1$ et $z_2$ sont dans $\mathcal{Z}(\mathcal{K})$, et sont dès lors $\mathbb{C}$-proportionnels : $z_1 = \mu z_2$. Donc, il existe $\lambda$ dans $\mathbb{C}$ tel que $z = \lambda z_1$.
Le choix d'une base $\beta$-orthonormée de $\boldsymbol{h}_{\mathbb{R}}$ dont le premier élément est $z_1$ induit alors une base compatible $\Delta$ de $\phi$. ∎

Notons $\phi^+$ les racines positives pour ce choix de $\Delta$.
$\Delta$ contient un système de racines simples $\Delta_\mathcal{K}$ de $\mathcal{K}$

$$\Delta_\mathcal{K} = \{\alpha \in \Delta \mid \alpha(z) = 0\}$$

Posons $\Delta_\mathcal{P} = \Delta \setminus \Delta_\mathcal{K} = \Delta \cap \phi_\mathcal{P}$, $\phi_\mathcal{K}^\pm = \pm\phi^+ \cap \phi_\mathcal{K}$ et $\phi_\mathcal{P}^\pm = \pm\phi^+ \cap \phi_\mathcal{P}$.

**Lemme 3.3.5.**

(i) $\Delta_\mathcal{P}$ ne contient qu'un seul élément : $\gamma$.

(ii) $\alpha' \in \phi_\mathcal{P} \iff \alpha' = \displaystyle\sum_{\alpha \in \Delta_\mathcal{K}} n_\alpha \cdot \alpha + n_\gamma \cdot \gamma$ où $n_\gamma \neq 0$.

(iii) Si on note, $\mathcal{P}^\pm = \displaystyle\sum_{\alpha \in \phi_\mathcal{P}^\pm} \mathcal{G}_\alpha$ alors

  1. $\mathcal{P} = \mathcal{P}^+ \oplus \mathcal{P}^-$;
  2. $[\mathcal{P}^+, \mathcal{P}^+] = [\mathcal{P}^-, \mathcal{P}^-] = 0$;
  3. $\mathcal{P}^\pm$ est un $\mathcal{K}$-module irréductible;
  4. $\mathcal{K} \oplus \mathcal{P}^+$ est une sous-algèbre parabolique de $\mathcal{G}$.

(iv) Dans (ii), $n_\gamma = \pm 1$.

PREUVE. $(i)$ et $(ii)$ sont deux conséquences immédiates du fait que

$$\dim \mathcal{Z}(\mathcal{K}) = 1$$

Comme $[\mathcal{P}\ \mathcal{P}] = \mathcal{K}$, on a $(iv)$ et $(iii)2.$.
Pour $(iii)3.$, on remarque que si $\beta \in \phi_\mathcal{K}$, $\alpha \in \phi_\mathcal{P}^+$ avec $\alpha + \beta \in \phi$, alors $[\mathcal{G}_\alpha, \mathcal{G}_\beta] = \mathcal{G}_{\alpha+\beta}$ où $(\alpha+\beta)(z_1) = \beta(z_1) > 0$ donc $\alpha + \beta \in \phi_\mathcal{P}^+$ et $\mathcal{P}^+$ est $\mathcal{K}$-invariant; comme $\mathcal{G}$ est simple, on conclut par la proposition 3.1.5.
Comme $\mathcal{K} \oplus \mathcal{P}^+ = \boldsymbol{h} \oplus \displaystyle\bigoplus_{\alpha \in \phi_\mathcal{K}} \mathcal{G}_\alpha \oplus \bigoplus_{\beta \in \phi_\mathcal{P}^+} \beta \supseteq \boldsymbol{h} \oplus \bigoplus_{\alpha \in \phi^+} \mathcal{G}_\alpha$, on a $(iii)\ 4.$. ∎

Remarquons que $\mathcal{P}^+$ admet la racine maximale $\mu$ relativement à $\Delta$ comme poid dominant.

**Définition 3.3.6.** *Une racine $\alpha \in \phi$ est dite admissible si il existe une base $\Delta'$ de $\phi$ contenant $\alpha$ et telle que si $\mu' \in \phi$ est la racine maximale relativement à $\Delta'$, on ait :*

$$\mu' = \alpha + \sum_{\alpha' \in \Delta' \setminus \{\alpha\}} n_{\alpha'} \cdot \alpha'$$

*où $n_{\alpha'} \in \mathbb{N}$. Une telle base $\Delta'$ sera dite $\alpha$-compatible et le couple $\{\alpha, \Delta'\}$ sera appelé système admissible.*

**Remarques**.

(i) Si $\alpha$ est admissible, toute racine positive $\beta \in \phi^+$ (où $\phi^+$ désigne les racines positives relativement à une base $\alpha$-compatible $\Delta'$) s'écrit

$$\beta = \frac{1}{2}(1+\varepsilon)\alpha + \sum_{\alpha' \in \Delta' \setminus \{\alpha\}} m_{\alpha'} \cdot \alpha'$$

où $m_{\alpha'} \in \mathbb{N}$ et $\varepsilon = \pm 1$.

En effet, soit $C = \{\psi \in \boldsymbol{h}_{\mathbb{R}}^\star \mid (\psi, \alpha) > 0;\ \forall \alpha \in \Delta'\}$ et supposons $\beta \in \phi^+$ avec

$$\beta = \gamma \cdot \alpha + \sum_{\alpha' \in \Delta' \setminus \{\alpha\}} m_{\alpha'} \cdot \alpha', \qquad \gamma \geqslant 2$$

Comme $\mu'$ est maximale, $\mu' - \beta \succ 0$ donc,

$$(\mu' - \beta, \psi) \geqslant 0 \qquad\qquad \psi \in \overline{C}$$

(cf. [Hu], pg. 52). En prenant $\psi = \omega_\alpha^{\Delta'}$ le poid fondamental associé à $\alpha$, on a

$$(\mu' - \beta, \omega_\alpha^{\Delta'}) = 1 - \gamma < 0$$

une contradiction.

(ii) Une base contenant une racine admissible $\alpha$ n'est pas toujours $\alpha$-compatible.

**Théorème 3.3.7.** *Il existe un système admissble $\{\gamma, \Delta\}$ dans $\phi$ et un nombre complexe $\lambda$ tels que*

(i) $\sigma = exp\ \pi\ i\ \mathrm{ad}(h_\gamma^\Delta)$ *où* $\{h_\alpha^\Delta\}$ *désigne la base de* $\boldsymbol{h}_{\mathbb{R}}$ $\beta$*-duale de* $\Delta$;

(ii) $\underline{\Omega_c} = \lambda \underline{h_\gamma^\Delta}$

PREUVE. Les lemmes 3.3.2., 3 et 4 nous fournissent un système admissible : $\{\gamma, \Delta\}$.
Clairement $\mathbb{C} \cdot h_\gamma^\Delta = \mathcal{Z}(\mathcal{K})$ d'où (ii).
Pour (i), on remarque, en décomposant

$$\mathcal{K} = \boldsymbol{h} \oplus \bigoplus_{\alpha \in \phi_\mathcal{K}} \mathcal{G}_\alpha \qquad \text{et} \qquad \mathcal{P} = \bigoplus_{\alpha \in \phi_\mathcal{P}} \mathcal{G}_\alpha$$

que

$$exp(\pi i \mathrm{ad}(h_\gamma^\Delta)) = id_\mathcal{K} \oplus (-id_\mathcal{P})$$

∎

Réciproquement,

**Théorème 3.3.8.** *Soient*

- $\mathcal{G}$ *une algèbre de Lie simple complexe,*
- $\boldsymbol{h}$ *une sous-algèbre de Cartan de $\mathcal{G}$,*
- $\phi$ *le système de racines associé,*
- $\Delta = \{\alpha_1, \ldots, \alpha_\ell\}$ *une base de $\phi$ et*
- $\{h_{\alpha_i}^\Delta\}$ *la base de $\boldsymbol{h}_{\mathbb{R}}$, $\beta$-duale de $\Delta$.*

*alors,*

(i) *Pour tout $i$, $\sigma_i = exp\ \pi\ i\ \mathrm{ad}(h_{\alpha_i}^\Delta)$ est un automorphisme involutif de $\mathcal{G}$.*

(ii) Si $\mathcal{G} = \mathcal{K}_i \oplus \mathcal{P}_i$ désigne la décomposition relativement à $\sigma_i$, on a $\mathcal{Z}(\mathcal{K}_i) \neq \{0\}$ si et seulement si $\{\alpha_i, \Delta\}$ est admissible.

(iii) Dans ce cas, $\dim \mathcal{Z}(\mathcal{K}_i) = 1$ et pour tout $\lambda \in \mathbb{C}^\star$,

$$t_{i,\lambda} = \left( \mathcal{G}, \sigma_i, \underline{\lambda(h_{\alpha_i}^\Delta)}\Big|_{\mathcal{P}_i \times \mathcal{P}_i} \right)$$

est un T.S.S. simple complexe.

PREUVE.

(i) est immédiat

(ii) (a) Si $\mathcal{Z}(\mathcal{K}_i) \neq \{0\}$, le théorème 3.3.7. nous livre un système admissible $\{\gamma, \Delta'\}$ avec

$$\sigma_i = exp\ \pi\ i \text{ad}(h_\gamma^{\Delta'})$$

et donc

$$h_{\alpha_i}^\Delta = \pm h_\gamma^{\Delta'}$$

On peut supposer $h_{\alpha_i}^\Delta = h_\gamma^{\Delta'}$ (car si $\{\gamma, \Delta'\}$ est admissible, $\{-\gamma, -\Delta'\}$ l'est aussi et $h_\gamma^{\Delta'} = -h_{-\gamma}^{-\Delta'}$).

Dès lors, si $\mu$ est la racine maximale pour $\Delta$, on a

$$\mu(h_{\alpha_i}^\Delta) = \mu(h_\gamma^{\Delta'}) = 1$$

(cf. lemme 3.3.5.).

(b) Supposons $\{\alpha_i, \Delta\}$ admissible.

Posons $\phi_0 = \{\alpha \in \phi \mid \alpha(h_{\alpha_i}^\Delta) = 0\}$ et $\phi_p = \phi \setminus \phi_0$. Soient alors $\mathcal{K} = \boldsymbol{h} \oplus \bigoplus_{\alpha \in \phi_0} \mathcal{G}_\alpha$ et $\mathcal{P} = \bigoplus_{\alpha \in \phi_p} \mathcal{G}_\alpha$. On vérifie alors que $\sigma_i = id_\mathcal{K} \oplus (-id_\mathcal{P})$ donc $\mathcal{K} = \mathcal{K}_i$ et $\mathcal{P} = \mathcal{P}_i$.

On remarque alors que $[h_{\alpha_i}^\Delta, \mathcal{K}] = 0$.

(iii) est immédiat.

∎

## 3.4 Structure des T.S.S. simples absolument simples

Le paragraphe précédent livre immédiatement le

**Théorème 3.4.1.** Soient $t = (\mathcal{G}_0, \sigma_0, \Omega)$ un T.S.S. simple absolument simple, $t_c = (\mathcal{G}, \sigma, \Omega_c)$ son T.S.S. "complexifié", $\boldsymbol{h}_0$ une sous-algèbre de Cartan de $\mathcal{G}_0$ contenant $\mathcal{Z}(\mathcal{K}_0)$ et $\boldsymbol{h} = \boldsymbol{h}_0^\mathbb{C}$. Si $\phi$ désigne le système de racines de $\mathcal{G}$ par rapport à $\boldsymbol{h}$, il existe un système admissible $\{\gamma, \Delta\}$ dans $\phi$ tel que $\sigma_0$ soit la restriction de l'automorphisme intérieur $exp\ \pi\ i\ \text{ad}(h_\gamma^\Delta)$ de $\mathcal{G}$ à $\mathcal{G}_0$ ($\sigma_0$ n'est pas nécessairement intérieur). De plus, on peut trouver un nombre complexe $\lambda$ tel que

$$\underline{\Omega_c} = \lambda\ \underline{h_\gamma^\Delta}$$

Remarquons que, dans le complexifié, $\mathcal{Z}(\mathcal{K}) = \mathbb{C} \cdot h_\gamma^\Delta$ et $\mathcal{Z}(\mathcal{K}_0) = \mathcal{Z}(\mathcal{K}) \cap \mathcal{K}_0$; il existe donc un nombre complexe $\lambda$ tel que $\lambda\, h_\gamma^\Delta \in \mathcal{Z}(\mathcal{K}_0)$. En outre, la forme de Killing, $\beta_0$, de $\mathcal{G}_0$ est la restriction de la forme de Killing, $\beta$ de $\mathcal{G}$ à $\mathcal{G}_0 \times \mathcal{G}_0$. Dès lors, on a

$$\begin{aligned}\beta(\lambda\, h_\gamma^\Delta, \lambda\, h_\lambda^\Delta) &= \lambda^2\, ||h_\lambda^\Delta||^2_{\boldsymbol{h}_\mathbb{R}} \\ &= \beta_0(\lambda\, h_\gamma^\Delta, \lambda\, h_\gamma^\Delta) \in \mathbb{R}\end{aligned}$$

donc $\lambda^2 \in \mathbb{R}$ et on a les deux possibilités :

(a) $i\, h_\gamma^\Delta \in \mathcal{Z}(\mathcal{K}_0)$

(b) $h_\gamma^\Delta \in \mathcal{Z}(\mathcal{K}_0)$

ou de manière équivalente

(a) $\mathcal{Z}(\mathcal{K}_0)$ est compact

(b) $\mathcal{Z}(\mathcal{K}_0)$ est non-compact.

Si $c_0$ est la conjugaison de $\mathcal{G}$ relativement à $\mathcal{G}_0$, on a $c_0(h_\gamma^\Delta) = \pm h_\gamma^\Delta$.

Remarquons encore que dans le cas $(a)$, l'élément $\mathcal{J}$ de $\mathrm{End}(\mathcal{P}_0)$ donné par

$$\mathcal{J} = \mathrm{ad}(i\, h_\gamma^\Delta)\big|_{\mathcal{P}_0}$$

définit une structure complexe $\mathcal{K}_0$-invariante sur $\mathcal{P}_0$.

**Lemme 3.4.2.** *Soient $\mathcal{G}$ une algèbre de Lie simple complexe et $\boldsymbol{h}$ une sous-algèbre de Cartan de $\mathcal{G}$. Toute conjugaison $\tau$ de $\mathcal{G}$ qui stabilise $\boldsymbol{h}$ stabilise $\boldsymbol{h}_\mathbb{R}$.*

<u>Preuve</u>. Soit $\mathcal{G}_\tau$ la forme réelle de $\mathcal{G}$ associée à $\tau$ et soit l'application $\mathbb{C}$-antilinéaire :

$$\tau^\star : \boldsymbol{h}^\star \to \boldsymbol{h}^\star,\ (\tau^\star \varphi)(H) = \overline{\varphi(\tau(H))}$$

$\forall \varphi \in \boldsymbol{h}^\star,\ H \in \boldsymbol{h}$.
Alors, $\forall X \in \mathcal{G}_\alpha : [\tau(H), \tau(X)] = (\tau^\star \alpha)(\tau H)\tau(X)$ et donc $\phi = \tau^\star \phi$.
Comme $<z, z'> = \beta(z, \tau z')$ définit une structure pseudo-hermitienne sur $\mathcal{G}$, on a

$$\beta(\tau H_\alpha, H) = \overline{\beta(H_\alpha, \tau H)} = (\tau^\star \alpha)(H)$$

$\forall H \in \boldsymbol{h}$; dès lors $\tau H_\alpha = H_{\tau^\star \alpha}$. ∎

Réciproquement au théorème 3.4.1., on a le

**Théorème 3.4.3.** *Soient $\mathcal{G}$ une algèbre de Lie simple complexe, $\boldsymbol{h}$ une sous-algèbre de Cartan de $\mathcal{G}$ et $\phi$ le système de racines associé. Soit $\{\alpha, \Delta\}$ un système admissible dans $\phi$. Soit $\tau$ une conjugaison de $\mathcal{G}$ stabilisant $\boldsymbol{h}$ et telle que $\tau h_\alpha^\Delta = \pm h_\alpha^\Delta$.
Alors, l'automorphisme intérieur de $\mathcal{G}$ donné par $\sigma = \exp \pi i\, \mathrm{ad}(h_\alpha^\Delta)$ est involutif et se restreint à $\mathcal{G}_\tau$ en un automorphisme involutif $\sigma_\tau$ (non nécessairement intérieur). De plus, si $\mathcal{G}_\tau = \mathcal{K}_\tau \oplus \mathcal{P}_\tau$ est la décomposition canonique du couple symétrique $(\mathcal{G}_\tau, \sigma_\tau)$, on a*

$$\mathcal{Z}(\mathcal{K}_\tau) \neq \{0\}$$

<u>Preuve</u>. Posons $\phi_0 = \{\beta \in \phi \mid \beta(h_\alpha^\Delta) = 0\}$,

$$\phi_p = \phi \setminus \phi_0$$

Comme $\tau\, h_\alpha^\Delta = \pm h_\alpha^\Delta$, $\tau^\star \phi_0 = \phi_0$ et $\tau^\star \phi_p = \phi_p$, et dès lors $\tau\sigma = \sigma\tau$. Le reste est immédiat. ∎

Remarquons que dans le théorème 3.4.1., il existe $r$ dans $\mathbb{R}_0$ avec $\underline{\Omega} = r\, e^{\frac{i\pi}{4}(1+\varepsilon)}\, h_\gamma^\Delta$ ($\varepsilon = \pm 1$).

## 3.5 Structures pseudokähleriennes — E.S.S. compacts

**Définition 3.5.1.** *Un quadruple $\kappa = (\mathcal{G}, \sigma, \Omega, \mathcal{J})$ est dit pseudokählerien si*

(i) $(\mathcal{G}, \sigma, \Omega)$ *est un T.S.S.*

(ii) $\mathcal{J}$ *est un endormorphisme de $\mathcal{P}$ tel que*

- $\mathcal{J}^2 = -id$
- $\Omega(\mathcal{J}X, \mathcal{J}Y) = \Omega(X, Y)$ *pour tous $X, Y$ dans $\mathcal{P}$*
- $\mathcal{J}$ *commute avec l'action de $\mathcal{K}$ sur $\mathcal{P}$.*

*Nous dirons que $\kappa$ est kählerien si la forme symétrique sur $\mathcal{P}$, $g(X, Y) = \Omega(\mathcal{J}X, Y)$, est définie.*

**Proposition 3.5.2.** *Soit $\kappa = (\mathcal{G}, \sigma, \Omega, \mathcal{J})$ un quadruple pseudokählerien simple. Alors*

$$\mathcal{J} = \pm \left.\frac{1}{\pi} \log \sigma\right|_{\mathcal{P}}$$

*En particulier, $\mathcal{J}$ est unique au signe près et la métrique pseudokählerienne est multiple de la forme de Killing restreinte à $\mathcal{P}$.*

PREUVE. On sait [Ca-Pa] que l'endormorphisme, $D$, de $\mathcal{G}$ extension de $\mathcal{J}$ par 0 sur $\mathcal{K}$ est une dérivation de $\mathcal{G}$; $D$ est dès lors réalisé par l'action adjointe d'un élément $z$ du centre de $\mathcal{K}$. La condition $\mathcal{J}^2 = -id$ livre alors la proposition. ∎

**Corollaire 3.5.3.** *Un T.S.S. simple admet une structure pseudokählerienne si et seulement si le centre de $\mathcal{K}$ contient un élément compact.*

**Proposition 3.5.4.** *Un T.S.S. réductif indécomposable et non plat est simple.*

PREUVE. Soit $t = (\widetilde{\mathcal{G}}, \widetilde{\sigma}, \widetilde{\Omega})$ un T.S.S. Avec les notations du paragraphe 3.1., on a

$$\widetilde{\Omega}(\mathcal{P}, \mathcal{Z}) = \widetilde{\Omega}([\mathcal{K}\ \mathcal{P}], \mathcal{Z}) = \widetilde{\Omega}(\mathcal{P}, [\mathcal{K}\ \mathcal{Z}]) = 0$$

∎

Tout idéal d'une algèbre réductive étant réductif, ceci classifie les T.S.S. réductifs comme sommes directes de T.S.S. simples et de facteurs plats.

**Théorème 3.5.5.** *Un E.S.S. compact indécomposable et non plat est simple.*

PREUVE. Par [Koh], tout espace symétrique connexe compact admet un groupe de transvections compact, donc réductif. ∎

**Corollaire 3.5.6.** *Tout E.S.S. compact est kählerien.*

## 3.6 Equivalences

Rappelons que si $q = (\mathcal{G}, \sigma)$ est un couple symétrique simple, il existe une involution de Cartan $\theta$ qui commute avec $\sigma$. En effet, soit $\widetilde{\theta}$ une involution de Cartan de $\mathcal{G}$ ($\widetilde{\theta} = id$ si $\mathcal{G}$ est compacte). Alors,

$$<X, Y> = \beta(X, \widetilde{\theta} Y)$$

définit un produit scalaire défini négatif sur $\mathcal{G}$. Soit alors $A = \sigma\widetilde{\theta}$ et $U = A^2$, on a :

$$<AX, Y> = \beta(\sigma\widetilde{\theta}X, \widetilde{\theta}Y) = \beta(\widetilde{\theta}X, \sigma\widetilde{\theta}Y) = <AY, X> = <X, AY>$$

On peut dès lors choisir une base $<,>$-orthonormée $\{X_i\}$ de $\mathcal{G}$ telle que
$U = diag[\lambda_1, \ldots, \lambda_n]$ avec $\lambda_i > 0$. Soit $\varphi_t = diag[\lambda_1^t, \ldots, \lambda_n^t]$, $t \in \mathbb{R}$.
Alors, on a :
$$[U\ X_i, U\ X_j] = C_{ij}^k\ \lambda_k X_k$$

c.à.d.
$$\lambda_i \lambda_j C_{ij}^k = C_{ij}^k \lambda_k$$

(sans sommation). Dès lors, $\forall t \in \mathbb{R}$ :
$$\lambda_i^t \lambda_j^t C_{ij}^t = C_{ij}^t\ \lambda_k^t$$

$\varphi_t$ est donc un groupe à un paramètre d'automorphismes de $\mathcal{G}$ tel que $\varphi_1 = U$.
Maintenant, $\varphi_t = e^{t\ \log U}$, donc,

$$\begin{aligned}\widetilde{\theta}\varphi_t\widetilde{\theta} &= \sum_{k=0}^\infty \frac{t^k}{k!}\widetilde{\theta}(\log U)^k\widetilde{\theta} \\ &= \sum_{k=0}^\infty \frac{t^k}{k!}(\widetilde{\theta}\log U\widetilde{\theta})^k \\ &= e^{t\widetilde{\theta}\log U\widetilde{\theta}}\end{aligned}$$

mais $\widetilde{\theta}A\widetilde{\theta} = A^{-1}$ donc $\widetilde{\theta}U\widetilde{\theta} = U^{-1}$, donc

$$\widetilde{\theta}e^{\log U}\widetilde{\theta} = e^{\widetilde{\theta}\log U\widetilde{\theta}} = e^{-\log U}$$

et dès lors
$$\widetilde{\theta}\varphi_t\widetilde{\theta} = \varphi_{-t}$$

Posons $\Theta_t = \varphi_t\widetilde{\theta}\varphi_{-t}$; comme $AU^{-1} = A^{-1}$ on a

$$\sigma\theta_t\ =\ \sigma\varphi_t\widetilde{\theta}\varphi_{-t}\ =\ \sigma\widetilde{\theta}\varphi_{-2t}\ =\ A\varphi_{-2t}$$

et
$$\begin{aligned}\theta_t\sigma\ =\ \varphi_t\widetilde{\theta}\varphi_{-t}\sigma\ =\ \varphi_{2t}\widetilde{\theta}\sigma\ &=\ \varphi_{2t}A^{-1} \\ &=\ A^{-1}\varphi_{2t} \\ &=\ A\varphi_{2t-1}\end{aligned}$$

Dès lors, $\sigma\theta_t - \theta_t\sigma = A(\varphi_{-2t} - \varphi_{2t-1})$.
Il suffit donc de choisir $\theta = \theta_{1/4}$

Rappelons encore que si $\mathcal{G} = \mathcal{K} \oplus \mathcal{P}$ est la décomposition relativement à $\sigma$, la restriction de $\theta$ à $\mathcal{K}$ ($[\theta, \sigma] = 0$) est encore une involution de Cartan de l'algèbre réductive $\mathcal{K}$; il en est de même pour sa restriction à l'algèbre semi-simple $[\mathcal{K}, \mathcal{K}]$.
On adoptera dans tout ce paragraphe les notations suivantes :

- $(\mathcal{G}_0, \sigma_0)$ est un C.S. absolument simple de décomposition canonique $\mathcal{G}_0 = \mathcal{K}_0 \oplus \mathcal{P}_0$ dont on suppose $\mathcal{Z}(\mathcal{K}_0) \neq \{0\}$.

- $\mathcal{G} = \mathcal{G}_0 \otimes \mathbb{C}$ en est l'algèbre (simple) complexifiée et $c_0$ est la conjugaison de $\mathcal{G}$ définissant $\mathcal{G}_0$.

- $\sigma$ est l'extension $\mathbb{C}$-linéaire de $\sigma_0$ à $\mathcal{G}$.

- $\theta_0$ est une involution de Cartan de $\mathcal{G}_0$ qui commute avec $\sigma_0$ et $\mathcal{G}_0 = \mathcal{K}_c \oplus \mathcal{P}_n$ est la décomposition de Cartan associée (Si $\mathcal{G}_0$ est compacte, $\theta_0 = id$ et $\mathcal{G}_0 = \mathcal{K}_c$).

- $\mathcal{G}_u = \mathcal{K}_c \oplus i\mathcal{P}_n$ est la forme réelle compacte de $\mathcal{G}$ définie par la conjugaison $c_u$ de $\mathcal{G}$, extension $\mathbb{C}$-antilinéaire de $\theta_o$ à $\mathcal{G}$.

- $\boldsymbol{h}_0$ est une sous-algèbre de Cartan de $\mathcal{G}_0$ contenue dans $\mathcal{K}_0$ et $\theta_0$-stable.

- $\boldsymbol{h} = \boldsymbol{h}_0 \otimes \mathbb{C}$ est la sous-algèbre de Cartan de $\mathcal{G}$, complexifiée de $\boldsymbol{h}_0$ et $\phi$ en est le système de racines associé.

- $\forall \alpha \in \phi$, $H_\alpha$ est l'unique élément de $\boldsymbol{h}$ tel que $\beta(H_\alpha, H) = \alpha(H)$ $\forall H \in \boldsymbol{h}$.

- $\boldsymbol{h}_\mathbb{R} = \sum_{\alpha \in \phi} \mathbb{R} \cdot H_\alpha$ et $i \, \boldsymbol{h}_\mathbb{R} = \boldsymbol{h}_u \subset \mathcal{G}_u$.

- $\mathcal{G} = \boldsymbol{h} \oplus \bigoplus_{\alpha \in \phi} \mathcal{G}_\alpha$ est la décomposition radicielle de $\mathcal{G}$ par rapport à $\boldsymbol{h}$.

- $(\,,\,)$ est la structure euclidienne induite par $\beta$ sur $\boldsymbol{h}_\mathbb{R}^\star$.

- $N(\boldsymbol{h}) = \{g \in Aut(\mathcal{G}) \mid g\boldsymbol{h} \subset \boldsymbol{h}\}$.

- $Aut(\phi)$ est le groupe d'isométries de $\boldsymbol{h}_\mathbb{R}^\star$ qui conservent $\phi$.

- $W$ est le groupe de Weyl de $\phi$.

- $\Gamma$ est le groupe d'automorphismes du diagramme de Dinkin de $\mathcal{G}$.

- Si $\tau$ est une conjugaison de $\mathcal{G}$ telle que $\tau\boldsymbol{h} \subset \boldsymbol{h}$, on note
$$\widetilde{N}_\tau = \{g \in N(\boldsymbol{h}) \mid g\tau = \tau g\}$$

- On note $\pi : N(\boldsymbol{h}) \to Aut(\phi)$ l'homomorphisme canonique; on rappelle que $\pi$ est un épimorphisme et on note $\pi\widetilde{N}_\tau = N_\tau$.

- Si $\alpha \in \phi$ et si $\Delta$ est une base de $\phi$ contenant $\alpha$, $\omega_\alpha^\Delta$ est le poid fondamental associé à $\alpha$ relativement à $\Delta$.

- Rappelons que pour tout $\alpha$ dans $\phi$, il existe $E_\alpha$ dans $\mathcal{G}_\alpha$ tel que

    (i) $i(E_\alpha - E_{-\alpha})$ et $E_\alpha + E_{-\alpha}$ sont dans $\mathcal{G}_u$.
    (ii) $[E_\alpha, E_{-\alpha}] = -\frac{2}{\alpha(H_\alpha)} H_\alpha$.
    (iii) $c_u E_\alpha = E_{-\alpha}$.
    (iv) $\beta(E_\alpha, E-\alpha) = -1$.

Comme corollaire du lemme 3.4.2., on a le

**Lemme 3.6.1.**

(i) $\boldsymbol{h}_{\mathbb{R}}$ est $c_0$ et $c_u$-stable.

(ii) Si on note $\boldsymbol{h}_{\mathbb{R}} = \boldsymbol{h}_{\mathbb{R}}^+ \oplus \boldsymbol{h}_{\mathbb{R}}^-$ avec
$$c_0|_{\boldsymbol{h}_{\mathbb{R}}} = id|_{\boldsymbol{h}_{\mathbb{R}}^+} \oplus (-id)|_{\boldsymbol{h}_{\mathbb{R}}^-}$$
et
$$\boldsymbol{h}_0 = \boldsymbol{h}_0^+ \oplus \boldsymbol{h}_0^-$$
avec
$$\theta_0|_{\boldsymbol{h}_0} = id|_{\boldsymbol{h}_0^+} \oplus (-id)|_{\boldsymbol{h}_0^-}$$
on a
$$\boldsymbol{h}_0^+ = i\boldsymbol{h}_{\mathbb{R}}^-$$
$$\boldsymbol{h}_0^- = \boldsymbol{h}_{\mathbb{R}}^+$$

**Remarque 3.6.2.**

- $\forall \alpha \in \phi$ et $H_r \in \boldsymbol{h}_{\mathbb{R}} : \alpha(H_r) \in \mathbb{R}$, donc, on a
$$\boldsymbol{h}_0^+ = \{H \in \boldsymbol{h}_0 \mid \mathrm{ad}(H) \text{ a ses valeurs propres imaginaires pures}\}$$
et
$$\boldsymbol{h}_0^- = \{H \in \boldsymbol{h}_0 \mid \mathrm{ad}(H) \text{ a ses valeurs propres réelles}\}$$
Ceci montre que deux involutions de Cartan de $\mathcal{G}_0$ qui stabilisent $\boldsymbol{h}_0$ induisent les mêmes décompositions : $\boldsymbol{h}_0 = \boldsymbol{h}_0^+ \oplus \boldsymbol{h}_0^-$.

- Si $\varphi \in Aut(\mathcal{G})$ avec $\varphi \boldsymbol{h} \subset \boldsymbol{h}$, alors
$$^\tau\varphi : \boldsymbol{h}^\star \to \boldsymbol{h}^\star$$
donné par
$$(^\tau\varphi.\eta)(H) = \eta(\varphi^{-1}(H))$$
induit l'élément $\pi(\varphi) = {^\tau_\varphi}\big|_\phi \in Aut(\phi)$.

- $\mathcal{Z}(\mathcal{K}_0)$ est contenu, soit dans $\boldsymbol{h}_0^+$, soit dans $\boldsymbol{h}_0^-$

**Définition 3.6.3.** *Une sous-algèbre de Cartan $\boldsymbol{h}_0$ de $\mathcal{G}_0$ sera dite toroïdale (resp. Iwasawa) si $\boldsymbol{h}_0^+$ (resp. $\boldsymbol{h}_0^-$) a la plus grande dimension possible.*

**Remarque 3.6.4.** *Deux sous-algèbres de Cartan toroïdales (resp. Iwasawa) sont conjuguées sous l'action du groupe adjoint [Su1, page 415].*

**Proposition 3.6.5.** *Si $\mathcal{Z}(\mathcal{K}_0) \subset \mathcal{K}_c$ (resp. $\mathcal{Z}(\mathcal{K}_0) \subset \mathcal{P}_n$), alors on peut supposer $\boldsymbol{h}_0$ toroïdale (resp. Iwasawa).*

PREUVE. On pose
$$\begin{array}{rcl}\mathcal{K}_0^c & = & \mathcal{K}_0 \cap \mathcal{K}_c \\ \mathcal{K}_0^n & = & \mathcal{K}_0 \cap \mathcal{P}_n \\ \mathcal{P}_0^c & = & \mathcal{P}_0 \cap \mathcal{K}_c \\ \mathcal{P}_0^n & = & \mathcal{P}_0 \cap \mathcal{P}_n\end{array}$$
Comme $[\theta_0, \sigma_0] = 0$, on a $\mathcal{K}_0 = \mathcal{K}_0^c \oplus \mathcal{K}_0^n$ et $\mathcal{P}_0 = \mathcal{P}_0^c \oplus \mathcal{P}_0^n$. Supposons $\mathcal{Z}(\mathcal{K}_0) \subset \mathcal{K}_c$ alors, $\mathcal{Z}(\mathcal{K}_0) \subset \mathcal{K}_0^c$. Soit $\boldsymbol{h}_c$ une sous-algèbre de Cartan de $\mathcal{K}_0^c$, alors, comme $[\mathcal{Z}(\mathcal{K}_0), \mathcal{K}_0] = 0$, on a
$$\boldsymbol{h}_c \supseteq \mathcal{Z}(\mathcal{K}_0)$$

Maintenant, $C_{\mathcal{G}_0}(\boldsymbol{h}_c) \subset C_{\mathcal{G}_0}(\mathcal{Z}(\mathcal{K}_0)) = \mathcal{K}_0$, dès lors, si $\mathcal{B}$ est maximale abélienne dans $C_{\mathcal{K}_0^n}(\boldsymbol{h}_c)$, $\boldsymbol{h}_0 = \boldsymbol{h}_c \oplus \mathcal{B}$ est une sous-algèbre abélienne maximale $\theta_0$-stable dans $\mathcal{K}_0$, donc de Cartan dans $\mathcal{K}_0$ et dans $\mathcal{G}_0$. Soit alors $p \in C_{\mathcal{P}_0^c}(\boldsymbol{h}_c)$ alors

$$[\mathcal{Z}(\mathcal{K}_0), p] = 0$$

et donc $p = 0$; $\boldsymbol{h}_c$ est donc un tore maximal de $\mathcal{K}_c$.

Supposons maintenant $\mathcal{Z}(\mathcal{K}_0) \subset \mathcal{P}_n$ et notons $\mathcal{D}_0 = [\mathcal{K}_0\,\mathcal{K}_0]$. $\theta_0|_{\mathcal{D}_0}$ est alors une involution de Cartan de $\mathcal{D}_0$ dont la décomposition de Cartan associée est

$$\mathcal{D}_0 = \mathcal{D}_0 \cap \mathcal{K}_c \oplus \mathcal{D}_0 \cap \mathcal{P}_n$$

Soit $\boldsymbol{a}_1$ une sous-algèbre maximale abélienne dans $\mathcal{D}_0 \cap \mathcal{P}_n$ et posons

$$\boldsymbol{a}_0 = \boldsymbol{a}_1 \oplus \mathcal{Z}(\mathcal{K}_0) \,;$$

$\boldsymbol{a}_0$ est maximale abélienne dans $\mathcal{K}_0^n$. Soit $c \in C_{\mathcal{P}_n}(\boldsymbol{a}_0)$ et écrivons $c = c_{\mathcal{K}_0} \oplus c_{\mathcal{P}_0}$ sa décomposition relativement à $\sigma_0$; alors, $[c, \mathcal{Z}(\mathcal{K}_0)] = 0$ livre $c_{\mathcal{P}_0} = 0$ et donc $c \in \boldsymbol{a}_0$ par maximalité dans $\mathcal{K}_0^n$. Dès lors, $\boldsymbol{a}_0$ est maximale abélienne dans $\mathcal{P}_n$. ■

**Lemme 3.6.6.** *Soient $\mathcal{D}_0$ une sous-algèbre semi-simple de $\mathcal{G}_0$ et $\boldsymbol{h}_1$ une sous-algèbre de Cartan de $\mathcal{D}_0$. Si $\boldsymbol{h}_0 = C_{\mathcal{G}_0}(\mathcal{D}_0) \oplus \boldsymbol{h}_1$ est une sous-algèbre de Cartan du type toroïdal (resp. Iwasawa) de $\mathcal{G}_0$, il en est de même pour $\boldsymbol{h}_1$ dans $\mathcal{D}_0$.*

PREUVE. Soit $\theta_0$ une involution de Cartan de $\mathcal{G}_0$ stabilisant $\mathcal{D}_0$ et soit $\theta_1$ sa restriction à $\mathcal{D}_0$ ($\theta_1$ est une involution de Cartan de $\mathcal{D}_0$). $\boldsymbol{h}_1$ est conjuguée dans $Int(\mathcal{G}_0)$ (et donc dans $Aut(\mathcal{G}_0)$) à une sous-algèbre de Cartan $\boldsymbol{h}_1^{(0)}$ $\theta_1$-stable. Comme $\theta_0 C_{\mathcal{G}_0}(\mathcal{D}_0) = C_{\mathcal{G}_0}(\mathcal{D}_0)$, $\boldsymbol{h}_0^{(0)} = C_{\mathcal{G}_0}(\mathcal{D}_0) \oplus \boldsymbol{h}_1^{(0)}$ est une sous-algèbre de Cartan de type toroïdal (resp. Iwasawa) $\theta_0$ stable dans $\mathcal{G}_0$.

(a) $\boldsymbol{h}_0$ est toroïdale.

Soit $k_0 \in \mathcal{D}_0$ tel que $\theta_1(k_0) = k_0$ et $[k_0, \boldsymbol{h}_1^{(0)+}] = 0$ où $\boldsymbol{h}_1^{(0)} = \boldsymbol{h}_1^{(0)+} \oplus \boldsymbol{h}_1^{(0)-}$ est la décomposition induite par $\theta_1$. Alors, $[k_0, \boldsymbol{h}_0^{(0)+}] = 0$ où $\boldsymbol{h}_0^{(0)} = \boldsymbol{h}_0^{(0)+} \oplus \boldsymbol{h}_0^{(0)-}$ est la décomposition induite par $\theta_0$. Comme $\boldsymbol{h}_0^{(0)}$ est toroïdale, $k_0 \in \boldsymbol{h}_0^{(0)+}$ c.à.d.

$$k_0 \in \boldsymbol{h}_0^{(0)+} \cap \mathcal{D}_0 = \boldsymbol{h}_1^{(0)+}$$

(b) $\boldsymbol{h}_0$ est Iwasawa.

Soit $p_0 \in \mathcal{D}_0$ tel que $\theta_1(p_0) = -p_0$ et $[p_0, \boldsymbol{h}_1^{(0)-}] = 0$. Alors un raisonnement analogue à celui tenu en (a) livre

$$p_0 \in \boldsymbol{h}_0^{(0)-}$$

■

**Lemme 3.6.7.** *Soit $\tau$ une conjugaison de $\mathcal{G}$ stabilisant $\boldsymbol{h}$. Soit $\rho$ dans $Aut(\phi)$ tel que $[\rho, \tau^\star] = 0$. Posons, pour tout $\alpha$ dans $\phi$,*

$$\tau E_\alpha = \nu_\alpha E_{\tau^\star(\alpha)} \qquad\qquad (\nu_\alpha \in \mathbb{C})$$

*Soit $\Delta = \{\alpha_1, \alpha_2, \ldots, \alpha_\ell\}$ une base de $\phi$. Alors, $\rho$ est un élément de $N_\tau$ si et seulement si le système d'équations*

$$z_{\alpha_i} = \overline{\nu_{\alpha_i}}\ \overline{z_{\tau^\star(\alpha_i)}}\ \nu_{\rho\tau^\star(\alpha_i)} \qquad (i = 1, 2, \ldots, \ell) \qquad (1)$$

*admet une solution du type*

$$(z_{\alpha_1}, \ldots, z_{\alpha_\ell}, z_{\tau^\star(\alpha_1)}, \ldots, z_{\tau^\star(\alpha_\ell)})$$

dans $\mathbb{C}_0^\ell \times \mathbb{C}^\ell$.

PREUVE. Il est bien connu ([Mu1,2,3, Su1]) que l'on peut supposer, pour tous $\alpha, \beta \in \phi$ :

$$[E_\alpha, E_\beta] = N_{\alpha,\beta} E_{\alpha+\beta}$$

avec

- $N_{\alpha,\beta} \neq 0 \iff \alpha + \beta \in \phi$,
- $N_{\alpha,\beta} \in \mathbb{R}$
- $N_{\alpha,\beta} = N_{-\alpha,-\beta}$

$\Delta^\rho = \{\rho(\alpha_1), \ldots, \rho(\alpha_\ell)\}$ est une base de $\phi$ et pour tout choix de $\underline{z} = (z_1, \ldots, z_\ell) \in \mathbb{C}_0^\ell$, il existe un unique automorphisme $\psi_{\underline{z}}$ de $\mathcal{G}$ tel que

$$\psi_{\underline{z}} \in N(\boldsymbol{h}); \quad \pi\psi_{\underline{z}} = \rho \quad \text{et} \quad \psi_{\underline{z}} E_{\alpha_i} = z_i\, E_{\rho(\alpha_i)}$$

Posons $\psi_{\underline{z}} E_\alpha = z_\alpha E_{\rho(\alpha)} \quad \forall \alpha \in \phi$.

On a $\psi_{\underline{z}}[E_\alpha, E_{-\alpha}] = z_\alpha \cdot z_{-\alpha} \left(\frac{-2}{\rho(\alpha)(H_{\rho(\alpha)})}\right) H_{\rho(\alpha)}$ et $\beta(\psi_{\underline{z}} H_\alpha,\, H) = \rho(\alpha)\,(H)$; ceci livre

$$\psi_{\underline{z}} H_\alpha = H_{\rho(\alpha)} \quad \text{et } z_\alpha\,\cdot\,z_{-\alpha} = 1 \tag{a}$$

De plus, $\psi_{\underline{z}}[E_\alpha, E_\beta] = N_{\alpha,\beta}\, \psi_{\underline{z}}(E_{\alpha+\beta})$ c'est-à-dire

$$z_\alpha\, z_\beta N_{\rho(\alpha),\rho(\beta)} = N_{\alpha,\beta}\, z_{\alpha+\beta} \tag{b}$$

Maintenant, $[\psi_{\underline{z}}, \tau] = 0$ si et seulement si pour tout $\alpha \in \phi$ : $(\tau\psi_{\underline{z}}\tau - \psi_{\underline{z}})E_\alpha = 0$
c'est-à-dire

$$\overline{\nu_\alpha z_{\tau^\star(\alpha)}}\nu_{\rho\tau^\star(\alpha)} = z_\alpha \tag{c}$$

On vérifie en utilisant $(a)$ et $(b)$, que si $(c)$ est vrai pour $\alpha$ et $\beta$ dans $\phi$, alors $(c)$ est vrai pour $\alpha + \beta$ et $-\alpha$ (si $\alpha + \beta \in \phi$).

Rappelons que pour tout $\alpha \in \phi$ :

$$\alpha = \alpha_1 + \ldots + \alpha_k \qquad (k \in \{1, \ldots, \ell\})$$

avec toute somme partielle

$$\alpha_1 + \ldots + \alpha_i \quad \text{dans } \phi \qquad [Hu,\ pg\ 50]$$

Dès lors, $\rho \in N_\tau$ si et seulement si il existe $\underline{z}$ dans $\mathbb{C}_0^\ell$ avec $\psi_{\underline{z}}$ tel que

$$z_i = \overline{\nu_{\alpha_i}}\ \overline{z_{\tau^\star(\alpha_i)}}\nu_{\rho\tau^\star(\alpha_i)}$$

∎

**Remarque**. Les conjugaisons d'une algèbre simple complexe sont très explicitement décrite dans [BodeSi, Mu1]. Ceci permet une description aisée du groupe $N_\tau$ pour toute conjugaison $\tau$.

**Corollaire** . *Si $\tau$ est une conjugaison de $\mathcal{G}$ stabilisant $\boldsymbol{h}$ et commutant avec $c_u$, alors $-id|_{\boldsymbol{h}_\mathbb{R}^\star}$ est un élément de $N_\tau$.*

PREUVE. Soit $A$ l'automorphisme de $\mathcal{G}_u$, restriction de $\tau$ à $\mathcal{G}_u$.

(a) Supposons que $A = exp\ \pi i\text{ad}(h_j)$, où $h_j \in \boldsymbol{h}_\mathbb{R}$ est donné par $\alpha_i(h_j) = \delta_{ij}$ pour $i = 1, \ldots, \ell$.
Alors, comme $c_u\ E_\alpha = E_{-\alpha}$, on a
$$U_\alpha = E_\alpha + E_{-\alpha} \in \mathcal{G}_u$$
et
$$V_\alpha = E_\alpha - E_{-\alpha} \in i\mathcal{G}_u$$
Tout élément $Z \in \mathcal{G}$ s'écrit $Z = X + iY$ avec $X, Y \in \mathcal{G}_u$ et $\tau$ est alors donné par $\tau(Z) = AX - iAY$ ou, ce qui revient au même, en notant $A^\mathbb{C}$ l'extension $\mathbb{C}$-linéaire de $A$ de $\mathcal{G}$ :
$$\tau(Z) = A^\mathbb{C} X - A^\mathbb{C}(iY)$$
dès lors, $\tau(E_\alpha) = \frac{1}{2}\tau(U_\alpha + V_\alpha) = A^\mathbb{C} E_{-\alpha}$; de là on tire
$$\nu_\alpha = e^{-\pi i \alpha(h_j)}$$
Donc $\nu_{\alpha_i} = (-1)^{\delta_{ij}}$ et le système (1), pour $\rho = -id$, prend la forme
$$z_{\alpha_i} = \overline{z_{\alpha_i}}$$

(b) Supposons que $A = \Theta$ où $\Theta$ induit un automorphisme $\theta$ du diagramme de Dinkin de $\mathcal{G}$ et est tel que $\Theta^\mathbb{C} E_{\alpha_i} = E_{\theta(\alpha_i)}\ \forall \alpha_i \in \Delta$.
Dans ce cas, le système (1) devient,
$$z_{\alpha_i} = \overline{z_{\theta(\alpha_i)}}$$
qui admet la solution $z_{\alpha_i} = 1\ \ \forall i$.

(c) $A = \Theta\ exp\ \pi i\text{ad}(\widetilde{h}_j)$, $1 \leqslant j \leqslant p \leqslant \ell$ où $\Theta \widetilde{h}_j \widetilde{h}_j$ est comme dans (b) et $\widetilde{h}_j$ est l'élement de $\boldsymbol{h}_\mathbb{R}$ donné par
$$\begin{aligned}\alpha_i(\widetilde{h}_j) &= \delta_{ij} & \forall i = 1; \ldots, p \\ \xi_k(\widetilde{h}_j) &= \xi_k^\star(\widetilde{h}_j) = 0 & \forall k = 1, \ldots, r\end{aligned}$$
où on écrit $\Delta$ sous la forme
$$\Delta = \{\alpha_1, \ldots, \alpha_p, \xi_1, \xi_1^\star, \ldots, \xi_r, \xi_r^\star\}$$
avec
$$\begin{aligned}\theta\alpha_i &= \alpha_i & \forall i = 1, \ldots, p \\ \theta\xi_k &= \xi_k^\star & \forall k = 1, \ldots, r\end{aligned}$$
Dans ce cas,
$$\begin{aligned}\nu_{\alpha_s} &= e^{-\pi\ i\alpha_s\ (\widetilde{h_j})} \\ &= (-1)^{\delta_{js}}\end{aligned} \qquad \forall s = 1, \ldots, \ell$$
et on se ramène au cas (b).

(d) On sait [Mu1,2,3] que $\tau|_{\mathcal{G}_u} = 1$ est conjugué par un élément $h \in exp\ \boldsymbol{h}_u$, à un élément $A_0$ d'une des formes décrites dans $(a)$ $(b)$ ou $(c)$. Dès lors, si $\tau_0$ désigne la conjugaison de $\mathcal{G}$ déduite de $A_0$ on a
$$\tau = h\ \tau_0\ h^{-1}$$
On sait qu'il existe $\varphi_0 \in \widetilde{N_{\tau_0}}$ tel que $\pi\varphi_0 = -id$ dès lors
$$[h\ \varphi_0\ h^{-1},\ \tau] = 0$$

et
$$\pi(h\ \varphi_0\ h^{-1}) = -id$$
donc
$$-id|_{\boldsymbol{h}_{\mathbb{R}}^\star} \in N_\tau$$

■

**Théorème 3.6.9.** *Supposons $\boldsymbol{h}_0$ du type toroïdal ou Iwasawa. Soient $\{\Delta_j, \gamma_j\}$ $j = 1, 2$ deux systèmes admissibles dans $\phi$ tels que $c_0(h_{\gamma_j}^{\Delta_j}) = \pm h_{\gamma_j}^{\Delta_j}$ ($j = 1$ ou $2$).*
*Soient $\sigma_j = exp\ \pi\ i\ ad(h_{\gamma_j}^{\Delta_j})\Big|_{\mathcal{G}_0}$ ($j = 1$ ou $2$) les automorphismes involutifs de $\mathcal{G}_0$ associés (cf. théorème 3.4.3.).*
*Alors, les couples symétriques $(\mathcal{G}_0, \sigma_i)$ ($i = 1, 2$) sont isomorphes si et seulement si les poids fondamentaux $\omega_{\gamma_i}^{\Delta_i}$ ($i = 1, 2$) sont dans la même orbite pour l'action de $N_{c_0}$ sur $\boldsymbol{h}_{\mathbb{R}}^\star$.*

PREUVE. Notons $\mathcal{G}_0 = \mathcal{K}_i \oplus \mathcal{P}_i$ la décomposition relativement à $\sigma_i$ $i = 1$ ou $2$. Alors, par construction de $\mathcal{K}_i$ (cf. théorème 3.4.3.) et par le Lemme 3.6.6., on a

(a) $\boldsymbol{h}_0$ est contenu dans $\mathcal{K}_1 \cap \mathcal{K}_2$ et est de Cartan dans $\mathcal{K}_i$ ($i = 1, 2$).

(b) $\boldsymbol{h}_i = \boldsymbol{h}_0 \cap [\mathcal{K}_i, \mathcal{K}_i]$ ($i = 1, 2$) est une sous-algèbre toroïdale ou Iwasawa dans $\mathcal{D}_i = [\mathcal{K}_i, \mathcal{K}_i]$. Supposons les deux couples isomorphes. Il existe alors $\varphi_0 \in Aut(\mathcal{G}_0)$ tel que $\sigma_2 = \varphi_0 \sigma_1 \varphi_0^{-1}$; donc $\mathcal{K}_2 = \varphi_0 \mathcal{K}_1$.
Soit $\boldsymbol{h}_0' = \varphi_0(\boldsymbol{h}_0)$, c'est une sous-algèbre de Cartan de $\mathcal{K}_2$ avec $\boldsymbol{h}_2' = \varphi_0(\boldsymbol{h}_1)$ toroïdale ou Iwasawa dans $\mathcal{D}_2$.
De plus, $\varphi_0(\mathcal{Z}(\mathcal{K}_1)) = \mathcal{Z}(\mathcal{K}_2)$.
Comme $\boldsymbol{h}_2$ et $\boldsymbol{h}_2'$ sont de même type dans $\mathcal{D}_2$, il existe $k_2 \in Int(\mathcal{G}_2)$ avec $k_2 \cdot \boldsymbol{h}_2' = \boldsymbol{h}_2$; dès lors comme $Int(\mathcal{K}_2) \circ \sigma_2 = \sigma_2$, on peut supposer $\varphi_0(\boldsymbol{h}_1) = \boldsymbol{h}_2$ et donc $\varphi_0(\boldsymbol{h}_0) = \boldsymbol{h}_0$.
Soit $\varphi$, l'extension $\mathbb{C}$-linéaire de $\varphi_0$ à $\mathcal{G}$. On a $\varphi(\boldsymbol{h}) = \boldsymbol{h}$ et $\varphi^\star \phi = \phi$; donc, comme $[\varphi, c_0] = 0$, $\rho = \pi(\varphi)$ est dans $N_{c_0}$. De plus, comme dim $\mathcal{Z}(\mathcal{K}_i) = 1$ et que $\varphi$ stabilise $\boldsymbol{h}_{\mathbb{R}}$, on a $\varphi\ h_{\gamma_1}^{\Delta_1} = \pm h_{\gamma_2}^{\Delta_2}$; c'est-à-dire, par $\beta$-dualité,
$$\rho \omega_{\gamma_1}^{\Delta_1} = \pm \omega_{\gamma_2}^{\Delta_2}$$
Comme $-id|_{\boldsymbol{h}_{\mathbb{R}}^\star}$ est dans $N_{c_0}$ on peut supposer
$$\rho \omega_{\gamma_1}^{\Delta_1} = \omega_{\gamma_2}^{\Delta_2}$$

La réciproque est immédiate.

■

**Définition 3.6.10.** *Dans $\boldsymbol{h}_{\mathbb{R}}^\star$, soient les ensembles finis de points*
$$\begin{aligned}
\square_\phi &= \{\omega_\alpha^\Delta \mid \{\alpha, \Delta\}\ \text{admissible}\} \\
\square_\phi^c &= \{x \in \square_\phi \mid x = -c_0^\star\ x\} \\
\square_\phi^n &= \{y \in \square_\phi \mid y = c_0^\star y\}
\end{aligned}$$

Alors, on a immédiatement la

**Proposition 3.6.11.**

(i) $Aut\ \phi$ agit sur $\square_\phi$.

(ii) $N_{c_0}$ agit sur $\square_\phi^c$.

(iii) $N_{c_0}$ agit sur $\square_\phi^n$.

## 3.7 Classifications

**Théorème 3.7.1.**

(a) $\mathcal{G}_0$ est l'algèbre sous-jacente à un T.S.S. si et seulement si $\mathcal{G}$ est du type $A, B, C, D, E_6$ ou $E_7$. Dans ce cas, l'automorphisme $\sigma_0$ est la restriction d'un automorphisme intérieur $\sigma$ de $\mathcal{G}$.

(b) En choisissant $\boldsymbol{h}_0$ toroïdale, on a que l'ensemble des classes d'isomorphie des couples symétriques sous-jacents aux T.S.S. pseudo-kähleriens est en bijection avec l'ensemble des orbites

$$\square_\phi^c / N_{c_0} \ ;$$

la structure complexe $\mathcal{J}$ est donnée par

$$\mathcal{J} = \pm \frac{1}{\pi} \log \sigma \bigg|_{\mathcal{P}_0}$$

(c) En choisissant $\boldsymbol{h}_0$ Iwasawa, on a que l'ensemble des classes d'isomorphie des couples symétriques sous-jacents aux T.S.S. non-pseudo-kähleriens est en bijection avec l'ensemble des orbites

$$\square_\phi^n / N_{c_0}$$

(d) L'ensemble des classes d'isomorphie des couples symétriques sous-jacents aux T.S.S. simples non-absolument simples du type $t = (\mathcal{G}, \sigma, \Omega)$ est paramétrisé par l'ensemble des orbites :

$$\square_\phi / Aut\ \phi$$

il est en bijection avec celui des T.S.S. simples compacts; en particulier les T.S.S. simples non-absolument simples sont pseudo-kähleriens et la structure complexe $\mathcal{J}$ est donnée par

$$\mathcal{J} = \pm \frac{1}{\pi} \log \sigma \bigg|_{\mathcal{P}}$$

**Remarques**

(i) Un T.S.S. simple $(\mathcal{G}_1, \sigma_1, \Omega_1)$ est pseudo-kählerien si et seulement si $\mathcal{Z}(\mathcal{K}_1)$ contient un élément compact.

(ii) Si $\mathcal{G}_0 = \mathcal{G}_u$ est compacte, $c_0 = c_u$, dès lors

$$\square_\phi = \square_\phi^c \qquad \text{et} \qquad N_{c_u} = Aut\ \phi$$

(iii) L'élément $\xi \in \mathcal{G}_0^\star$ (resp. $\mathcal{G}^\star$) tel que $\delta\xi = \underline{\Omega_0}$ (resp. $\underline{\Omega}$) est un multiple réel (resp. complexe) de

- $i\omega_\alpha^\Delta$ dans le cas pseudo-kählerien
- $\omega_\alpha^\Delta$ dans le cas non-pseudo-kählerien.

Deux multiples différents en valeur absolue donnent lieu à des T.S.S. non isomorphes.

En choisissant dans la liste de Berger (Tableau II [Be]) les couples symétriques dont l'holonomie linéaire admet un centre, on obtient une liste des T.S.S. simples.

Dans les tableaux suivants, $so(2)$ signifie que $\mathcal{Z}(\mathcal{K}_0)$ est compact, $\mathbb{R}$ non compact.
Les nombres $p$ et $q$ peuvent prendre la valeur 0.

    ***TABLEAU A***   (pseudo-kähleriens)

1. **Algèbres classiques**

| $\mathcal{G}$ | $\mathcal{K}$ |
|---|---|
| $su(p,q)$ | $su(r,s) \oplus su(p-r, q-s) \oplus so(2)$ |
|  | $p-r \geq 0; q-s \geq 0; (r,s) \neq (p-r, q-s)$ |
| $sl(p+q, \mathbb{C})$ | $sl(p, \mathbb{C}) \oplus sl(q, \mathbb{C}) \oplus \mathbb{C}$ |
| $sl(2n, \mathbb{R})$ | $sl(n, \mathbb{C}) \oplus so(2)$ |
| $su^\star(2n)$ | $sl(n, \mathbb{C}) \oplus so(2)$ |
| $so^\star(2n)$ | $su(p, n-p) \oplus so(2)$ $\quad 0 \leq p \leq n$ |
| $so^\star(2n)$ | $so^\star(2n-2) \oplus so(2)$ |
| $so(2n)$ | $su(n) \oplus so(2)$ |
| $so(n)$ | $so(n-2) \oplus so(2)$ |
| $so(2n, \mathbb{C})$ | $sl(n, \mathbb{C}) \oplus \mathbb{C}$ |
| $so(2, n)$ | $so(n-2) \oplus so(2)$ |
| $so(n, \mathbb{C})$ | $so(n-2, \mathbb{C}) \oplus \mathbb{C}$ |
| $so(p, q)$ | $so(p-2, q) \oplus so(2)$ |
| $so(2p, 2q)$ | $su(p, q) \oplus so(2)$ |
| $sp(n, \mathbb{R})$ | $su(p, n-p) \oplus so(2)$ $\quad 0 \leq p \leq n$ |
| $sp(p, q)$ | $su(p, q) \oplus so(2)$ |
| $sp(n, \mathbb{C})$ | $sl(n, \mathbb{C}) \oplus \mathbb{C}$ |

2. **Algèbres exceptionnelles**

| | |
|---|---|
| $e_6^3$ | $so(10) \oplus so(2)$ |
| $e_6^3$ | $so(2,8) \oplus so(2)$ |
| $e_6^3$ | $so^\star \oplus so(2)$ |
| $e_6$ | $so(10) \oplus so(2)$ |
| $e_6^{\mathbb{C}}$ | $so(10, \mathbb{C}) \oplus \mathbb{C}$ |
| $e_6^2$ | $so^\star(10) \oplus so(2)$ |
| $e_6^2$ | $so(4,6) \oplus so(2)$ |
| $e_7^2$ | $so(12) \oplus so(2)$ |
| $e_7$ | $so(12) \oplus so(2)$ |
| $e_7^2$ | $e_6 \oplus so(2)$ |
| $e_7$ | $e_6 \oplus so(2)$ |
| $e_7^{\mathbb{C}}$ | $e_6^{\mathbb{C}} \oplus \mathbb{C}$ |
| $e_7^1$ | $e_6^2 \oplus so(2)$ |
| $e_7^2$ | $e_6^3 \oplus so(2)$ |
| $e_7^2$ | $e_6^2 \oplus so(2)$ |
| $e_7^3$ | $e_6^3 \oplus so(2)$ |

**_TABLEAU B_**   (non-pseudo-kähleriens)

1. **_Algèbres classiques_**

| $\mathcal{G}$ | $\mathcal{K}$ |  |
|---|---|---|
| $sl(n,\mathbb{R})$ | $sl(p,\mathbb{R}) \oplus sl(n-p,\mathbb{R}) \oplus \mathbb{R}$ | $0 < p < n$ |
| $su^\star(2n)$ | $su^\star(2p) \oplus su^\star(2n-2p) \oplus \mathbb{R}$ | $0 < p < n$ |
| $su(n,n)$ | $sl(n,\mathbb{C}) \oplus \mathbb{R}$ | |
| $so^\star(4n)$ | $su^\star(2n) \oplus \mathbb{R}$ | |
| $so(p,q)$ | $so(p-2,q) \oplus \mathbb{R}$ | |
| $so(n,n)$ | $sl(n,\mathbb{R}) \oplus \mathbb{R}$ | |
| $sp(n,\mathbb{R})$ | $sl(n,\mathbb{R}) \oplus \mathbb{R}$ | |
| $sp(n,n)$ | $su^\star(2n) \oplus \mathbb{R}$ | |

2. **_Algèbres exceptionnelles_**

| | |
|---|---|
| $e_6^1$ | $so(5,5) \oplus \mathbb{R}$ |
| $e_6^4$ | $so(1,9) \oplus \mathbb{R}$ |
| $e_7^1$ | $e_6^1 \oplus \mathbb{R}$ |
| $e_7^3$ | $e_6^4 \oplus \mathbb{R}$ |

**_TABLEAU C_**   (T.S.S. compacts, indécomposables, non plats)

1. **_Algèbres classiques_**

| $\mathcal{G}$ | $\mathcal{K}$ |
|---|---|
| $su(p+q)$ | $su(p) \oplus su(q) \oplus so(2)$ |
| $so(n)$ | $so(n-2) \oplus so(2)$ |
| $so(2n)$ | $su(n) \oplus so(2)$ |
| $sp(n)$ | $su(n) \oplus so(2)$ |

2. **_Algèbres exceptionnelles_**

| | |
|---|---|
| $e_6$ | $so(10) \oplus so(2)$ |
| $e_7$ | $e_6 \oplus so(2)$ |

# Chapitre 4

# Résultats relatifs aux T.S.S. ni résolubles, ni semi-simples

Dans ce chapitre, $t = (\mathcal{G}, \sigma, \Omega)$ est un T.S.S. tel que $\mathcal{G}$ n'est ni résoluble, ni semi-simple. Nous adoptons les mêmes notations qu'au chapitre précédent où $\mathcal{S}$ désigne un facteur de Levi $\sigma$-stable de $\mathcal{G}$.

## 4.1 Facteurs semi-simples

**Définition 4.1.1.** *Soit $t = t_{\mathcal{L}} \oplus t_1$ une décomposition où le T.S.S.*
$$t_{\mathcal{L}} = (\mathcal{L}, \sigma|_{\mathcal{L}}, \Omega_{\mathcal{L}})$$
*est semi-simple. On dira que $t_{\mathcal{L}}$ est un facteur semi-simple de $t$. Si $t_1$ n'admet pas de facteur semi-simple, on dira que $t_{\mathcal{L}}$ est un facteur semi-simple maximal.*

**Proposition 4.1.2.** *Un T.S.S. $t = (\mathcal{G}, \sigma, \Omega)$ admet un unique facteur semi-simple maximal $t_{\mathcal{L}}$. En d'autres termes, si $t = t_{\mathcal{L}_1} \oplus t_1 = t_{\mathcal{L}_2} \oplus t_2$ sont deux décompositions où $t_{\mathcal{L}_i}$ est semi-simple maximal ($i = 1, 2$), alors*
$$\mathcal{L}_1 = \mathcal{L}_2$$

PREUVE. Posons, pour $i = 1$ ou $2$,
$$t_{\mathcal{L}_i} = \bigoplus_{\alpha_i} t_{\mathcal{L}_i^{(\alpha_i)}}$$
où
$$t_{\mathcal{L}_i^{(\alpha_i)}} = \left(\mathcal{L}_i^{(\alpha_i)}, \sigma_{\mathcal{L}_i^{(\alpha_i)}}, \Omega_{\mathcal{L}_i^{(\alpha_i)}}\right)$$
est un T.S.S. simple avec la décomposition
$$\mathcal{L}_i^{(\alpha_i)} = \mathcal{K}_i^{(\alpha_i)} \oplus \mathcal{P}_{\mathcal{L}_i^{(\alpha_i)}}$$
Posons encore $t_j = (\mathcal{G}_j, \sigma_j, \Omega_j)$ avec
$$\mathcal{G}_j = \mathcal{K}_j \oplus \mathcal{P}_j \qquad (j = 1, 2)$$
$\forall \alpha_1$, $t_{\mathcal{L}_1^{(\alpha_1)}}$ est non plat; dès lors on a
$$\begin{array}{cccl}\mathcal{P}_{\mathcal{L}_1^{(\alpha_1)}} & \cap & \mathcal{P}_2 & \neq \{0\} \qquad \text{ou bien} \\ \mathcal{P}_{\mathcal{L}_1^{(\alpha_1)}} & \cap & \mathcal{P}_{\mathcal{L}_2^{(\alpha_2)}} & \neq \{0\}\end{array}$$

pour un certain $\alpha_2$.

(a) $\mathcal{P}_{\mathcal{L}_1^{(\alpha_1)}} \cap \mathcal{P}_2 = L \neq \{0\}$. $L$ est lagrangien dans $\mathcal{P}_{\mathcal{L}_1^{(\alpha_1)}}$ et irréductible en tant que $\mathcal{K}_1^{(\alpha_1)}$-module.
Soit alors $L'$ un sous-espace lagrangien de $\mathcal{P}_{\mathcal{L}_1^{(\alpha_1)}}$, $\mathcal{K}$-invariant et en dualité avec $L$. On a alors $[L', L] = \mathcal{K}_1^{(\alpha_1)} = [\pi_2 L', L]$ où $\pi_2 : \mathcal{P} \to \mathcal{P}_2$ est la projection parallèlement à $\mathcal{P}_2^\perp$.
Comme $[\mathcal{K}, L'] = [\mathcal{K}_1^{(\alpha_1)}, L'] = L'$ et $\mathcal{K}_1^{(\alpha_1)} \subset \mathcal{G}_2$, on a $L' = \pi_2(L')$, une contradiction; on a donc

(b) $\mathcal{P}_{\mathcal{L}_1^{(\alpha_1)}} \cap \mathcal{P}_{\mathcal{L}_2^{(\alpha_2)}} \neq \{0\}$
Si $\mathcal{P}_{\mathcal{L}_1^{(\alpha_1)}} \cap \mathcal{P}_{\mathcal{L}_2^{(\alpha_2)}} = L \neq \mathcal{P}_{\mathcal{L}_2^{(\alpha_2)}}$, $L$ est lagrangien. Soit $L'$ comme ci-dessus. De nouveau, $[L', L] = \mathcal{K}_1^{(\alpha_1)} = [\widetilde{\pi_2} L', L]$ où $\widetilde{\pi_2} : \mathcal{P} \to \mathcal{P}_{\mathcal{L}_2^{(\alpha_2)}}$ est la projection parallèlement à $\mathcal{P}_{\mathcal{L}_2^{(\alpha_2)}}^\perp$.
Et donc $L' = \widetilde{\pi_2}(L')$, une contradiction. Dès lors, $\mathcal{P}_{\mathcal{L}_1^{(\alpha_1)}} = \mathcal{P}_{\mathcal{L}_2^{(\alpha_2)}}$; donc $\mathcal{L}_1 \subset \mathcal{L}_2$ et par symétrie $\mathcal{L}_1 = \mathcal{L}_2$.

∎

Dès maintenant nous noterons $\mathcal{L}$ le facteur semi-simple maximal de $t$.

**Lemme 4.1.3.** *Si $[\mathcal{K}_s, \mathcal{P}_\mathcal{R}] = 0$, alors, $\mathcal{S} = \mathcal{L}$; on a la décomposition $t = t_\mathcal{L} \oplus t_\mathcal{R}$ où*

$$t_\mathcal{R} = (\mathcal{R}, \sigma|_\mathcal{R}, \Omega|_{\mathcal{P}_\mathcal{R} \times \mathcal{P}_\mathcal{R}})$$

PREUVE.

(a) $\Omega(\mathcal{P}_s, \mathcal{P}_\mathcal{S}) = \Omega([\mathcal{K}_s, \mathcal{P}_s], \mathcal{P}_\mathcal{R}) = \Omega(\mathcal{P}_s, [\mathcal{K}_s \ \mathcal{P}_\mathcal{R}]) = 0$

(b) $[\mathcal{K}_\mathcal{R}, \mathcal{P}_s] = 0$ car
$$\begin{aligned} \Omega([\mathcal{K}_\mathcal{R} \ \mathcal{P}_s], \mathcal{P}_s \oplus \mathcal{P}_\mathcal{R}) &= \Omega([\mathcal{K}_\mathcal{R} \ \mathcal{P}_s], \mathcal{P}_\mathcal{R}) \\ &= \Omega(\mathcal{P}_s, [\mathcal{K}_\mathcal{R} \ \mathcal{P}_\mathcal{R}]) \\ &= 0 \end{aligned}$$

(c) $[\mathcal{K}_\mathcal{R}, \mathcal{K}_s] = 0$ car
$$\begin{aligned} {[[\mathcal{K}_\mathcal{R} \ \mathcal{K}_s] \ \mathcal{P}_s]} &= [[\mathcal{P}_s \ \mathcal{K}_\mathcal{R}] \ \mathcal{K}_s] + [[\mathcal{K}_s \ \mathcal{P}_s] \ \mathcal{K}_\mathcal{R}] \\ &= 0 + [\mathcal{P}_s \ \mathcal{K}_\mathcal{R}] = 0 \end{aligned}$$

et
$$\begin{aligned} {[[\mathcal{K}_\mathcal{R} \ \mathcal{K}_s] \ \mathcal{P}_\mathcal{R}]} &= [[\mathcal{P}_\mathcal{R} \ \mathcal{K}_\mathcal{R}] \ \mathcal{K}_s] + [[\mathcal{K}_s \ \mathcal{P}_\mathcal{R}] \ \mathcal{K}_\mathcal{R}] \\ &\subset [\mathcal{K}_s \ \mathcal{P}_\mathcal{R}] = 0 \end{aligned}$$

(d) $[\mathcal{P}_\mathcal{R} \ \mathcal{P}_s] = 0$ car
$$[\mathcal{P}_\mathcal{R} \ \mathcal{P}_s] = [\mathcal{P}_\mathcal{R} \ [\mathcal{K}_s \ \mathcal{P}_s]] = 0$$

par Jacobi.

∎

**Remarque 4.1.4.**

Si $[\mathcal{P}_\mathcal{R} \ \mathcal{P}_s] = 0$ ou $[\mathcal{K}_\mathcal{R} \ \mathcal{P}_s] = 0$, alors, comme $[\mathcal{P}_s \ \mathcal{P}_s] = \mathcal{K}_s$, Jacobi implique dans chaque cas : $[\mathcal{K}_s, \mathcal{P}_\mathcal{R}] = 0$.

**Lemme 4.1.5.** *Si $\mathcal{P}_\mathcal{R}$ est symplectique, il existe un facteur de Levi $\mathcal{S}_1$ de $\mathcal{G}$, $\sigma$-stable tel que*

$$\mathcal{P}_\mathcal{R}^\perp = \mathcal{P}_{\mathcal{S}_1} = \mathcal{P} \cap \mathcal{S}_1$$

PREUVE. Soit $\mathcal{P}_1 = \mathcal{P}_\mathcal{R}^\perp$.
$\mathcal{P}_1$ est $\mathcal{K}$-invariant et $\mathcal{P} = \mathcal{P}_1 \oplus \mathcal{P}_\mathcal{R}$. $[\mathcal{K}_\mathcal{R} \ \mathcal{P}_1] = 0$ car

$$\Omega([\mathcal{K}_\mathcal{R} \ \mathcal{P}_1], \ \mathcal{P}_\mathcal{R} \oplus \mathcal{P}_s) = \Omega(\mathcal{P}_1 \ \mathcal{P}_\mathcal{R}) = 0$$

Soit alors $\mathcal{K}_1 = [\mathcal{P}_1 \ \mathcal{P}_1]$ et $\mathcal{S}_1 = \mathcal{K}_1 \oplus \mathcal{P}_1$.
On a $[\mathcal{K}_1 \ \mathcal{P}_\mathcal{R}] = 0$ car

$$\begin{aligned}\Omega([[\mathcal{P}_1 \ \mathcal{P}_1] \ \mathcal{P}_\mathcal{R}], \mathcal{P}_1 \oplus \mathcal{P}_\mathcal{R}) &= \Omega([[\mathcal{P}_1 \ \mathcal{P}_1]\mathcal{P}_\mathcal{R}], \mathcal{P}_\mathcal{R}) \\ &= \Omega([[\mathcal{P}_\mathcal{R} \ \mathcal{P}_1] \ \mathcal{P}_1], \mathcal{P}_\mathcal{R}) \\ &= 0\end{aligned}$$

par invariance de $\mathcal{P}_1$.
Dès lors,

$$[\mathcal{K}_1 \cap \mathcal{K}_\mathcal{R}, \ \mathcal{P}_1 \oplus \mathcal{P}_\mathcal{R}] = [\mathcal{K}_1 \cap \mathcal{K}_\mathcal{R}, \ \mathcal{P}_1] = 0$$

et donc

$$\mathcal{K}_1 \cap \mathcal{K}_\mathcal{R} = 0$$

$\mathcal{P}_1$ étant invariant, $\mathcal{S}_1$ est une sous-algèbre. Comme $\mathcal{G} = \mathcal{S}_1 \oplus \mathcal{R}$, $\mathcal{S}_1$ est un facteur de Levi. ∎

**Proposition 4.1.6.** *Si $\mathcal{P}_\mathcal{R}$ est symplectique, $\mathcal{L}$ est facteur de Levi de $\mathcal{G}$.*

PREUVE. On choisit $\mathcal{S} = \mathcal{S}_1$, on a vu dans la preuve du lemme 4.1.5. que $[\mathcal{K}_1 \ \mathcal{P}_\mathcal{R}] = 0$ c'est-à-dire $[\mathcal{K}_s \ \mathcal{P}_\mathcal{R}] = 0$. On conclut alors par le lemme 4.1.3. ∎

**Lemme 4.1.7.** *Soit $(\mathcal{S}, \sigma)$ un couple symétrique semi-simple, de décomposition canonique $\mathcal{S} = \mathcal{K}_s \oplus \mathcal{P}_s$ Soit $\omega$ une forme antisymétrique, $\mathcal{K}_s$-invariante sur $\mathcal{P}_s$. Alors, il existe un unique $z \in \mathcal{Z}(\mathcal{K}_s)$ avec*

$$\underline{\omega} = \underline{z} \qquad \text{cf. définition} 3.1.7.$$

*En particulier, si $\mathcal{S}$ est simple, $\omega$ est symplectique ou identiquement nulle.*

PREUVE. On a $\delta\underline{\omega} = 0$; donc il existe un unique $X \in \mathcal{S}$ tel que $\underline{\omega} = \delta^\flat X = \underline{X}$ $(H^1(\mathcal{S}) = H^2(\mathcal{S}) = 0)$.
Comme $\underline{\omega}(\mathcal{K}, \mathcal{K}) = \underline{\omega}(\mathcal{K}, \mathcal{P}) = 0$,

$$\beta(X, [\mathcal{K} \ \mathcal{K}] \oplus \mathcal{P}) = 0$$

Donc, $X \in \mathcal{Z}(\mathcal{K})$. On conclut par le théorème 3.1.8. . ∎

**Lemme 4.1.7.** *Soit $(V, \omega)$ un espace vectoriel symplectique. Soient $X$ un sous-espace isotrope de $V$ et $Y$ un supplémentaire de $X$ dans $V$. Alors*

$$\dim \ rad_\omega(Y) \leqslant \dim X$$

PREUVE. Par non-dégénérescence de $\omega$, on a

$$X^\perp \cap rad_\omega(Y) = 0$$

on conclut alors par

$$\dim(X^\perp \oplus rad_\omega(Y)) \leqslant \dim V$$

∎

**Proposition 4.1.8.**
$$2 \leqslant \dim \mathcal{P}_\mathcal{R} \leqslant \dim \mathcal{P} - 2$$

PREUVE. Supposons $\dim \mathcal{P}_\mathcal{R} = 1$; dès lors par le lemme 4.1.7., le rang de $\Omega_s = \Omega|_{\mathcal{P}_s \times \mathcal{P}_s}$, $rg\ \Omega_s$, est tel que
$$rg\ \Omega_s \geqslant \dim \mathcal{P} - 2$$

Si $rg\ \Omega_s = \dim \mathcal{P} - 1$, $\mathcal{P}_s$ est symplectique et $\dim \mathcal{P} - 1$ est pair, une contradiction.
On a donc $\dim \mathcal{P} - 2 = rg\ \Omega_s = \dim \mathcal{P}_s - 1$. Comme $\dim \mathcal{P}_s \geqslant 2$, il existe $z \in \mathcal{Z}(\mathcal{K}_s) \setminus \{0\}$ tel que $\Omega_s = \underline{z}|_{\mathcal{P}_s \times \mathcal{P}_s}$; en particulier $\mathcal{Z}(\mathcal{K}_s) \neq \{0\}$.
Notons $\mathcal{S} = \bigoplus_\alpha \mathcal{S}^\alpha$ la décomposition de $\mathcal{S}$ induite par celle du couple semi-simple $(\mathcal{S}, \sigma|_\mathcal{S})$. On a
$$z = \sum_\alpha z^\alpha \qquad \text{où } z^\alpha \in \mathcal{Z}(\mathcal{K}_{\mathcal{S}^\alpha})$$

Si $z^\alpha \neq 0$, $\mathcal{P}_{\mathcal{S}^\alpha}$ est symplectique.
Si $z^\alpha = 0$, $\mathcal{P}_{\mathcal{S}^\alpha}$ est isotrope et $\Omega(\mathcal{P}_{\mathcal{S}^\alpha}, \mathcal{P}_{\mathcal{S}^\beta}) = 0$ si $\alpha \neq \beta$; dès lors, comme $\dim \mathcal{P}_{\mathcal{S}^\alpha} \geqslant 2\ \forall \alpha$, on a
$$rg\ \Omega_s = \dim \mathcal{P}_s - 2p \qquad \text{où } p \in \mathbb{N}$$

une contradiction.
On conclut par le fait que
$$\dim \mathcal{P}_\mathcal{S} \geqslant 2$$

∎

Nous classifierons dans la suite les T.S.S. avec $\dim \mathcal{P}_\mathcal{R} = 2$.

## 4.2 Résultats relatifs aux cas où $\mathcal{P}_\mathcal{R}$ est isotrope - $\dim \mathcal{P}_\mathcal{R} = 2$

Dans cette section, on suppose $\mathcal{L} = 0$ et $\mathcal{P}_\mathcal{R}$ isotrope.

**Lemme 4.2.1.** $\mathcal{R}$ *est abélien ($\mathcal{P}_\mathcal{R}$ isotrope).*

PREUVE.

- $[\mathcal{K}_\mathcal{R}\ \mathcal{P}_\mathcal{R}] = 0$ car
$$\Omega([\mathcal{K}_\mathcal{R}\ \mathcal{P}_\mathcal{R}],\ \mathcal{P}_s \oplus \mathcal{P}_\mathcal{R}) = \Omega(\mathcal{P}_\mathcal{R},\ [\mathcal{K}_\mathcal{R}\ \mathcal{P}_s]) = 0$$

- $[\mathcal{K}_\mathcal{R}\ \mathcal{K}_\mathcal{R}] = 0$ car
$$\begin{aligned}
[[\mathcal{K}_\mathcal{R}\ \mathcal{K}_\mathcal{R}]\ \mathcal{P}_\mathcal{R} \oplus \mathcal{P}_s] &= [[\mathcal{K}_\mathcal{R}\ \mathcal{K}_\mathcal{R}]\ \mathcal{P}_s] \\
&= [[\mathcal{P}_s\ \mathcal{K}_\mathcal{R}]\ \mathcal{K}_\mathcal{R}] \\
&\subset [\mathcal{K}_\mathcal{R}\ \mathcal{P}_\mathcal{R}] \\
&= 0
\end{aligned}$$

- $[\mathcal{P}_\mathcal{R}\ \mathcal{P}_\mathcal{R}] = 0$ car
$$\begin{aligned}
[[\mathcal{P}_\mathcal{R}\ \mathcal{P}_\mathcal{R}]\ \mathcal{P}_\mathcal{R} \oplus \mathcal{P}_s] &= [[\mathcal{P}_\mathcal{R}\ \mathcal{P}_\mathcal{R}\ \mathcal{P}_s]] \\
&= [[\mathcal{P}_s\ \mathcal{P}_\mathcal{R}]\ \mathcal{P}_\mathcal{R}] \\
&\subset [\mathcal{K}_\mathcal{R}\ \mathcal{P}_\mathcal{R}] \\
&= 0
\end{aligned}$$

∎

**Lemme 4.2.2.** *Si $\mathcal{Z}(\mathcal{K}_s)$ contient un élément $z$ tel que $(\mathrm{ad}(z)|_{\mathcal{P}_s})^2 = \varepsilon\ I|_{\mathcal{P}_s}$ ($\varepsilon = \pm 1$); alors, en posant $\mathcal{J} = \mathrm{ad}(z)|_\mathcal{P}$, on a*

(i) $\mathcal{J}^2 = \varepsilon I$

(ii) Si $\varepsilon = -1$, $\mathcal{J}$ est symplectique.

PREUVE. Pour tous $p, p' \in \mathcal{P}$, on a, en notant $p = p_r + p_s$ et $p = p'_r + p'_s$ les décompositions relativement à $\mathcal{P} = \mathcal{P}_\mathcal{R} \oplus \mathcal{P}_s$ :
$$\begin{aligned}
\Omega(\mathcal{J}^2 p, p') &= \Omega(\varepsilon p_s + \mathcal{J}^2 p_r, p'_r + p'_s) \\
&= \varepsilon \Omega(p_s, p') + \Omega(\mathcal{J}^2 p_r, p'_s) \\
&= \varepsilon \Omega(p_s, p') + \varepsilon \Omega(p_r, p'_s) \\
&= \varepsilon \Omega(p_s, p') + \varepsilon \Omega(p_r, p') \\
&= \varepsilon \Omega(p, p')
\end{aligned}$$

Maintenant, $\Omega(\mathcal{J}p, \mathcal{J}p') = -\Omega(\mathcal{J}^2 p, p') = -\varepsilon \Omega(p, p')$. ∎

Maintenant, notons
$$\mathcal{S} = \bigoplus_{\alpha=1}^{A+I} \mathcal{S}^\alpha = \bigoplus_{a=1}^{A} \mathcal{S}^a \oplus \bigoplus_{i=A+1}^{A+I} \mathcal{S}^i$$

où
$$\mathcal{S}^\alpha = \mathcal{K}_s^\alpha \oplus \mathcal{P}_s^\alpha$$

est la décomposition induite par celle du couple symétrique semi-simple $(\mathcal{S}, \sigma|_\mathcal{S})$ où l'on suppose

$$\mathcal{Z}(\mathcal{K}_\mathcal{S}^a) = \{0\} \qquad\qquad \forall a = 1, \ldots, A$$

et

$$\mathcal{Z}(\mathcal{K}_\mathcal{S}^i) \neq \{0\} \qquad\qquad \forall i = A+1, \ldots, A+I$$

avec la convention
$$I = 0 \quad \Longleftrightarrow \quad \mathcal{Z}(\mathcal{K}_s) = \{0\}$$

Remarquons que $\mathcal{S}^i$ est simple pour tout $i \geqslant A + I$.

**Lemme 4.2.3.** $\Omega(\mathcal{P}_s^\alpha, \mathcal{P}_\mathcal{R}) = 0 \iff [\mathcal{K}_s^\alpha, \mathcal{P}_\mathcal{R}] = 0$.

PREUVE.
$$\begin{aligned}
& \Omega(\mathcal{P}_s^\alpha,\ \mathcal{P}_\mathcal{R}) = 0 \\
\iff & \Omega([\mathcal{K}_s^\alpha\ \mathcal{P}_s^\alpha],\ \mathcal{P}_\mathcal{R}) = 0 \\
\iff & \Omega([\mathcal{K}_s^\alpha\ \mathcal{P}_s],\ \mathcal{P}_\mathcal{R}) = 0 \\
\iff & \Omega(\mathcal{P}_s \oplus \mathcal{P}_\mathcal{R},\ [\mathcal{K}_s^\alpha\ \mathcal{P}_\mathcal{R}]) = 0 \qquad (\mathcal{P}_\mathcal{R} \text{ est isotrope}) \\
\iff & [\mathcal{K}_s^\alpha\ \mathcal{P}_\mathcal{R}] = 0
\end{aligned}$$

∎

**Lemme 4.2.4.** $\Omega(\mathcal{P}_s^\alpha,\ \mathcal{P}_s^\beta) = 0$ si $\alpha \neq \beta$.

PREUVE.
$$\begin{aligned}
\Omega(\mathcal{P}_s^\alpha,\ \mathcal{P}_s^\beta) &= \Omega([\mathcal{K}_s^\alpha\ \mathcal{P}_s^\alpha],\ \mathcal{P}_s^\beta) \\
&= \Omega(\mathcal{P}_s^\alpha,\ [\mathcal{K}_s^\alpha\ \mathcal{P}_s^\beta])
\end{aligned}$$

∎

**Lemme 4.2.5.** *L'application*
$$\bigoplus_{a=1}^A \mathcal{P}_s^a \to \mathcal{P}_\mathcal{R}^\star : s \mapsto \Omega(s, \cdot)$$

*est injective.*

PREUVE. Soit $s \in \bigoplus_a \mathcal{P}_s^a$.

Si $\Omega(s, \mathcal{P}_\mathcal{R}) = 0$, alors
$$\begin{aligned}
\Omega(s, \mathcal{P}) &= \Omega\left(s, \bigoplus_{a=1}^A \mathcal{P}_s^a\right) \\
&= \xi\left[s, \bigoplus_{a=1}^A \mathcal{P}_s^a\right]
\end{aligned}$$

où $\xi \in \mathcal{Z}\left(\bigoplus_{a=1}^A \mathcal{K}_s^a\right)^\star$

Donc $\Omega(s, \mathcal{P}) = 0$ et $s = 0$.

∎

Posons
$$\begin{aligned}
N_s^\alpha &= \mathcal{P}_\mathcal{R}^\perp \cap \mathcal{P}_s^\alpha \\
\text{et} \quad N_\mathcal{R}^\alpha &= \mathcal{P}_s^{\alpha\perp} \cap \mathcal{P}_\mathcal{R}
\end{aligned}$$

**Lemme 4.2.6.** $N_s^\alpha$ et $N_\mathcal{R}^\alpha$ sont $\mathcal{K}$-invariants $\forall \alpha$.

PREUVE. $N_s^\alpha$ (resp. $N_\mathcal{R}^\alpha$) est clairement $\mathcal{K}_s$-invariant (resp. $\mathcal{K}_\mathcal{R}$-invariant).
On a
$$\begin{aligned}
\Omega([\mathcal{K}_\mathcal{R},\ N_s^\alpha],\ \mathcal{P}_\mathcal{R} \oplus \mathcal{P}_s) &= \Omega(N_s^\alpha,\ [\mathcal{K}_\mathcal{R}\ \mathcal{P}_s] \subset \mathcal{P}_\mathcal{R}) \\
&= 0
\end{aligned}$$

$$\implies [\mathcal{K}_\mathcal{R},\ N_s^\alpha] = 0$$

et
$$\begin{aligned}
\Omega([\mathcal{K}_s,\ N_\mathcal{R}^\alpha],\ \mathcal{P}_s^\alpha) &= \Omega(N_\mathcal{R}^\alpha,\ [\mathcal{K}_s^\alpha\ \mathcal{P}_s^\alpha] = \mathcal{P}_s^\alpha) \\
&= 0
\end{aligned}$$

∎

Choisissons $L_s^\alpha$ un sous-espace $\mathcal{K}_s$-invariant dans $\mathcal{P}_s$ tel que $L_s^\alpha \oplus N_s^\alpha = \mathcal{P}_s^\alpha$ et $L_\mathcal{R}^\alpha$ un sous-espace $\mathcal{K}$-invariant supplémentaire de $N_\mathcal{R}^\alpha$ dans $\mathcal{P}_\mathcal{R}$ ($\mathcal{K}_s$ agit complètement ré- ductiblement sur $\mathcal{P}$.

**Lemme 4.2.7.** $L_\mathcal{R}^\alpha \neq \{0\} \quad \forall \alpha$.

PREUVE. Si $L_\mathcal{R}^\beta = \{0\}$, alors $N_\mathcal{R}^\beta = \mathcal{P}_\mathcal{R}$ et $\Omega(\mathcal{P}_s^\beta, \mathcal{P}_\mathcal{R}) = 0$.
Donc,
$$\Omega([\mathcal{K}_s^\beta \ \mathcal{P}_\mathcal{R}], \ \mathcal{P}) = \Omega([\mathcal{K}_s^\beta, \ \mathcal{P}_\mathcal{R}], \ \mathcal{P}_s^\beta) = 0$$
($\mathcal{P}_\mathcal{R}$ est isotrope), c.à.d. $[\mathcal{K}_s^\beta \ \mathcal{P}_\mathcal{R}] = 0$.
De plus,
$$\Omega([\mathcal{K}_\mathcal{R} \ \mathcal{P}_s^\beta], \ \mathcal{P}_s) = \Omega(\mathcal{P}_s^\beta, \ [\mathcal{K}_\mathcal{R} \ \mathcal{P}_s]) = 0$$
c.à.d. $[\mathcal{K}_\mathcal{R} \ \mathcal{P}_s^\beta] = 0$.
Un raisonnement analogue à celui tenu dans la preuve du lemme 4.1.3. montre alors que $\mathcal{S}^\beta$ est un idéal de $\mathcal{G}$. Ceci contredit $\mathcal{L} = 0$. ∎

Remarquons qu'on a les injections $L_\mathcal{R}^\alpha \hookrightarrow \mathcal{P}_s^{\alpha\star}$, $L_s^\alpha \hookrightarrow \mathcal{P}_\mathcal{R}^\star$ et que $\forall a$, $N_s^a = 0$ et $L_s^a = \mathcal{P}_s^a$.

**Lemme 4.2.8.** $\mathcal{P}_\mathcal{R}$ est un $\mathcal{K}_s^\alpha$-module fidèle $\forall \alpha$.

PREUVE. Soit $k \in \mathcal{K}_s^\alpha$ tel que $[k \ \mathcal{P}_\mathcal{R}] = 0$; alors

- si $N_s^\alpha \neq 0$, $\Omega([k \ L_s^\alpha], \ \mathcal{P}_\mathcal{R}) = 0$
  donc $[k \ L_s^\alpha] = 0$ et $k = 0$
  (voir chapitre 3.).

- si $N_s^\alpha = 0, L_s^\alpha = \mathcal{P}_s^\alpha$ et $\Omega(\mathcal{P}_s^\alpha, \ [k \ \mathcal{P}_\mathcal{R}]) = 0$ donc $[k \ \mathcal{P}_s^\alpha] = 0$ et $k = 0$.

∎

**Lemme 4.2.9.** Les sous-$\mathcal{K}$-modules $\{L_\mathcal{R}^\alpha\}$ sont en somme directe dans $\mathcal{P}_\mathcal{R}$ et $[\mathcal{K}_s^\beta \ L_\mathcal{R}^\alpha] = 0$ si $\alpha \neq \beta$.

PREUVE. Posons $L^{\alpha_1,...,\alpha_p} = L_\mathcal{R}^{\alpha_1} \cap [L_\mathcal{R}^{\alpha_2} + \ldots + L_\mathcal{R}^{\alpha_p}]$
On a $\Omega([\mathcal{K}_s^\beta \ L_\mathcal{R}^\alpha], \ \mathcal{P}_s^\alpha) = 0$ donc
$$[\mathcal{K}_s^\beta \ L_\mathcal{R}^\alpha] = 0$$
Dès lors, $[\mathcal{K}_s^{\alpha_1}, \ L^{\alpha_1,...,\alpha_p}] = 0$ donc
$$\Omega([\mathcal{K}_s^{\alpha_1}, \ L^{\alpha_1,...,\alpha_p}], \ \mathcal{P}_s^{\alpha_1}) = 0 = \Omega(L^{\alpha_1,...,\alpha_p}, \ \mathcal{P}_s^{\alpha_1})$$
et $L^{\alpha_1,...,\alpha_p} = 0$. ∎

**Lemme 4.2.10.** $\Omega(\mathcal{P}_s^\alpha, \ L_\mathcal{R}^\beta) = 0$ si $\alpha \neq \beta$.

PREUVE.
$$\begin{array}{rcl}\Omega(\mathcal{P}_s^\alpha, \ L_\mathcal{R}^\beta) & = & \Omega([\mathcal{K}_s^\alpha \ \mathcal{P}_s^\alpha], \ L_\mathcal{R}^\beta) \\ & = & \Omega(\mathcal{P}_s^\alpha, \ [\mathcal{K}_s^\alpha \ L_\mathcal{R}^\beta]) \\ & = & 0 \quad \text{si } \alpha \neq \beta\end{array}$$
∎

On pose $L_\mathcal{R} = \displaystyle\bigoplus_{a=1}^{A} L_\mathcal{R}^a \oplus \bigoplus_{i=A+1}^{A+I} L_\mathcal{R}^i$.

**Proposition 4.2.11.**
  (i) $\dim \mathcal{P}_\mathcal{R} \geqslant 2A + I$
  (ii) $I \neq \dim \mathcal{P}_\mathcal{R}$.

Preuve.

(i) On a $[\mathcal{K}_s^a \ L_\mathcal{R}^a] \neq \{0\}$ car $\Omega([\mathcal{K}_s^a \ L_\mathcal{R}^a], \mathcal{P}_s^a) = \Omega(L_\mathcal{R}^a, \mathcal{P}_s^a)$ et $L_\mathcal{R}^a \neq \{0\}$.
$\mathcal{K}_s^a$ étant simple, $\dim L_\mathcal{R}^a \geqslant 2$.
Comme $L_\mathcal{R}^i \neq \{0\}$, on a $\dim L_\mathcal{R}^i \geqslant 1$ donc
$$\dim \mathcal{P}_\mathcal{R} \geqslant \dim L_\mathcal{R} \geqslant 2A + I$$

(ii) Si $I = \dim \mathcal{P}_\mathcal{R}$, $A = 0$.
On a $\dim L_\mathcal{R}^i = 1 \ \forall i$ et
$$\mathcal{P}_\mathcal{R} = L_\mathcal{R}$$
Dès lors,
$$[\mathcal{K}_s', \mathcal{P}_\mathcal{R}] = \bigoplus_i [\mathcal{K}_s^{i'}, L_\mathcal{R}^i] = 0$$
car $\mathcal{K}_s^{i'}$ est semi-simple.
Comme $\mathcal{P}_\mathcal{R}$ est $\mathcal{K}_s^i$-fidèle, on a $\dim_\mathbb{R} \mathcal{K}_s^i = 1 \ \forall i$ et donc
$$\dim \mathcal{P}_s^i = 2 \qquad\qquad \forall i$$
Mais alors, $\Omega(\mathcal{P}_s^i \oplus L_\mathcal{R}^i, \mathcal{P}_s^j \oplus L_\mathcal{R}^j) = 0 \ \forall i \neq j$; donc $\mathcal{P}_s^i \oplus L_\mathcal{R}^i$ est symplectique $\forall i$, ce qui est impossible car $\dim(\mathcal{P}_s^i \oplus L_\mathcal{R}^i) = 3$.

∎

**Remarque**. Si $\dim \mathcal{P}_\mathcal{R} = 2$, on a $(A = 1, I = 0)$ ou $(A = 0, I = 1)$; le lemme 4.2.5. montre que le premier cas $(A = 1)$ est impossible.

**Proposition 4.2.12.** *Si $A = 0$ et $I = 1$, alors $\mathcal{P}_\mathcal{R}$ et $\mathcal{P}_s$ sont deux $\mathcal{K}_s$-modules en dualité; en particulier, on a $\dim \mathcal{P} = 2 \dim \mathcal{P}_s = 2 \dim \mathcal{P}_\mathcal{R}$.*

D'où le

**Théorème 4.2.13.** *Si $t$ est un T.S.S. non résoluble, non semi-simple et sans facteur semi-simple alors $\dim t = 4$ si et seulement si $\dim \mathcal{P}_\mathcal{R} = 2$.*

Preuve. Comme $A = 0$ et $I = 1$, on a $\mathcal{P}_\mathcal{R} = N_\mathcal{R} \oplus L_\mathcal{R}$ et $\mathcal{P}_s = N_s \oplus L_s$ où $N_s = \mathcal{P}_s \cap \mathcal{P}_\mathcal{R}^\perp$ et $N_\mathcal{R} = \mathcal{P}_s^\perp \cap \mathcal{P}$; mais $\Omega(N_\mathcal{R}, \mathcal{P}_\mathcal{R} \oplus \mathcal{P}_s) = 0$ donc $N_\mathcal{R} = 0$ et $\mathcal{P}_\mathcal{R} = L_\mathcal{R}$. Par le lemme 4.1.7. $\mathcal{P}_s$ est symplectique ou isotrope; $N_s$ lui ne peut être égal à $\mathcal{P}_s$ (par non dégénérescence de $\Omega$). Quant à $L_s$, il est en $\Omega$-dualité avec $\mathcal{P}_\mathcal{R}$; en effet, si $r \in \mathcal{P}_\mathcal{R}$ est tel que $\Omega(r, L_s) = 0$, on a, $\Omega(r, \mathcal{P}_s) = 0$ et $r \in N_\mathcal{R} = 0$.
Supposons $N_s \neq 0$. Dans ce cas, $\mathcal{P}_s$ est symplectique et, $N_s$ et $L_s$ sont des $\mathcal{K}_s$-modules irréductibles en $\Omega$-dualité totalement isotropes ($\mathcal{S}$ est simple car $\mathcal{Z}(\mathcal{K}_s) \neq 0$). Soit $\{\ell_i\}$ une base de $L_s$; soient $\{p_i\}$ et $\{n_i\}$ les bases duales de $\mathcal{P}_\mathcal{R}$ et $N_s$ respectivement. Dans la base $\{n_i\} \cup \{\ell_i\} \cup \{p_i\}$, $\Omega$ prend la forme

$$\Omega = \begin{pmatrix} 0 & I_n & 0 \\ -I_n & 0 & -I_n \\ 0 & I_n & 0 \end{pmatrix} \qquad \text{où } n = \dim N_s$$

Cette matrice étant singulière, on a $N_s = 0$. Maintenant soit $\{s_i\}$ une base de $\mathcal{P}_s$ et $\{r_i\}$ la base duale de $\mathcal{P}_\mathcal{R}$. Notons alors $\forall k \in \mathcal{K}_s$,

$$[k, s_i] = S(k)^j{}_i s_j \quad \text{et} \quad [k, r_i] = R(k)^j{}_i r_j$$

On a alors
$$\Omega([k\ s_i],\ r_j) = \Omega([k\ r_j],\ s_i)$$

c'est-à-dire
$$S(k)^a{}_i\ \Omega(s_a, r_j) = R(k)^b{}_j \Omega(r_b, s_i)$$

c'est-à-dire
$$S(k)^a{}_i \delta_{aj} = R(k)^b{}_j\ \delta_{bi}$$

c'est-à-dire
$$-S(k)^j{}_i = R(k)^i{}_j$$

ou bien
$$R(k) = -{}^\tau S(k)$$

■

## 4.3 $\dim t = 4$

Dans cette section, on suppose que $t = (\mathcal{G}, \sigma, \Omega)$ est un T.S.S. de dimension 4, non résoluble, non semi-simple et sans facteur semi-simple. On a alors $\dim \mathcal{S} = 3$, $\dim \mathcal{P}_\mathcal{R} = 2$, $\mathcal{P}_\mathcal{R}$ est isotrope et en dualité avec $\mathcal{P}_s$ en tant que $\mathcal{K}_s$-module.

La table de $\mathcal{G}$ est alors du type :

$$\begin{array}{rcl} [k\ s_i] & = & \kappa^j{}_i\ s_j \quad (i, j = 0, 1) \\ [s_0\ s_1] & = & k \\ [k\ r_i] & = & -\kappa^i{}_j\ s_j \\ [r_i\ s_j] & = & x_{ij} \\ [x_{ij}\ k] & = & A^{k\ell}_{ij} x_{k\ell} \\ [x_{ij}\ s_k] & = & X^{\ell}_{ij,k} r_\ell \end{array}$$

où
$$\mathcal{K}_s = \mathbb{R}k \qquad \mathcal{P}_s = \rangle s_i \langle \qquad \mathcal{P}_\mathcal{R} = \rangle r_i \langle \quad \text{et} \qquad \mathcal{K}_\mathcal{R} = \rangle x_{ij} \langle$$

On remarque alors que les relations
$$[[k\ s_i]\ r_j] + [[r_j\ k]\ s_i] + [[s_i\ r_j]\ k] = 0$$

et
$$[[s_0\ s_1]\ x_{k\ell}] + [[x_{k\ell}\ s_0]\ s_1] + [[s_1\ x_{kl}]\ s_0] = 0$$

fixent les $A^{k\ \ell}_{ij}$ et les $X^\ell_{ij,k}$ en fonction des $\kappa^i{}_j$. Dès lors, 3 cas se présentent :

(a) $\mathcal{S} = sl(2, \mathbb{R})$ et $k$ est compact,

(b) $\mathcal{S} = su(2)$,

(c) $\mathcal{S} = sl(2, \mathbb{R})$ et $k$ n'est pas compact.

Ces trois cas sont réalisés et on a les trois tables :

$$\begin{array}{rcl} [k\ s_i] & = & (-1)^{i+\varepsilon+1} s_{(i+1)\ \mathrm{mod}(2)} \\ [s_0\ s_1] & = & k \\ [k\ r_i] & = & (-1)^{i+\varepsilon} r_{(i+1)\ \mathrm{mod}(2)} \\ [x\ s_i] & = & (-1)^\varepsilon r_i \\ [r_i\ s_j] & = & \delta_{ij} x \end{array} \qquad \text{où } i = 0, 1;\ \varepsilon = \pm 1$$

et

$$\begin{array}{rcl}
[k\ s_i] & = & (-1)^{i+1}2s_i \\
[s_0\ s_1] & = & k \\
[k\ r_i] & = & (-1)^i 2r_i \\
[x\ s_i] & = & 2\ r_{(i+1)\ \mathrm{mod}(2)} \\
[r_i\ s_j] & = & -\delta_{ij}x
\end{array} \qquad \text{où } i = 0,1$$

Les deux premières tables correspondent aux cas $(a)$ et $(b)$ avec respectivement $\varepsilon = -1$ et $\varepsilon = 1$; la troisième au cas $(c)$.
Par un automorphisme de $(\mathcal{G}, \sigma)$ du type

$$\begin{array}{rcl rcl}
r_i & \leftarrow & \alpha\ r_i & \quad s_i & \leftarrow & s_i \\
x & \leftarrow & \alpha x & \quad k & \leftarrow & k
\end{array}$$

on ramène toute forme symplectique $\mathcal{K}$-invariante sur $\mathcal{P}$ à une des formes suivantes :

$$\Omega_a = \begin{pmatrix} a\mathcal{J} & I_2 \\ -I_2 & 0 \end{pmatrix}$$

(exprimée dans la base $\{s_0, s_1, r_0, r_1\}$ de $\mathcal{P}$), où

$$\mathcal{J} = \begin{pmatrix} 0 & 1 \\ -1 & 0 \end{pmatrix}$$

et $a \in \mathbb{R}$.
Par un calcul brutal, on vérifie alors que $(\mathcal{G}, \sigma, \Omega_a)$ n'est pas isomorphe à $(\mathcal{G}, \sigma, \Omega_b)$ si $a \neq b$.
**Remarque**. Les E.S.S. associées aux triples obtenus dans les cas $(a)$, $(b)$ et $(c)$ sont les fibrés cotangents à respectivement

$$SL_2(\mathbb{R})\big/SO(2), \qquad SU(2)\big/SO(2) \qquad \text{et} \qquad SL_2(\mathbb{R})\big/\mathbb{R}$$

(cf. [3, pg. 102]).

## 4.4 L'unique facteur compact

**Théorème 4.4.1.** *Soit $(M, \omega, s)$ un E.S.S. simplement connexe. Soient $o$ un point de $M$ et $K$ l'holonomie linéaire en $o$ relativement à la connexion affine canonique $\nabla$.*
*Soit $M_c$ une sous-variété de $M$ maximale pour les propriétés suivantes :*

(a) $M_c$ passe par $o$

(b) $M_c$ est connexe

(c) $M_c$ est totalement géodésique

(d) $T_o(M_c)$ est invariant par $K$

(e) $M_c$ est compacte

*Alors,*

(i) $M_c$ *est simplement connexe et symplectique.*

(ii) $M_c$ *est la sous-variété de $M$ sous-jacente à un sous-espace symétrique compact $(M_c, \omega_c, s_c)$ facteur direct de $(M, \omega, s)$.*

(iii) $(M_c, \omega_c, s_c)$ *admet un groupe des transvections semi-simple compact.*

(*iv*) $M_c$ est l'unique sous-variété de $M$ maximale pour les propriétés (a)–(e)

Ceci résultera d'une succession de lemmes.

Notons $G$ le groupe des transvections de $(M, \omega, s)$, $(\mathcal{G}, \sigma, \Omega)$ le T.S.S. associé et $\widetilde{\sigma}$ l'automorphisme canonique de $G$. Notons $(\pi_{\star_e}|_{\mathcal{P}})^{-1}(T_o(M_c)) = \mathcal{P}_c$.
$T_o(M_c)$ étant $\mathcal{K}$-invariant, il en est de même pour $\mathcal{P}_c$ et dès lors $\boldsymbol{h} = [\mathcal{P}_c \, \mathcal{P}_c] \oplus \mathcal{P}_c$ est une sous-algèbre ($\sigma$-stable) de $\mathcal{G}$; notons $H$ le sous-groupe connexe de $G$ d'algèbre $\boldsymbol{h}$.

**Lemme 4.4.2.** *$H$ est compact.*

PREUVE. $M_c$ étant compacte, son stabilisateur $L$ dans $G$ est fermé et la restriction de l'action de $L$ sur $M_c$ définit un homomorphisme de groupes de Lie :

$$L \xrightarrow{r} Aut(M_c)$$

$H$ est engendré par $Exp_G(\mathcal{P}_c)$ c'est-à-dire par

$$\{s_{Exp_o(X)} \circ s_o \mid X \in T_o(M_c)\}$$

$M_c$ étant totalement géodésique et stabilisé par les symétries centrées en ses points, $r$ induit un homomorphisme injectif de groupes de Lie :

$$H \longrightarrow Aut(M_c)$$

Cet homomorphisme est surjectif sur le groupe $G_c$ des transvections de $M_c$.
Par [Koh], $G_c$ est compact. ∎

Notons $\mathcal{R}$ le radical de $\mathcal{G}$.

**Lemme 4.4.3.** *Tout sous-groupe compact $C$ stable par $\widetilde{\sigma}$ stabilise une sous-algèbre de Levi $\mathcal{S}$ $\sigma$-stable dans $\mathcal{G}$. De plus, la partie semi-simple de l'algèbre $\mathcal{C}$ de $C$ est contenue dans $\mathcal{S}$.*

PREUVE. Soit le groupe fini d'automorphismes de $\mathcal{G}$ : $\mathbb{Z}_2 = \{id, \sigma\}$. $C$ étant $\widetilde{\sigma}$-stable, on peut construire le produit semi-direct $C_\sigma = \mathbb{Z}_2 \times C$ où

$$(\varepsilon, u) \cdot (\varepsilon', u') = (\varepsilon \varepsilon', \varepsilon'(u) \, u')$$

pour tous $\varepsilon, \varepsilon'$ dans $\mathbb{Z}_2$ et $u, u'$ dans $C$ (si $\varepsilon = \sigma$, on définit $\varepsilon(u) = \widetilde{\sigma}(u)$, $u$ dans $C$).
$C_\sigma$ est alors un groupe compact d'automorphisme de $\mathcal{G}$ (($\varepsilon, u) \cdot X = \varepsilon(Ad(u)(X))$) et le théorème de Taft nous livre un tel facteur de Levi $\mathcal{S}$. Soit $\mathcal{D} = [\mathcal{C}, \mathcal{C}]$. Notons $\mathcal{D}_{(\mathcal{R})}$ la projection de $\mathcal{D}$ sur $\mathcal{R}$ parallèlement à $\mathcal{S}$ ($\mathcal{G} = \mathcal{S} \oplus \mathcal{R}$) et $\mathcal{D}_{(\mathcal{S})}$ la projection sur $\mathcal{S}$ parallèlement à $\mathcal{R}$. Soit $d \in \mathcal{D}$ alors, relativement à la décomposition $\mathcal{G} = \mathcal{S} \oplus \mathcal{R}$ on a $d = d_\mathcal{S} + d_\mathcal{R}$. Comme $C$ stabilise $\mathcal{S}$, on a $[\mathcal{D}, \mathcal{S}] \subset \mathcal{S}$ et donc

$$[d, \mathcal{S}] = [d_\mathcal{S}, \mathcal{S}] + [d_\mathcal{R}, \mathcal{S}] \subset \mathcal{S}$$

comme $[d_\mathcal{R}, \mathcal{S}] \subset \mathcal{R}$ on a $[\mathcal{D}_{(\mathcal{R})}, \mathcal{S}] = 0$; en particulier, $[\mathcal{D}_{(\mathcal{R})}, \mathcal{D}_{(\mathcal{S})}] = 0$.
Dès lors, $\mathcal{D}_{(\mathcal{R})}$ et $\mathcal{D}_{(\mathcal{S})}$ sont des sous-algèbres et pour tout $N \in \mathbb{N}$, on a $\mathcal{D}^{(N)} \subset \mathcal{D}_{(\mathcal{R})}^{(N)} \oplus \mathcal{D}_{(\mathcal{S})}^{(N)}$. Comme $\mathcal{D}_{(\mathcal{R})}$ est résoluble, il existe $N$ avec $\mathcal{D}_{(\mathcal{R})}^{(N)} = \{0\}$; mais $\mathcal{D}$ est semi-simple donc pour tout $N$ : $\mathcal{D}^{(N)} = \mathcal{D}$ dès lors $\mathcal{D} \subset \mathcal{S}$. ∎

Notons $\mathcal{Z}$ le centre de $\boldsymbol{h}$. Soit $U$ un sous-groupe compact maximal de $G$, stable par $\widetilde{\sigma}$ ($\widetilde{\sigma}(g) = s_o \, g \, s_o \; \forall g \in G$) et tel que $M$ se fibre sur l'espace symétrique compact $M_U = U/U \cap K$, les fibres étant homéomorphes à des espaces euclidiens ([Koh]).

**Lemme 4.4.4.** *H est semi-simple.*

PREUVE. Soit $\mathcal{U}$ l'algèbre de $U$; on a
$$\mathcal{U} = \mathcal{U}_\mathcal{K} \oplus \mathcal{U}_\mathcal{P}$$
où $\mathcal{U}_\mathcal{K} = \mathcal{U} \cap \mathcal{K}$ et $\mathcal{U}_\mathcal{P} = \mathcal{U} \cap \mathcal{P}$.
Comme $M$ est simplement connexe, $M_U$ l'est aussi. Dès lors, $[\mathcal{U}_\mathcal{P}\,\mathcal{U}_\mathcal{P}] \oplus \mathcal{U}_\mathcal{P} = \mathcal{A}$ est un idéal semi-simple de $\mathcal{U}$. Choisissons un supplémentaire $\mathcal{B}$ de $\mathcal{A}$ dans $\mathcal{U}$ avec $\mathcal{B} \subset \mathcal{U}_\mathcal{K}$.
Soit $\mathcal{S}_1$ une sous-algèbre de Levi de $\mathcal{G}$, stable par $\sigma$ et par $U$. Alors $\mathcal{A} \subset \mathcal{S}_1$.
Par ailleurs, soit $\mathcal{S}$ une sous-algèbre de levi de $\mathcal{G}$, stable par $\sigma$ et $H$. Notons $\mathcal{Z}$ le centre de $\boldsymbol{h}$. Alors $\mathcal{Z} \subset \mathcal{P}$ et $[\mathcal{Z}, \mathcal{K}_s] \subset \mathcal{Z} \cap \mathcal{S}$ ($\mathcal{P}_c$ est $\mathcal{K}$-stable). $\mathcal{Z} \cap \mathcal{K}_s$ est un sous-espace $\mathcal{K}_s$-invariant compact et abélien dans $\mathcal{P}_s$, il est donc nul (utiliser une décomposition en "indécomposables" du couple $(\mathcal{S}, \sigma\big|_\mathcal{S})$ et la proposition 3.1.5. Maintenant, tout élément $z$ de $\mathcal{Z}$ s'écrit $z = z_\mathcal{R} + z_\mathcal{S}$ avec $z_\mathcal{R} \in \mathcal{R}$ et $z_s \in \mathcal{S}$.
Comme $[\mathcal{Z}, \mathcal{S}] \subset \mathcal{S}$, $[z_\mathcal{R}, \mathcal{S}] = 0$ mais $[\mathcal{K}_s, z] = [\mathcal{K}_s,\ z_s] = 0$ donc $z_s = 0$ (car $z_s \in \mathcal{P}_s$ sur lequel $\mathcal{K}_s$ a une action effective). Dès lors, on a $[\mathcal{Z}\,\mathcal{S}] = 0$ et $\mathcal{Z} \subset \mathcal{R}$.
Soit $\varphi$ un automorphisme de $(\mathcal{G}, \sigma)$ avec
$$\varphi(\mathcal{S}_1) = \mathcal{S}$$
On a $\varphi(\mathcal{B}) \subset \mathcal{K}$ et $\varphi(\mathcal{A}) \subset \mathcal{S}$ donc $[\mathcal{Z}, \varphi(\mathcal{A})] = 0$ et $[\varphi(\mathcal{B}), \mathcal{Z}] \subset \mathcal{Z}$. Dès lors, $\mathcal{Z} + \varphi(\mathcal{U})$ est une sous-algèbre compacte dans $\mathcal{G}$; par maximalité
$$\mathcal{Z} + \varphi(\mathcal{U}) = \varphi(\mathcal{U})$$
donc $\mathcal{Z} \subset \varphi(\mathcal{U})$.
Mais $\mathcal{Z} \subset \mathcal{P}$ donc $\mathcal{Z} \subset \varphi(\mathcal{U}_\mathcal{P})$; comme $\varphi(\mathcal{U}_\mathcal{P}) \subset \mathcal{S}$ et $\mathcal{Z} \subset \mathcal{R}$, $\mathcal{Z} = 0$. ∎

Soit $\mathcal{S}$ une sous-algèbre de Levi $\sigma$-stable qui contient $\boldsymbol{h}$.

**Lemme 4.4.5.** $\boldsymbol{h}$ *est un idéal de* $\mathcal{G}$ *et* $\mathcal{P}_c$ *est symplectique.*

PREUVE. On a $[\mathcal{K}_\mathcal{R},\ \mathcal{P}_c] = 0$ donc $[\mathcal{K}_{\boldsymbol{h}},\ \mathcal{P}_\mathcal{R}] = 0$ où $\mathcal{K}_{\boldsymbol{h}} = [\mathcal{P}_c, \mathcal{P}_c]$, dès lors, un raisonnement identique à celui utilisé dans les points (c) et (d) du lemme 4.1.3., livre
$$[\boldsymbol{h},\ \mathcal{R}] = 0$$
La condition $[\mathcal{K}_s,\ \mathcal{P}_c] = \mathcal{P}_c$ implique que $\boldsymbol{h}$ est un idéal de $\mathcal{S}$.
Comme $[\mathcal{K}_{\boldsymbol{h}},\ \mathcal{P}_c] = \mathcal{P}_c$, $\mathcal{P}_c$ est symplectique. ∎

PREUVE DU THÉORÈME 4.4.1.. Il reste à voir $(iv)$. Ce qui est clair car $\boldsymbol{h}$ est contenu dans le facteur semi-simple $\mathcal{L}$ de $\mathcal{G}$ qui est unique. ∎

# Chapitre 5

# Exemples résolubles

## 5.1 Quadruples admissibles

### 5.1.1 Définition

Soient $\boldsymbol{a}$ et $\boldsymbol{h}$ deux algèbres de Lie abéliennes. Soit $\mathcal{G}$ une extension de $\boldsymbol{h}$ par $\boldsymbol{a}$ :

$$\boldsymbol{a} \xrightarrow{i} \mathcal{G} \xrightarrow{p} \boldsymbol{h}$$

$\boldsymbol{a}$ est alors munie canoniquement d'une structure de $\boldsymbol{h}$-module. Notons $\rho : \boldsymbol{h} \mapsto \text{End}(\boldsymbol{a})$ l'homomorphisme de représentation associé à cette structure.

**Définition 5.1.1.1.** *Un quadruple $q = (\boldsymbol{a}, \boldsymbol{h}, \mathcal{G}, s)$ est admissible si*

*(i) $s$ est un automorphisme involutif de $\boldsymbol{a}$;*

*(ii) $\rho$ anticommute avec $s$ : pour tout $h \in \boldsymbol{h}$ :*

$$\rho(h) \circ s + s \circ \rho(h) = 0$$

*(iii) $[\mathcal{G}, \ \mathcal{G}] = i(\boldsymbol{a})$;*

*(iv) $\mathcal{G}$ est une extension inessentielle de $\boldsymbol{h}$ par $\boldsymbol{a}$.*

  *$\boldsymbol{a}$, $\boldsymbol{h}$ et $s$ étant fixés, deux tels quadruples $q_i = (\boldsymbol{a}, \boldsymbol{h}, \mathcal{G}_i, s)$ $(i = 1, 2)$ sont isomorphes si il existe un automorphisme $\alpha$ de $\boldsymbol{h}$ et un automorphisme $\phi$ de $\boldsymbol{a}$ tels que*

  *(a) $[\phi, s] = 0$*

  *(b) $\rho_1(h) = \phi \ \rho_2(\alpha(h)) \ \phi^{-1}$ pour tout $h$ dans $\boldsymbol{h}$ et où $\rho_i$ désigne l'homomorphisme de représentation associée à l'extension $\mathcal{G}_i$ $(i = 1, 2)$.*

**Définition 5.1.1.2.** *Une paire $p = (\mathcal{G}, \sigma)$ est dite symétrique si $\mathcal{G}$ est une algèbre de Lie réelle de dimension finie et $\sigma$ en est un automorphisme involutif.*
*Si $\mathcal{G} = \mathcal{K} \oplus \mathcal{P}$ est la décomposition de $\mathcal{G}$ relativement à $\sigma$, étant donné un quadruple admissible $q = (\boldsymbol{a}, \boldsymbol{h}, \mathcal{G}, s)$, on a le produit semi-direct $\mathcal{G} = \boldsymbol{a} \times_\rho \boldsymbol{h}$; en définissant l'endomorphisme $\sigma$ de $\mathcal{G}$ par*

$$\sigma = s \oplus (-id)_{\boldsymbol{h}}$$

*on voit que $p(q) = (\mathcal{G}, \sigma)$ est une paire symétrique; on l'appelle la paire symétrique associée au quadruple $q$.*

**Définition 5.1.1.3.** *On note*
$$Hom^{(s)}(\boldsymbol{h}, \text{End}(\boldsymbol{a}))$$
*l'ensemble des homomorphismes de $\boldsymbol{h}$ vers $\text{End}(\boldsymbol{a})$ qui anticommutent avec $s$ et $Aut_{(s)}(\boldsymbol{a})$ le centralisateur de $s$ dans $Aut(\boldsymbol{a})$. On voit alors que l'on a deux actions qui commutent :*

$$Hom^{(s)}(\boldsymbol{h}, \text{End}(\boldsymbol{a})) \times GL(\boldsymbol{h}) \to Hom^{(s)}(\boldsymbol{h}, \text{End}(\boldsymbol{a}))$$

$$\rho.\alpha = \rho \circ \alpha^{-1} \qquad (\alpha \text{ dans } GL(\boldsymbol{h}))$$

*et*

$$Aut_{(s)}(\boldsymbol{a}) \times Hom^{(s)}(\boldsymbol{h}, \text{End}(\boldsymbol{a})) \to Hom^{(s)}(\boldsymbol{h}, \text{End}(\boldsymbol{a}))$$

$$(\phi.\rho)(h) = \phi\rho(h)\ \phi^{-1}$$

On a alors la

**Proposition 5.1.1.4.** *Fixons $\boldsymbol{a}, \boldsymbol{h}$ et $s$ comme précédemment.*

(i) *Deux quadruples $q_i = (\boldsymbol{a}, \boldsymbol{h}, \mathcal{G}_i, s)$ ($i = 1, 2$) isomorphes livrent des paires symétriques isomorphes.*

(ii) *L'ensemble des classes d'isomorphie des paires symétriques associées aux quadruples $q = (\boldsymbol{a}, \boldsymbol{h}, \mathcal{G}, s)$ ($\boldsymbol{a}$, $\boldsymbol{h}$ et $s$ fixés) est en bijection avec celui des classes d'isomorphie de ces quadruples et est paramétrisé par le double quotient*

$$Aut_{(s)}(\boldsymbol{a}) \left\backslash Hom^{(s)}(\boldsymbol{h}, \text{End}(\boldsymbol{a})) \right/ GL(\boldsymbol{h})$$

PREUVE.

(i) Soit $(\alpha, \phi)$ un isomorphisme entre $q_1$ et $q_2$.
En considérant $\mathcal{G}_i = \boldsymbol{a} \times_{\rho_i} \boldsymbol{h}$ ($i = 1, 2$), on a l'isomorphisme linéaire

$$\mathcal{G}_1 \xrightarrow{\varphi} \mathcal{G}_2$$

$\varphi = \phi^{-1} \oplus \alpha$.
On vérifie que $\varphi$ est un isomorphisme d'algèbres de Lie. Comme $[\varphi, s] = 0$, on a $\varphi \circ \sigma = \sigma \circ \varphi$ et donc $\varphi$ induit un isomorphisme de $p(q_1)$ sur $p(q_2)$.

(ii) Notons $\boldsymbol{a} = \mathcal{K} \oplus V$ la décomposition relativement à $s$, c'est-à-dire $s = id_{\mathcal{K}} \oplus (-id)_V$.
Soit alors $\varphi$ un isomorphisme de $p(q_1)$ sur $p(q_2)$. On a $\varphi(V) = V$; en effet : $V \subset [\mathcal{G}_1, \mathcal{G}_1]$ donc $\varphi(V) \subset \mathcal{K} \oplus V$; donc $\forall v \in V$

$$\sigma\ \varphi(v) = \varphi\ \sigma(v) = -\varphi(v)$$

Posons $\varphi_V = pr_V \circ \varphi \circ pr_{\boldsymbol{h}}$ et $\varphi_{\boldsymbol{h}} = pr_{\boldsymbol{h}} \circ \varphi \circ pr_{\boldsymbol{h}}$.
Posons encore $\mu = \varphi - \varphi_V$. On a $\forall a \in \boldsymbol{a}$, $h \in \boldsymbol{h}$

$$\begin{aligned} \mu[h,a] &= \varphi([h,a]) - \varphi_V[h,a] = \varphi[h,\ a] \\ &= [\varphi(h), \varphi(a)] = [\varphi_V(h) + \varphi_{\boldsymbol{h}}(h),\ \mu(a)] \\ &= [\varphi_{\boldsymbol{h}}(h),\ \mu(a)] = [\mu(h),\ \mu(a)] \end{aligned}$$

et

$$\forall h' \in \boldsymbol{h} \qquad\qquad \mu[h,\ h'] = 0 = [\mu(h),\ \mu(h')]$$

Donc $\mu$ est un homomorphisme.
Si $0 \neq h \in \boldsymbol{h}$ est tel que $\mu(h) = 0$ alors $\varphi(h) \in V$ ce qui est impossible car $\varphi$ est un isomorphisme et $\varphi(V) = V$.

Dès lors, $\mu$ est un isomorphisme et on a

$$\mu = \varphi|_{\boldsymbol{a}} \oplus \varphi_{\boldsymbol{h}}$$

Le couple $(\alpha = \varphi_{\boldsymbol{h}}, \phi = (\varphi|_{\boldsymbol{a}})^{-1})$ définit alors un isomorphisme entre $q_1$ et $q_2$.

∎

### 5.1.2  $\dim \mathcal{P} = 4$

Dans ce paragraphe, $\boldsymbol{a}$ et $\boldsymbol{h}$ sont des algèbres de Lie réelles abéliennes et $s$ un automorphisme involutif de $\boldsymbol{a}$ tels que si $\boldsymbol{a} = \mathcal{K} \oplus V$ est la décomposition relativement à $s$ ($s = id_\mathcal{K} \oplus (-id)_V$) on ait

(i) $\dim V = \dim \boldsymbol{h} = \dim \mathcal{K} = 2$.

(ii) $\boldsymbol{a}$, $\boldsymbol{h}$, $s$ étant fixés, $q = (\boldsymbol{a}, \boldsymbol{h}, \mathcal{G}, s)$ désigne un quadruple admissible tel que si $p(q) = (\mathcal{G}, \sigma)$ est la paire symétrique associée et si $\mathcal{G} = \mathcal{K} \oplus \mathcal{P}$ est la décomposition relativement à $\sigma$ on ait une action effective de $\mathcal{K}$ sur $\mathcal{P}$.

On fixe alors une base $\{e_1, e_2\}$ de $V$, $\{k_1, k_2\}$ de $\mathcal{K}$ et $\{f_1, f_2\}$ de $\boldsymbol{h}$. Si $\rho : \boldsymbol{h} \to \text{End}(\boldsymbol{a})$ est l'homomorphisme induit par $q$, on note $F_i = \rho(f_i)$ $i = 1, 2$. Dans la base $\{k_1, k_2, e_1, e_2\}$ de $\boldsymbol{a}$, on a

$$F_i = \begin{pmatrix} 0 & Y_i \\ X_i & 0 \end{pmatrix} \qquad (i = 1, 2)$$

où $X_i, Y_i \in Mat(2 \times 2, \mathbb{R})$ (matrices $2 \times 2$ réelles). Si on décrit tout automorphisme $\alpha$ de $\boldsymbol{h}$ et tout automorphisme $\phi$ de $\mathcal{K} \oplus V$ matriciellement dans ces bases par :

$$\alpha = \begin{pmatrix} a & b \\ c & d \end{pmatrix}$$

et

$$\phi = \begin{pmatrix} A & 0 \\ 0 & B \end{pmatrix}$$

$a, b, c, d \in \mathbb{R}$, $ad - bc \neq 0$; $A, B \in GL(\mathbb{R}^2)$, la proposition 5.1.1.4. nous dit que décrire les classes d'isomorphie des paires $p(q)$ équivaut à décrire les classes d'équivalences des couples $(F_1, F_2)$ sous la relation

$$(F_1, F_2) \sim (F_1', F_2') \quad \Leftrightarrow \quad \begin{array}{l} F_1' = \phi(aF_1 + cF_2)\phi^{-1} \\ F_2' = \phi(bF_1 + dF_2)\phi^{-1} \end{array}$$

c'est-à-dire, si

$$F_i' = \begin{pmatrix} 0 & Y_i' \\ X_i' & 0 \end{pmatrix} \qquad (i = 1, 2)$$

$$\begin{array}{rcl} Y_1' &=& A(aY_1 + cY_2)\, B^{-1} \\ X_1' &=& B(aX_1 + cX_2)\, A^{-1} \\ Y_2' &=& A(bY_1 + dY_2)\, B^{-1} \\ X_2' &=& B(bX_1 + dX_2)\, A^{-1} \end{array}$$

**Lemme 5.1.2.1.** *Tout couple $(F_1, F_2)$ est équivalent à un couple*

$$F(x, y, k) = (F_1(x, y, k), F_2(x, y, k))$$

où
$$F_1(x,y,k) = \begin{pmatrix} 0 & x + y\,U_k \\ I & 0 \end{pmatrix}$$
$$F_2(x,y,k) = \begin{pmatrix} U_k & 0 \\ 0 & U_k \end{pmatrix} F_1(x,y,k)$$

où $x, y \in \mathbb{R}$, $k \in \{1, 2, 3, 4\}$ avec

$$U_1 = \begin{pmatrix} 0 & 1 \\ 0 & 0 \end{pmatrix} \quad U_2 = \begin{pmatrix} 0 & 1 \\ 1 & 0 \end{pmatrix} \quad U_3 = \begin{pmatrix} 0 & 1 \\ -1 & 0 \end{pmatrix} \quad U_4 = 0$$

(Si $k = 4$, on pose $y = 0$).

PREUVE. Par effectivité de l'action de $\mathcal{K}$ sur $\mathcal{P}$, on a

$$Ker(F_1\big|_{\mathcal{K}}) \cap Ker(F_2\big|_{\mathcal{K}}) = \{0\}$$

comme $[\boldsymbol{h}\ \mathcal{K}] = V$ ($\mathcal{K} \oplus V = [\mathcal{G}\ \mathcal{G}]$), on peut supposer $\det X_1 \neq 0$.
Dès lors, en choisissant $\alpha = id$ ; $A = I$ ; $B = X_1^{-1}$, on se ramène à $X_1 = I$.
La condition $[F_1,\ F_2] = 0$ livre alors

$$Y_2 = Y_1\ X_2 = X_2\ Y_1,$$

en écrivant

$$X_2 = x + X \qquad Y_1 = y + Y$$

($x, y \in \mathbb{R}$, $X, Y \in sl(2, \mathbb{R})$), on a $[X,\ Y] = 0$ donc $Y = \lambda X$ ; $\lambda \in \mathbb{R}$.
En remplaçant $F_2$ par $F_2 - x\ F_1$, on se ramène à

$$X_1 = I \qquad Y_1 = y + \lambda X \qquad X_2 = X \qquad Y_2 = r\lambda + yX$$

où $X^2 = r\ I$, $r \in \mathbb{R}$,

(a) $r = 0$, dans ce cas, $\exists A \in SL_2(\mathbb{R})$ avec $X = \varepsilon A^{-1} U_1 A$ ou $X = 0$ ($\varepsilon = \pm 1$).
En choisissant $\phi = \begin{pmatrix} A & 0 \\ 0 & A \end{pmatrix}$ et en remplaçant $F_2$ par $\varepsilon F_2$, on se ramène à

$$F_1 = \begin{pmatrix} 0 & a + b\,U_1 \\ I & 0 \end{pmatrix} \qquad F_2 = \begin{pmatrix} U_1 & 0 \\ 0 & U_1 \end{pmatrix} F_1 \quad \text{ou} \quad 0$$

($a, b \in \mathbb{R}$).

(b) $r > 0$; dans ce cas, il existe $A \in SL_2(\mathbb{R})$ avec

$$X = \sqrt{r}\ A^{-1}\ U_2\ A$$

on se ramène alors à

$$F_1 = \begin{pmatrix} 0 & a + b\,U_2 \\ I & 0 \end{pmatrix} \qquad F_2 = \begin{pmatrix} U_2 & 0 \\ 0 & U_2 \end{pmatrix} F_1 \qquad (a, b \in \mathbb{R})$$

(c) $r < 0$; dans ce cas, il existe $A \in SL_2(\mathbb{R})$ avec

$$X = \varepsilon\sqrt{|r|}\ A^{-1} U_3 A \qquad\qquad \varepsilon = \pm 1$$

on se ramène alors à

$$F_1 = \begin{pmatrix} 0 & a + b\,U_3 \\ I & 0 \end{pmatrix} \qquad F_2 = \begin{pmatrix} U_3 & 0 \\ 0 & U_3 \end{pmatrix} F_1 \qquad (a, b \in \mathbb{R})$$

∎

Si on note $p(x, y, k)$ la paire symétrique associée au couple $(F_1(x, y, k),\ F_2(x, y, k))$, on a la

**Proposition 5.1.2.2.** *L'ensemble des classes d'isomorphie des paires $p(q)$ est paramétrisé par l'ensemble suivant :*
$$\{p(\varepsilon, 0, k);\ k = 1, 2, 4;\ \varepsilon = \pm 1\} \cup \{p(0, \varepsilon, \ell);\ \ell = 2, 3;\ \varepsilon = \pm 1\}$$

PREUVE. Soit
$$F_1 = \begin{pmatrix} 0 & x + y\ U \\ I & 0 \end{pmatrix} \qquad F_2 = \underline{U}\ F_1$$

où $\underline{U} = \begin{pmatrix} U & 0 \\ 0 & U \end{pmatrix}$, $U \in \{U_k;\ k = 1, 2, 3, 4\}$ et de même :

$$F'_1 = \begin{pmatrix} 0 & x' + y'\ U \\ I & 0 \end{pmatrix} \qquad F'_2 = \underline{U}\ F'_1$$

On suppose
$$\begin{cases} \phi\ \alpha(F_1)\ \phi^{-1} & = & F'_1 \qquad (1) \\ \phi\ \alpha(F_2)\ \phi^{-1} & = & F'_2 \qquad (2) \end{cases}$$

(a) $U = U_4 = 0$

On voit alors immédiatement que
$$(F_1(x, 0, 4),\ F_2(x, 0, 4)) \sim (F_1(\eta, 0, 4), 0)$$

avec $\eta = -1, 0, +1$

(b) $U = U_k;\ k \neq 4$.

Alors, (1) s'écrit

$$(**) \qquad \begin{cases} A(ax + cy(-1)^{\delta_{3k}}(1 - \delta_{1k}) + (ay + cx)U_k)B^{-1} = x' + y'U_k \\ B(a + c\ U_k)A^{-1} = I \end{cases}$$

donc $B^{-1} = (a + cU_k)A^{-1}$ et on remarque que $[\phi, \underline{U}] = 0$ dès lors :

$$\begin{pmatrix} x' \\ y' \end{pmatrix} = \begin{pmatrix} a^2 + (-1)^{\delta_{3k}}(1 - \delta_{1k})c^2 & 2ac(-1)^{\delta_{3k}}(1 - \delta_{1k}) \\ 2ac & a^2 + (1 - \delta_{1k})(-1)^{\delta_{3k}}c^2 \end{pmatrix} \begin{pmatrix} x \\ y \end{pmatrix}$$

On pose $\alpha(F_2) = \underline{U}\ \alpha(F_1)$.

- $k = 1$. On a $a \neq 0$ et on se ramène à

$$\begin{array}{ll} (F_1(0, \varepsilon, 1),\ F_2(0, \varepsilon, 1)) & \text{ou} \\ (F_1(\varepsilon, 0, 1),\ F_2(\varepsilon, 0, 1)) & \text{ou} \\ (F_1(0, 0, 1),\ F_2(0, 0, 1)) \quad \varepsilon = \pm 1 & \end{array}$$

- $k = 2$. On a $a^2 \neq c^2$ et on se ramène à

$$\begin{array}{ll} (F_1(\varepsilon, 0, 2),\ F_2(\varepsilon, 0, 2)) & \text{ou} \\ (F_1(0, 0, 2),\ F_2(0, 0, 2)) & \text{ou} \\ (F_1(0, \varepsilon, 2),\ F_2(0, \varepsilon, 2)) & \text{ou} \\ (F_1(\varepsilon, \varepsilon', 2),\ F_2(\varepsilon, \varepsilon', 2)) \qquad \varepsilon = \pm 1, \varepsilon' = \pm 1 & \end{array}$$

- $k = 3$. On a $a \neq 0$ ou $c \neq 0$; et on se ramène à

$$\begin{array}{ll} (F_1(0, \varepsilon, 3),\ F_2(0, \varepsilon, 3)) & \text{ou} \\ (F_1(0, 0, 3),\ F_2(0, 0, 3)) \qquad \varepsilon = \pm 1 & \end{array}$$

Soit $S_k(x,y)$ défini par
$$S_k(x,y) = \{g \in Aut(\mathcal{K} \oplus V) \mid g\ F_i(x,y,k)\ g^{-1} = F_i(x,y,k)\ , i = 1,2\}$$
Clairement $S_k(x,y)$ ne dépend pas de $(x,y)$ et en notant $S_k = S_k(x,y)$, on a
$$S_k = \mathbb{R} \cdot Stab_{SL_2(\mathbb{R})}(U_k)$$
Maintenant, si on note $S(F_1, F_2) = \{g \in Aut(\mathcal{K} \oplus V) \mid g\ F_i = F_i\ g\}$, on a $\phi\ S(F_1, F_2)\phi^{-1} = S(\phi\alpha F_1\ \phi^{-1},\ \phi\alpha F_2\phi^{-1})$. De plus, le rang et la signature de la forme bilinéaire $\beta^\rho(X,Y) = tr(\rho(x)\ \rho(y))$ sont invariants sous les transformations $X \mapsto \phi\ \alpha(x)\ \phi^{-1}$. Ceci et la relation $[\boldsymbol{h}\ V] = \mathcal{K}$ permettent de conclure. ∎

Les formes symplectiques sur $\mathcal{P}$ $\mathcal{K}$-invariantes dans la base $\{e_1, e_2, f_1, f_2\}$ s'écrivent :

- pour $p(\varepsilon, 0, 1)$ et $p(0, \varepsilon, 1)$ :

$$\Omega_1(a,b,x) = \begin{pmatrix} 0 & 0 & a & 0 \\ 0 & 0 & b & a \\ -a & -b & 0 & x \\ 0 & -a & -x & 0 \end{pmatrix} \qquad a,b,x \in \mathbb{R}, a \neq 0$$

- pour $p(\varepsilon, 0, 2)$ et $p(0, \varepsilon, 2)$ :

$$\Omega_2(a,b,x) = \begin{pmatrix} 0 & 0 & a & b \\ 0 & 0 & b & a \\ -a & -b & 0 & x \\ -b & -a & -x & 0 \end{pmatrix} \qquad a,b,x, \in \mathbb{R}, a^2 \neq b^2$$

- pour $p(0, \varepsilon, 3)$ :

$$\Omega_3(a,b,x) = \begin{pmatrix} 0 & 0 & a & b \\ 0 & 0 & -b & a \\ -a & b & 0 & x \\ -b & -a & -x & 0 \end{pmatrix} \qquad a,b,x \in \mathbb{R}, a^2 + b^2 \neq 0$$

$p(\varepsilon, 0, 4)$ n'admet pas de telle structure.

Un automorphisme de $p(q)$ agit sur $\mathcal{K} \oplus V \oplus \boldsymbol{h}$ par $(\phi \oplus \alpha) + \lambda$ où $\lambda : \boldsymbol{h} \to V$ et $[\phi, \underline{U}] = 0$.

Par un tel automorphisme, on ramène, dans chaque cas, $\Omega_i(a, b, x)$ à

$$\Omega_x = \Omega_i(1, 0, x) = \begin{pmatrix} 0 & I & \\ -I & 0 & x \\ & -x & 0 \end{pmatrix}$$

on a alors la

**Proposition 5.1.2.3.** *Tout T.S.S. résoluble de dimension 4 de la forme $t = (p(q), \Omega)$ est isomorphe à un des triples non isomorphes suivants :*

$$\begin{array}{rcl}
(p(\varepsilon, 0, 1), \Omega_x) & = & t^{4\ (1)}_{1,\varepsilon,x} \\
(p(\varepsilon, 0, 2), \Omega_x) & = & t^{4\ (1)}_{2,\varepsilon,0,x} \\
(p(0, \varepsilon, 2), \Omega_x) & = & t^{4\ (1)}_{2,0,\varepsilon,x} \\
(p(0, \varepsilon, 3), \Omega_x) & = & t^{4\ (1)}_{3,\varepsilon,x}
\end{array}$$

*où $\varepsilon = \pm 1$, $x \in \mathbb{R}$.*

## 5.2 Quintuples admissibles

### 5.2.1 Définition

Soient $\boldsymbol{a}$ et $\boldsymbol{h}$ deux algèbres de Lie abéliennes. Soit $\mathcal{G}$ une extension de $\boldsymbol{h}$ par $\boldsymbol{a}$

$$\boldsymbol{a} \xrightarrow{i} \mathcal{G} \xrightarrow{p} \boldsymbol{h}$$

Soit $\rho : \boldsymbol{h} \to \text{End}(\boldsymbol{a})$ comme plus haut (cf. 5.1.1), et soit $s$ un automorphisme involutif de $\boldsymbol{a}$.

**Définition 5.2.1.1.** *Un quintuple $Q = (\boldsymbol{a}, \boldsymbol{h}, \mathcal{G}, s, \mu)$ est admissible si*

(i) $\mu$ est une application linéaire de $\boldsymbol{h}$ vers $\mathcal{G}$ telle que $p \circ \mu = id_{\boldsymbol{h}}$

(ii) $[\mathcal{G}\ \mathcal{G}] = i(\boldsymbol{a})$

(iii) $\rho$ anticommute avec $s$

(iv) Si on note $r : \Lambda^2(\boldsymbol{h}) \to \boldsymbol{a}$ l'application définie par $r(h \wedge h') = i^{-1}[\mu(h),\ \mu(h')]$, alors, on a
- $r\ \Lambda^2(\boldsymbol{h}) \cap \rho(\boldsymbol{h})\ \boldsymbol{a} = \{0\}$ et
- $s\ r(h \wedge h') = r(h \wedge h')$ pour tous $h, h'$ dans $\boldsymbol{h}$.

$\boldsymbol{a}$, $\boldsymbol{h}$ et $s$ étant fixés, deux tels quintuples $Q_i = (\boldsymbol{a}, \boldsymbol{h}, \mathcal{G}_i, s, \mu_i)$ $(i = 1, 2)$ sont isomorphes si il existe un automorphisme $\alpha$ de $\boldsymbol{h}$ et un automorphisme $\phi$ de $\boldsymbol{a}$ tels que

(a) $[\phi,\ s] = 0$

(b) $\rho_1(h) = \phi\ \rho_2(\alpha(h))\ \phi^{-1}$ pour tout $h$ dans $\boldsymbol{h}$

(c) $\phi^{-1} \circ r_1 = \alpha \cdot r_2$

Etant donné un quintuple admissible $Q$, en notant $W = \mu(\boldsymbol{h})$ on a la somme directe vectorielle

$$\mathcal{G} = \boldsymbol{a} \oplus W$$

en définissant $\sigma = s \oplus (-id)_W$, on voit que $p(Q) = (\mathcal{G}, \sigma)$ est une paire symétrique, on l'appelle la paire symétrique associée au quintuple $Q$.

**Lemme 5.2.1.2.** *Deux quintuples isomorphes $Q_i = (\boldsymbol{a}, \boldsymbol{h}, \mathcal{G}_i, s, \mu_i)$ $(i = 1, 2)$ livrent des paires symétriques isomorphes.*

PREUVE. Soit $(\alpha, \phi)$ un isomorphisme entre $Q_1$ et $Q_2$. Alors $\varphi = \phi^{-1} \oplus \widetilde{\alpha}$, où

$$\widetilde{\alpha} = \mu_2 \circ \alpha \circ p_1$$

livre un isomorphisme de $p(Q_1)$ sur $p(Q_2)$. ∎

Si $Q = (\boldsymbol{a}, \boldsymbol{h}, \mathcal{G}, s, \mu)$ est un quituple admissible, on note $\mathcal{G} = \mathcal{K} \oplus \mathcal{P}$ la décomposition relativement à $\sigma$ et $\boldsymbol{a} = \mathcal{K} \oplus V$ la décomposition relativement à $s$ (mêmes notations que pour les quadruples).

### 5.2.2   $\dim \mathcal{P} = 4$

On suppose ici $\dim W = \dim V = 2$ et $\dim \mathcal{K} = 2$.

Soient $\{e_1, e_2\} \cup \{f_1, f_2\}$ une base de $V \oplus W$ et $\{k_1, k_2\}$ une base de $\mathcal{P}$, posons $p(f_i) = h_i$.
Soient, dans ces bases,
$$\alpha = \begin{pmatrix} a & b \\ c & d \end{pmatrix}$$
et
$$\phi = A \oplus B$$
$a, b, c, d \in \mathbb{R}$, $A, B \in GL(\mathbb{R}^2)$ avec $ad - bc \neq 0$.
Posons
$$F_i = \rho(h_i) \qquad (i = 1, 2)$$
dans la base $\{k_1, k_2, e_1, e_2\}$ on a
$$F_i = \begin{pmatrix} 0 & Y_i \\ X_i & 0 \end{pmatrix}$$
où $X_i, Y_i$ sont des matrices réelles $2 \times 2$.

Nous supposerons, de plus, que l'action de $\mathcal{K}$ sur $\boldsymbol{h}$ est effective et que $\dim \rho(\boldsymbol{h})V = 1$; donc
$$\begin{array}{rcll} Im\ Y_1 & \subseteq & Im\ Y_2 = \rho(\boldsymbol{h})\ V & \text{ou} \\ Im\ Y_2 & \subseteq & Im\ Y_1 = \rho(\boldsymbol{h})\ V & \end{array}$$

Dès lors, pour classer les quintuples, on peut commencer par paramétriser les classes d'équivalences des couples $F = (F_1, F_2)$ sous la relation $(F_1, F_2) \sim (F_1', F_2')$ si et seulement si
$$\begin{array}{rcl} Y_1' & = & a\ A\ Y_1\ B^{-1} + c\ A\ Y_2\ B^{-1} \\ Y_2' & = & b\ A\ Y_1\ B^{-1} + d\ A\ Y_2\ B^{-1} \\ X_1' & = & a\ B\ X_1\ A^{-1} + c\ B\ X_2\ A^{-1} \\ X_2' & = & b\ B\ X_1\ A^{-1} + d\ B\ X_2\ A^{-1} \end{array}$$
où
$$F_i' = \begin{pmatrix} 0 & Y_i' \\ X_i' & 0 \end{pmatrix} \qquad (i = 1, 2)$$

Un raisonnement identique à celui utilisé pour les quadruples et le fait que $\dim \rho(\boldsymbol{h})V = 1$ livrent le

**Lemme 5.2.2.1.** *Tout couple $F$ est équivalent à un des couples suivants : (mêmes notations que pour les quadruples )*
$$F(0, \varepsilon, 1), \qquad F(\varepsilon, \varepsilon', 2)$$
où $\varepsilon = \pm 1$, $\varepsilon = \pm 1$.

On a alors les tables pour $\mathcal{G}$ :
$$\begin{array}{rcll} [f_1,\ k_1] & = & e_1 & \\ [f_1,\ k_2] & = & e_2 & \\ [f_2,\ k_2] & = & e_1 & \\ [f_1,\ e_1] & = & \varepsilon k_1 & \\ [f_1,\ f_2] & = & a\ k_1 + b\ k_2 & b \neq 0 \end{array}$$

et
$$\begin{aligned}
[f_1,\ k_1] &= e_1 \\
[f_1,\ k_2] &= e_2 \\
[f_2,\ k_1] &= e_2 \\
[f_2,\ k_2] &= e_1 \\
[f_1,\ e_1] &= \varepsilon\, k_1 + \varepsilon'\, k_2 \\
[f_1,\ e_2] &= \varepsilon'\, k_1 + \varepsilon\, k_2 \\
[f_2,\ e_1] &= \varepsilon'\, k_1 + \varepsilon\, k_2 \\
[f_2,\ e_2] &= \varepsilon\, k_1 + \varepsilon'\, k_2 \\
[f_1,\ f_2] &= a(\varepsilon k_1 + \varepsilon'\, k_2) + b(\varepsilon k_1 - \varepsilon'\, k_2) \qquad b \neq 0
\end{aligned}$$
$(a, b \in \mathbb{R})$

Une transformation sur $f_2$ du type
$$f_2 \mapsto f_2 + x\, e_2 \qquad\qquad x \in \mathbb{R}$$
annulle les paramètres "$a$" et ne modifie pas le reste des deux tables.

Pour la première table, les transformations :
$$\begin{array}{lll}
f_1 \leftarrow \varepsilon |b|^{-1/4}\, f_1 & e_1 \leftarrow \varepsilon\varepsilon'|b|^{-3/4} e_1 & k_1 \leftarrow \varepsilon'|b|^{-1/2} k_1 \\
f_2 \leftarrow |b|^{-3/4}\, f_2 & e_2 \leftarrow \varepsilon'|b|^{-1/4}\, e_2 & k_2 \leftarrow \varepsilon\, \varepsilon'\, k_2
\end{array}$$

livre
$$\begin{aligned}
[f_1,\ k_1] &= e_1 \\
[f_1,\ k_2] &= e_2 \\
[f_2,\ k_2] &= e_1 \\
[f_1,\ e_2] &= k_1 \\
[f_1,\ f_2] &= k_2
\end{aligned}$$

Pour la seconde, la transformation
$$\begin{array}{lll}
f_1 \leftarrow \eta\, f_1 & e_1 \leftarrow |b|\, e_1 & k_1 \leftarrow b\, k_1 \\
f_2 \leftarrow \eta\, f_2 & e_2 \leftarrow |b|\, e_2 & k_2 \leftarrow |b|\, k_2
\end{array}$$

(où $\eta$ est le signe de $b$) livre
$$\begin{aligned}
[f_1,\ k_1] &= e_1 \\
[f_1,\ k_2] &= e_2 \\
[f_2,\ k_1] &= e_2 \\
[f_2,\ k_2] &= e_1 \\
[f_1,\ e_1] &= \varepsilon\, k_1 + \varepsilon'\, k_2 \\
[f_1,\ e_2] &= \varepsilon'\, k_1 + \varepsilon\, k_2 \\
[f_2,\ e_1] &= \varepsilon'\, k_1 + \varepsilon\, k_2 \\
[f_2,\ e_2] &= \varepsilon\, k_1 + \varepsilon'\, k_2 \\
[f_1,\ f_2] &= \varepsilon k_1 - \varepsilon'\, k_2
\end{aligned}$$

Notons $p_1$ la paire symétrique associée à la première table et $p_2(\varepsilon, \varepsilon')$ celle associée à la seconde. Dans les deux cas, soit la forme symplectique
$$\Omega_{a,x} = a \begin{pmatrix} 0 & I & \\ -I & 0 & x \\ & -x & 0 \end{pmatrix} \qquad a \in \mathbb{R}_0,\ x \in \mathbb{R}$$

exprimée dans la base $\{e_1, e_2, f_1, f_2\}$ de $\mathcal{P}$.

Comme dans le cas des quadruples, le calcul explicte de $Aut(\mathcal{G}, \sigma)$ livre la

**Proposition 5.2.2.2.** *Tout T.S.S. de dimension 4 de la forme $t = (p(Q), \Omega)$ avec $\dim V = \dim W = \dim \mathcal{K} = 2$ est isomorphe un des triples non isomorphes :*

$$t^{4\ (2)}_{1,a,x} = (p_1, \Omega_{a,x})$$

$$t^{4\ (2)}_{2,\varepsilon,\varepsilon',a,x} = (p_2(\varepsilon, \varepsilon'), \Omega_{a,x})$$

*où $\varepsilon = \pm 1$, $\varepsilon' = \pm 1$, $a \in \mathbb{R}_0$, $x \in \mathbb{R}$.*

# Chapitre 6

# Classification en dimension 2 ou 4

## 6.1 Le cas résoluble en dimension 4

### 6.1.1 Liste exhaustive

Dans ce paragraphe, $(\mathcal{G}, \sigma, \Omega)$ désigne un T.S.S. avec $\mathcal{G}$ résoluble et $\dim \mathcal{P} = 4$.

**Proposition 6.1.1.1.**

(i) $\mathcal{K}$ agit sur $\mathcal{P}$ par endomorphismes nilpotents

(ii) $\mathcal{K}$ stabilise un lagrangien $\mathcal{L}$ de $\mathcal{P}$

(iii) $\dim \mathcal{K} \leqslant 4$

(iv) $[\mathcal{K}, [\mathcal{K}, \mathcal{K}]] = 0$

PREUVE.

(i) $\mathcal{G}$ étant résoluble $\mathcal{DG}$ est nilpotente, et est un idéal de $\mathcal{G}$ contenant $\mathcal{K} = [\mathcal{P}, \mathcal{P}]$. Dès lors, $ad(k)\big|_{\mathcal{P}}$ est un endomorphisme nilpotent.

(ii) Par $(i)$, on peut trouver une base $\{e_1, \ldots, e_4\}$ de $\mathcal{P}$ telle que matriciellement on ait

$$ad(k)\big|_{\mathcal{P}} =: k = \begin{pmatrix} 0 & x & y & z \\ 0 & 0 & a & t \\ 0 & 0 & 0 & u \\ 0 & 0 & 0 & 0 \end{pmatrix}$$

on remarque que $\mathcal{K}\, e_1 = 0$.

1. Supposons $a$ nul
    (a) Supposons $t$ ou $u$ non nul; alors il existe $k_0 \in \mathcal{K}$ tel que $k_0\, e_4 =: p_0 \neq 0$ et $p_0 \neq e_1$; comme $\mathcal{K}\, e_1 = 0$, on a
    $$\Omega(e_1, p_0) = 0$$
    de plus, $\forall k \in \mathcal{K} : k\, p_0 = (tx + uy)\, e_1$ et dès lors $\mathcal{L} = \rangle e_1, p_0 \langle$ est un lagrangien $\mathcal{K}$-stable.
    (b) Supposons $t = u = 0$.
    $\mathcal{L} = \rangle e_1, e_2 \langle$ est $\mathcal{K}$-stable et lagrangien.

2. $a$ est non nul. Alors il existe $k_0 \in \mathcal{K}$ tel que $e_2 = k_0\, e_3 + \lambda\, e_1$ ($\lambda \in \mathbb{R}$) et

$$\begin{array}{rcl} \Omega(e_1, e_2) & = & \Omega(e_1, k\, e_3 + \lambda\, e_1) \\ & = & 0 \quad (\mathcal{K}\, e_1 = 0) \end{array}$$

de plus, $k\, e_2 = x\, e_1$, dès lors $\mathcal{L} = \rangle e_1,\ e_2 \langle$ est un lagrangien $\mathcal{K}$-stable.

(iii) Par (i) et (ii), on peut trouver une base symplectique $\{e_1,\ e_2,\ f_1,\ f_2\}$ de $\mathcal{P}$ telle que

$$\Omega = \begin{pmatrix} 0 & I \\ -I & 0 \end{pmatrix}$$

et $\forall k \in \mathcal{K}$,

$$k = \begin{pmatrix} A & S \\ 0 & -{}^\tau A \end{pmatrix}$$

où $S = {}^\tau S$. $k$ étant nilpotent, $A$ l'est aussi et on peut donc supposer

$$A = \begin{pmatrix} 0 & a \\ 0 & 0 \end{pmatrix} \hspace{4cm} a \in \mathbb{R}$$

Dès lors, $\dim \mathcal{K} \leqslant 4$ et $\mathcal{K}'' = 0$.

∎

Définissons $\mathcal{N}$, la sous-algèbre de $sp(2, \mathbb{R})$ formée des éléments

$$n = \begin{pmatrix} 0 & a & x & y \\ 0 & 0 & y & z \\ 0 & 0 & 0 & 0 \\ 0 & 0 & -a & 0 \end{pmatrix} \quad (\Omega = \begin{pmatrix} 0 & I \\ -I & 0 \end{pmatrix})$$

En notant $n =: a\,A + \frac{x}{2}\,X + y\,Y + z\,Z$ la table de $\mathcal{N}$ est donnée par

$$\begin{array}{rcl} [A\ Y] & = & X \\ [A\ Z] & = & Y \end{array}$$

**Proposition 6.1.1.2.**

(i) $\mathcal{N}$ ne contient qu'une seule sous-algèbre abélienne de dimension 3 : $\boldsymbol{a} = \rangle X, Y, Z \langle$.

(ii) Les autres sous-algèbres de dimension 3 de $\mathcal{N}$ sont isomorphes à l'algèbre de Heisenberg $\mathcal{H}_1$, et sont données par

$$\mathcal{H}^{(\lambda)} := \rangle A + \lambda Z, X, Y \langle \hspace{4cm} (\lambda \in \mathbb{R})$$

(iii) Toute sous-algèbre de dimension $\leqslant 2$ de $\mathcal{N}$ est abélienne.

PREUVE. Remarquons d'abord que le centre $\mathcal{Z}(\mathcal{N})$ de $\mathcal{N}$ est le sous-vectoriel $\mathbb{R}\,X$; et que $\mathcal{N}\big/\mathcal{Z}(\mathcal{N})$ est isomorphe à $\mathcal{H}_1$.

en notant $\pi : \mathcal{N} \to \mathcal{H}_1$ la projection canonique et

$$\begin{array}{rcl} \pi(A) & = & U \\ \pi(Z) & = & V \\ \pi(Y) & = & E \end{array}$$

on a $[U, V] = E$, $E$ central. Soit $\mathcal{S}$ une sous-algèbre de dimension 3 de $\mathcal{N}$.

1. $\dim \pi(\mathcal{S}) = 3$.
   Alors, $\pi(\mathcal{S}) = \mathcal{H}_1$ et l'extension
   $$0 \longrightarrow \mathcal{Z}(\mathcal{N}) \longrightarrow \mathcal{N} \longrightarrow \mathcal{H}_1 \longrightarrow 0$$
   est inessentielle, dès lors
   $$\mathcal{N} = \mathcal{H}_1 \oplus \mathcal{Z}(\mathcal{N})$$
   et $\dim [\mathcal{N}, \mathcal{N}] = 1$ ce qui est faux.

2. $\dim \pi(\mathcal{S}) = 2$, alors $\pi(\mathcal{S}) = \rangle W, E \langle$ avec $W \in \rangle U, V \langle$

   (a) $\mathcal{S}$ est abélienne. Alors
   $$\mathcal{S} = \rangle aA + zZ + \mu Y, yY + \nu Z, X \langle \qquad (a, y, z, \mu, \nu \in \mathbb{R})$$
   avec $ay = 0 = a\nu$.
   La dimension 3 livre $a = 0$ et dès lors $\mathcal{S} = \boldsymbol{a}$.

   (b) $\mathcal{S}$ est non abélienne. Alors
   $$\mathcal{S} = \rangle aA + zZ + \mu X, yY + \nu X, X \langle$$
   c'est-à-dire
   $$\mathcal{S} = \rangle A + \lambda Z, Y, X \langle = \mathcal{H}^{(\lambda)}$$

   Soit maintenant $\mathcal{S}$ de dimension 2

   ($a$) $\dim \pi(\mathcal{S}) = 2$. Donc
   $$\mathcal{S} = \rangle aA + zZ + \mu X, Y + \nu X \langle$$
   Si $a$ est non nul, $\dim \mathcal{S} = 3$ donc $\mathcal{S}$ est abélien.

   ($b$) $\dim \pi(\mathcal{S}) = 1$ livre uniquement des sous-algèbres abéliennes.

■

La proposition 6.1.1.1. nous dit que $\mathcal{K}$ est une sous-algèbre de $\mathcal{N}$. En classifiant les sous-algèbres $\mathcal{K}$ de $\mathcal{N}$ admissibles – c'est-à-dire qui peuvent être vues comme espace de points fixes de "$\sigma$" pour certains triples $(\mathcal{G}, \sigma, \Omega)$ – on obtiendra liste exhaustive des triples.

**Lemme 6.1.1.3.** $\mathcal{N}$ n'est pas admissible.

PREUVE. Si $(\mathcal{G}, \sigma, \Omega)$ est un T.S.S. tel que $\mathcal{G} = \mathcal{N} \oplus \mathcal{P}$ avec $\mathcal{P} = \rangle \{e_1, e_2, f_1, f_2\} \langle$, la table de $\mathcal{G}$ commence par :
$$\begin{array}{rcl}
[A\ Y] & = & X \\
[A\ Z] & = & Y \\
[A\ e_2] & = & e_1 \\
[A\ f_1] & = & -f_2 \\
[X\ f_1] & = & 2e_1 \\
[Y, f_1] & = & e_2 \\
[Y\ f_2] & = & e_1 \\
[Z\ f_2] & = & e_2 \\
& \text{etc...} &
\end{array}$$

Jacobi va nous donner une contradiction.

- $[e_1\ e_2] = 0$ car
$$\begin{array}{rcl}[e_1\ e_2] & = & \frac{1}{2}\left[[X\ f_1],\ e_2\right]\\ & = & -\frac{1}{2}\left[[e_2\ X]\ f_1\right] - \frac{1}{2}\left[[f_1\ e_2]\ X\right]\\ & = & 0\end{array}$$
comme $X$ est central et $[X\ e_2] = 0$.

- $[e_1\ f_2] = 0$ par un argument identique.

- $[e_1\ f_1] = \mu X$ ($\mu \in \mathbb{R}$) car
$$[Z\ [e_1\ f_1]] + [f_1\ [Z\ e_1]] + [e_1\ [f_1\ Z]] = 0$$
Or, $[Z\ e_1] = [f_1\ Z] = 0$ donc $[[e_1\ f_1]\ Z] = 0$.
D'autre part, $[A\ [e_1\ f_1]] = 0$ car
$$[A\ [e_1\ f_1]] + [f_1\ [A\ e_1]] + [e_1\ [f_1\ A]] = 0$$
Or, $[A\ e_1] = 0$ et $[e_1\ [f_1\ A]] = [e_1\ f_2] = 0$. Dès lors, $[e_1\ f_1]$ est central.

- $[e_2\ f_1] = \eta Y$ ($\eta \in \mathbb{R}$) car
$$\begin{array}{rcl}[e_2\ f_1] & = & [[Z\ f_2]\ f_1]\\ & = & -[[f_1\ Z]\ f_2] - [[f_2\ f_1]\ Z]\end{array}$$
Or, $[f_1\ Z] = 0$ et $[\mathcal{N},\ Z] = \mathbb{R}\ Y$.

- $[e_2\ f_2] = \nu X$ ($\nu \in \mathbb{R}$) car
$$[e_2\ [f_1\ A]] = -[A\ [e_2\ f_1]] - [f_1\ [A\ e_2]]$$
mais $[e_2\ f_1] = \eta Y$, $[A\ e_2] = e_1$ et $[e_1\ f_1] = \mu X$. Comme $[A\ Y] = X$ on a
$$[e_2\ [f_1\ A]] = \nu X$$
mais $[e_2\ [f_1\ A]] = [e_2\ f_2]$. Dès lors dim$[\mathcal{P},\ \mathcal{P}] \leqslant 3$, contradiction.

∎

**Lemme 6.1.1.4.** *$a$ n'est pas admissible.*

PREUVE. La table de $\mathcal{G} = a \oplus \mathcal{P}$ commence par
$$\begin{array}{rcl}[X\ f_1] & = & 2e_1\\ [Y\ f_1] & = & e_2\\ [Y\ f_2] & = & e_1\\ [Z\ f_2] & = & e_2\end{array}$$
comme dans le lemme 6.1.1.3., on a :

- $[e_1\ e_2] = [e_1\ f_2] = 0$

- $[e_2\ f_2] = [e_1\ f_1]$ car
$$[Y\ [f_1\ f_2]] = 0 = [f_2\ [Y\ f_1]] + [f_1\ [f_2\ Y]]$$
mais $[Y\ f_1] = e_2$ et $[f_2\ Y] = -e_1$

- $[e_2\ f_1] = 0$ car
$$[Z\ [f_1\ f_2]] = 0 = [f_2\ [Z\ f_1]] + [f_1\ [f_2\ Z]]$$
mais $[Z\ f_1] = 0$ et $[Z\ f_2] = e_2$ dès lors dim$[\mathcal{P},\ \mathcal{P}] \leqslant 2$, une contradiction.

∎

**Lemme 6.1.1.5.**

(i) $\mathcal{H}^{(\lambda)}$ n'est pas admissible si $\lambda \neq 0$

(ii) $\mathcal{H}^{(0)}$ est admissible

PREUVE. Soit $(\mathcal{G}, \sigma, \Omega)$ un T.S.S. tel que

$$\mathcal{G} = \mathcal{H}^{(\lambda)} \oplus \mathcal{P} \qquad \mathcal{H}^{(\lambda)} \simeq \mathcal{H}_1$$

Alors la table de $\mathcal{G}$ commence par

$$\begin{array}{rcl}
[U, e_2] & = & e_1 \\
[U, f_1] & = & -f_2 \\
[U, f_2] & = & \lambda e_2 \\
[V, f_1] & = & e_2 \\
[V, f_2] & = & e_1 \\
[E, f_1] & = & 2e_1
\end{array}$$

où $U := A + \lambda Z$ ; $V := Y$ et $E := X$.

- $[e_1, e_2] = [e_1, f_2] = 0$ par un argument identique à celui du lemme 6.1.1.3.

- $[e_1, f_1] = \mu E$ $(\mu \in \mathbb{R})$ car

$$\begin{array}{rcl}
[U\,[e_1\ f_1]] & = & [[U\ e_1]\ f_1] + [e_1\,[U\ f_1]] \\
& = & 0 - [e_1\ f_2] \\
& = & 0
\end{array}$$

et

$$\begin{array}{rcl}
[V\,[e_1\ f_1]] & = & [[V\ e_1]\ f_1] + [e_1\,[V\ f_1]] \\
& = & 0 + [e_1\ e_2] \\
& = & 0
\end{array}$$

donc $[e_1\ f_1]$ est central dans $\mathcal{K}$; un argument identique livre

- $[e_2\ f_2] = \nu E$ $(\nu \in \mathbb{R})$

Supposons $\lambda \neq 0$.
Alors, $[e_2\ f_1] = \gamma E$ $(\gamma \in \mathbb{R})$ car

$$\begin{array}{rcl}
[U\,[f_1\ f_2]] & = & [[U\ f_1]\ f_2] + [f_1\,[U\ f_2]] \\
& = & 0 + [f_1,\ \lambda f_2]
\end{array}$$

donc $[e_2,\ f_1] = \frac{1}{\lambda}[[f_1,\ f_2]\,U] \in [\mathcal{K}\ \mathcal{K}] = \mathbb{R}\,E$
mais alors $\dim[\mathcal{P}, \mathcal{P}] \leqslant 2$ ce qui rend $\mathcal{H}^{(\lambda)}$ non admissible pour $\lambda \neq 0$.

On a donc $\lambda = 0$ et

- $[[f_1,\ f_2]\ U] = 0$

- $[U,\ f_2] = 0$

- $[[e_2,\ f_1]\ V] = 0$ car

$$\begin{array}{rcl}
[V\,[e_2\ f_1]] & = & [[V\ e_2]\ f_1] + [e_2\,[V,\ f_1]] \\
& = & 0 + 0 \\
& = & 0
\end{array}$$

$\mathcal{G}$ a donc une table de la forme

$$\begin{aligned}
[U,\ V] &= E \\
[U,\ e_2] &= e_1 \\
[U,\ f_1] &= -f_2 \\
[V,\ f_1] &= e_2 \\
[V,\ f_2] &= e_1 \\
[E,\ f_1] &= 2\,e_1 \\
[e_2,\ f_2] &= \nu E \\
[e_1,\ f_1] &= \mu E \\
[f_1,\ f_2] &= u\,U + \varepsilon\,E \\
[e_2,\ f_1] &= vV + \eta U
\end{aligned}$$

Maintenant

$$\begin{aligned}
[U\,[e_2\,f_1]] + [f_1\,[U\,e_2]] + [e_2\,[f_1\,U]] &= 0 \\
[V\,[f_1\,f_2]] + [f_2\,[V\,f_1]] + [f_1\,[f_2\,V]] &= 0 \\
[e_2\,[f_1,\ f_2]] + [f_2\,[e_2\,f_1]] + [f_1\,[f_2\,e_2]] &= 0
\end{aligned}$$

fournit le système

$$\begin{pmatrix} 0 & 1 & -1 \\ -1 & 0 & 1 \\ 1 & 1 & 0 \end{pmatrix} \begin{pmatrix} u \\ v \\ \mu \end{pmatrix} = \begin{pmatrix} -\nu \\ \nu \\ 2\nu \end{pmatrix}$$

dont l'unique solution est

$$\mu = 2\nu\ ; \qquad u = v = \nu$$

ceci impose en particulier $\nu \in \mathbb{R}_0$.

On vérifie alors aisément que la table (I)

$$\begin{aligned}
[U,\ V] &= E \\
[U,\ e_2] &= e_1 \\
[U,\ f_1] &= -f_2 \\
[U,\ f_1] &= e_2 \\
[V,\ f_2] &= e_1 \\
[E,\ f_1] &= 2e_1 \\
[e_1,\ f_1] &= 2\nu E \\
[e_2,\ f_1] &= \nu V + \eta E \\
[e_2,\ f_2] &= \nu E \\
[f_1,\ f_2] &= \nu U + \widetilde{\varepsilon}\,E \qquad \widetilde{\varepsilon}, \eta \in \mathbb{R}, \nu \in \mathbb{R}_0
\end{aligned}$$

définit sur $\mathcal{G} = \mathcal{H}^{(0)} \oplus \mathcal{P}$ une structure d'algèbre de Lie. ∎

**Lemme 6.1.1.6.** *Toute sous-algèbre de dimension 2 de $\mathcal{N}$ est admissible.*

PREUVE. Dans ce cas, $\mathcal{K}$ est abélienne (cf. proposition 6.1.1.2.) et on a vu dans la preuve de la proposition 6.1.1.2. que deux cas se présentent :

1. $\mathcal{K} \not\subset \boldsymbol{a}$ et $X \in \mathcal{K}$ ou bien

2. $\mathcal{K} \subset \boldsymbol{a}$

1. Alors, on peut supposer que $\mathcal{K}$ est engendrée par deux éléments de $\mathcal{N}$ de la forme

$$U = \begin{pmatrix} 0 & 1 & 0 & y \\ 0 & 0 & y & z \\ 0 & 0 & 0 & 0 \\ 0 & 0 & -1 & 0 \end{pmatrix} \qquad X = \begin{pmatrix} 0 & 0 & 1 & 0 \\ 0 & 0 & 0 & 0 \\ 0 & 0 & 0 & 0 \\ 0 & 0 & 0 & 0 \end{pmatrix}$$

(où $y, z \in \mathbb{R}$) dans la base $\{e_1, e_2, f_1, f_2\}$. c'est-à-dire

$$\begin{aligned}
[U, e_2] &= e_1 \\
[U, f_1] &= y\, e_2 - f_2 \\
[U, f_2] &= y\, e_1 + z\, e_2 \\
[X, f_1] &= e_1
\end{aligned}$$

On a alors, comme plus haut, $[e_1, e_2] = [e_1\, f_2] = 0$ et, $[e_1, f_1] = [e_2, f_2] = \alpha X$ ($\alpha \in \mathbb{R}$) car comme $[\mathcal{K}, \mathcal{K}] = 0$, $[[U\ e_2]\ f_1] = [[U\ f_1]\ e_2]$ donc

$$[e_1\ f_1] = -[f_2\ e_2]$$

et,

- si $y$ ou $z$ est non nul,

$$[[e_1\ f_1]\ f_2] + [[f_2\ e_1]\ f_1] + [[f_1\ f_2]\ e_1] = 0$$

  livre $[[e_1\ f_1]\ f_2] = 0$ (car $[\mathcal{K}\ e_1] = 0$)

- si $y = z = 0$,

$$[[e_2\ f_1]\ f_2] + [[f_2\ e_2]\ f_1] + [[f_1\ f_2]\ e_2] = 0$$

  livre $[[f_2\ e_2]\ f_1] = \xi e_1$ (car $[\mathcal{K}\ e_2] = \mathbb{R}\ e_1$) ce qui impose à la composante en $U$ de $[e_2\ f_2]$ de s'annuler.

Dès lors, en posant $[f_1\ f_1] = uU + xX$ et $[e_2\ f_1] = \widetilde{u}\, U + \widetilde{x}\, X$

$$[[e_2\ f_1]\ f_2] + [[f_2\ e_2]\ f_1] + [[f_1\ f_2]\ e_2] = 0$$

livre

$$\begin{cases} \widetilde{u}\, z = 0 \\ \widetilde{u}\, y - \alpha + u = 0 \end{cases}$$

(a) $\underline{z \neq 0}$. Alors, $\widetilde{u} = \widetilde{x} = 0$ car, comme $[\mathcal{K}\ \mathcal{K}] = 0$,

$$[[U\ f_1]\ f_2] = [[U\ f_2]\ f_1]$$

c'est-à-dire

$$y\,[e_2\ f_2] = [ye_1 + ze_2, f_1]$$

c'est-à-dire

$$z\,[e_2\ f_1] = 0$$

et on vérifie que la table (II)

$$\begin{aligned}
[U\ e_2] &= e_1 \\
[U\ f_1] &= y\, e_2 - f_2 \\
[U\ f_2] &= y\, e_1 + z\, e_2 \\
[X\ f_1] &= e_1 \\
[e_1\ f_1] &= [e_2\ f_2] = u\, X \qquad (\alpha = u) \\
[f_1\ f_2] &= uU + xX
\end{aligned}$$

(avec $u, z \in \mathbb{R}_0$; $x, y \in \mathbb{R}$) fournit une structure d'algèbre de Lie à $\mathcal{G} = \mathcal{K} \oplus \mathcal{P}$.

(b) $\underline{z = 0}$. Alors, la table (III) :

$$\begin{array}{rcl}
[U\ e_2] & = & e_1 \\
[U\ f_1] & = & y\ e_2 - f_2 \\
[U\ f_2] & = & y\ e_1 \\
[X\ f_1] & = & e_1 \\
[e_1\ f_1] & = & [e_2\ f_2] = (u + \widetilde{u}y)\ X \\
[f_1\ f_2] & = & uU + xX \\
[e_2\ f_1] & = & \widetilde{u}\ U + \widetilde{x}\ X
\end{array}$$

(avec $u, \widetilde{u}, x, \widetilde{x}, y \in \mathbb{R}$) fournit une structure d'algèbre de Lie à $\mathcal{G} = \mathcal{K} \oplus \mathcal{P}$.

2. ($\mathcal{K} \subset \boldsymbol{a}$). Alors, on peut supposer que $\mathcal{K}$ est engendré par deux éléments de $\mathcal{N}$ de la forme

$$U = aX + bZ; \qquad V = xX + yY + zZ$$

c'est-à-dire

$$\begin{array}{rcl}
[U\ f_1] & = & a\ e_1 \\
[U\ f_2] & = & b\ e_2 \\
[V\ f_1] & = & x\ e_1 + y\ e_2 \\
[V\ f_2] & = & y\ e_1 + z\ e_2
\end{array}$$

en particulier, en notant $\mathcal{L} = \rangle\{e_1, e_2\}\langle$, on a

$$[\mathcal{K}\ \mathcal{L}] = 0 \quad \text{et} \quad [\mathcal{K}\ \mathcal{P}] \subseteq \mathcal{L}$$

Supposons $[\mathcal{K}\ \mathcal{P}] \subsetneq \mathcal{L}$ c'est-à-dire $\dim[\mathcal{K}\ \mathcal{P}] = 1$.
Dans ce cas, comme $[U\ \mathcal{P}] \neq 0$, on a deux possibilités.

  (a) $a = 0$, $b \neq 0$ ou bien
  (b) $a \neq 0$, $b = 0$.

Considérons le cas (a). On peut dès lors supposer $b = 1$, $z = 0$ et $[\mathcal{K}\ \mathcal{P}] = \mathbb{R}e_2$, donc $y = x = 0$. Mais alors $[V\ \mathcal{P}] = 0$, une contradiction. Le cas (b) se traite de la même manière et on a

$$\begin{array}{rcl}
[\mathcal{K}\ \mathcal{P}] & = & \mathcal{L} \\
[\mathcal{K}\ \mathcal{L}] & = & 0 \\
[\mathcal{K}\ \mathcal{K}] & = & 0 \quad \text{et par Jacobi} \\
[\mathcal{L}\ \mathcal{L}] & = & 0
\end{array}$$

On a

- $a[e_1\ f_2] = b[e_2\ f_1]$ car
$$[[U\ f_1]\ f_2] = [[U\ f_2]\ f_1]$$

- $x[e_1\ f_2] - z[e_2\ f_1] = y([e_1\ f_1] - [e_2\ f_2])$ car
$$[[V\ f_1]\ f_2] = [[V\ f_2]\ f_1]$$

Deux cas se présentent :

(a) $[\mathcal{L}, \rangle f_1, f_2\langle] = \mathcal{K}$. Alors on peut trouver $\ell_1, \ell_2 \in \mathcal{L}$ tels que
$$[f_1 + \ell_1,\ f_2 + \ell_2] = 0$$

Dès lors, l'extension
$$[\mathcal{G}\ \mathcal{G}] \to \mathcal{G} \to \mathcal{G}\big/_{[\mathcal{G}\ \mathcal{G}]}$$

est inessentielle; donc, en posant
$$\boldsymbol{a} = [\mathcal{G}\ \mathcal{G}] = \mathcal{K} \oplus \mathcal{L}, \qquad s = \sigma\big|_{\boldsymbol{a}}$$

et
$$\boldsymbol{h} = \mathcal{G}\big/_{\boldsymbol{a}};$$

$q = (\boldsymbol{a}, \boldsymbol{h}, \mathcal{G}, s)$ est un quadruple admissible avec $\dim \mathcal{K} = 2$.

(b) $[\mathcal{L}, \rangle f_1, f_2\langle] \neq \mathcal{K}$.
Soit l'extension
$$\boldsymbol{a} \to \mathcal{G} \to \boldsymbol{h}$$

où $\boldsymbol{a} = \mathcal{K} \oplus \mathcal{P} = [\mathcal{G}\ \mathcal{G}]$ et $\boldsymbol{h} = \mathcal{G}/_{\boldsymbol{a}}$. Soit $\mu : \boldsymbol{h} \to \mathcal{G}$ une application linéaire injective avec $\mu(\boldsymbol{h}) = \rangle f_1, f_2\langle$, alors en notant $s = \sigma\big|_{\boldsymbol{a}}$ on a que $Q = (\boldsymbol{a}, \boldsymbol{h}, \mathcal{G}, s, \mu)$ est un quintuple admissible avec $\dim \mathcal{K} = 2$.

**Lemme 6.1.1.7.** *Toute sous-algèbre de dimension 1 de $\mathcal{N}$ est admissible.*

PREUVE. On a $\mathcal{K} = \mathbb{R}\ U$. Dès lors, deux cas se présentent :

1. $[U\ e_2] \neq 0$
2. $[U\ e_2] = 0$

1. Dans ce cas, on peut poser
$$\begin{array}{rcl}
[U\ e_2] &=& e_1 \\
[U\ f_1] &=& x\ e_1 + y\ e_2 - f_2 \\
[U\ f_2] &=& y\ e_1 + z\ e_2 \\
[e_1\ e_2] &=& a\ U \\
[e_1\ f_1] &=& b\ U \\
[e_1\ f_2] &=& c\ U \\
[e_2\ f_1] &=& d\ U \\
[e_2\ f_2] &=& e\ U \\
[f_1\ f_2] &=& f\ U
\end{array}$$

Dès lors
$$\oint [[U\ p]p'] = 0\ et \qquad \oint [[p\ p']p''] = 0 \qquad\qquad (p, p', p'' \in \mathcal{P})$$

fournit la table :
$$\begin{array}{rcl}
[U\ e_2] &=& e_1 \\
[U\ f_2] &=& y\ e_1 \\
[U\ f_1] &=& x\ e_1 + y\ e_2 - f_2 \\
[e_2\ f_1] &=& dU \qquad (d \neq 0) \\
[f_1\ f_2] &=& -dy\ U
\end{array}$$

dont on vérifie qu'elle munit $\mathcal{G} = \mathbb{R}\ U \oplus \mathcal{P}$ d'une structure d'algèbre de Lie.

2. ($[U\ e_2] = 0$). Dans ce cas, on ramène, par un élément de $O(2) \subset Sp(2,\mathbb{R})$, $U$ à la forme
$$U = x\ X + z\ Z$$
$$\oint [[p\ p']p''] = 0 = \oint [[U\ p]p''] \qquad (p, p', p'' \in \mathcal{P})$$

livrent les tables

(a) $[U\ f_2] = e_2$
$[e_1\ f_2], [e_2\ f_2], [f_1\ f_2] \in \mathcal{K}$

(b) $[U\ f_1] = e_1$
$[U\ f_2] = z\ e_2 \quad (z \neq 0)$
$[f_1\ f_2] \in \mathcal{K}$

(c) $[U\ f_1] = e_1$
$[e_1\ f_1], [e_2\ f_1], [f_1\ f_2] \in \mathcal{K}$

dont on vérfie qu'elles munissent $\mathcal{G} = \mathcal{K} \oplus \mathcal{P}$ de structures d'algèbres de Lie.

∎

### 6.1.2 Isomorphismes

Dans cette section, on déduit une liste non redondante de la liste exhaustive obtenue dans la section 6.1.1. Nous reprenons les mêmes notations et numérotations que dans cette section.

**La table (I) dans la preuve du lemme 6.1.1.5.**

En posant
$$\mu = \sqrt{|\nu|}, \varepsilon = \frac{\nu}{\mu^2} \qquad (\varepsilon = \pm 1)$$

et en effectuant la transformation
$$U \longleftarrow \mu\varepsilon U - \frac{\varepsilon}{3\mu}\ E$$
$$V \longleftarrow \mu\varepsilon V - \frac{\eta}{3\mu}\ E$$
$$E \longleftarrow \nu\varepsilon E$$
$$e_1 \longleftarrow \mu\ e_1 \qquad e_2 \longleftarrow \varepsilon e_2 - \varepsilon\frac{2\eta}{3\nu}\ e_1$$
$$f_1 \longleftarrow \mu^{-1}\ f_1 \qquad f_2 \longleftarrow \varepsilon f_2 + \frac{2}{3\nu}\ e_1$$

la table (I) prend la forme
$$\begin{array}{rcl}
[U\ V] & = & E \\
[U\ e_2] & = & e_1 \\
[U\ f_1] & = & -f_2 \\
[V\ f_1] & = & e_2 \\
[V\ f_2] & = & e_1 \\
[E\ f_1] & = & 2\ e_1 \\
[e_1\ f_1] & = & 2\varepsilon E \\
[e_2\ f_1] & = & \varepsilon V \\
[e_2\ f_2] & = & \varepsilon E \\
[f_1\ f_2] & = & \varepsilon U + aE
\end{array} \qquad (a \in \mathbb{R})$$

On remarque que si $\mathcal{D} = [\mathcal{G}\ \mathcal{G}]$ on a le produit semi-direct

$$\mathcal{G} = \mathcal{D} \times_{\varphi_a} \mathbb{R}$$

où $\varphi_a(1) = ad(f_1)\big|_{\mathcal{D}}$.

Dans la base $d = \{U, V, E, e_2, f_2, e_1\}$ de $\mathcal{D}$, on a

$$\varphi_a = \begin{pmatrix} 0 & 0 & 0 & 0 & \varepsilon & 0 \\ 0 & 0 & 0 & -\varepsilon & 0 & 0 \\ 0 & 0 & 0 & 0 & a & -2\varepsilon \\ 0 & -1 & 0 & 0 & 0 & 0 \\ 1 & 0 & 0 & 0 & 0 & 0 \\ 0 & 0 & -2 & 0 & 0 & 0 \end{pmatrix}$$

**Lemme 6.1.2.1.** *Soient $\mathcal{G}_i$ $(i = 1, 2)$ deux extensions inessentielles de $\mathbb{R}$ par une algèbre de Lie $\mathcal{D}$ :*

$$\begin{array}{ccccccccc}
& & & & \mathcal{G}_1 & & & & \\
& & \nearrow^{j_1} & & & \searrow^{\pi_1} & & & \\
0 & \longrightarrow & \mathcal{D} & & & & \mathbb{R} & \longrightarrow & 0 \\
& & \searrow_{j_2} & & & \nearrow_{\pi_2} & & & \\
& & & & \mathcal{G}_2 & & & &
\end{array}$$

*Supposons $[\mathcal{G}_i, \mathcal{G}_i] = j_i(\mathcal{D})$ $(i = 1, 2)$. Notons $\varphi_i : \mathbb{R} \to \mathcal{D}er(\mathcal{D})$ l'homomorphisme déterminé par $\mathcal{G}_i$ et $D_i = \varphi_i(1)$. Si $\mathcal{G}_1$ est isomorphe à $\mathcal{G}_2$, il existe $A \in Aut(\mathcal{D})$, $\xi \in \mathbb{R}_0$ et $X \in \mathcal{D}$ tels que*

$$D_2 = \xi\ A\ D_1\ A^{-1} + ad(X)$$

Peut-on, par une transformation $\varphi$ de $\mathcal{G}$ préservant la décomposition $\mathcal{G} = \mathcal{K} \oplus \mathcal{P}$, modifier le paramètre "$a$" en conservant le reste de la table ?

On a $\varphi \mathcal{D} = \mathcal{D}$ et comme $\mathcal{Z}(\mathcal{D}) = \rangle E, e_1 \langle$, on a dans la base $\{U, V, E\} \cup \{e_2, f_2, e_1\}$ :

$$\varphi\big|_{\mathcal{D}} = A \oplus B$$

où

$$A = \begin{pmatrix} \alpha & 0 \\ & 0 \\ {}^\tau\beta & k \end{pmatrix} \quad \text{et} \quad B = \begin{pmatrix} \phi & 0 \\ & 0 \\ {}^\tau v & r \end{pmatrix}$$

avec $\alpha, \beta \in GL_2(\mathbb{R})$, $k, r \in \mathbb{R}_0$ et $\beta, v \in \mathbb{R}^2$.

On a

$$-ad(x\ e_2 + y\ f_2)\big|_{\mathcal{D}} = \begin{pmatrix} 0 & 0 & 0 & 0 & 0 & 0 \\ 0 & 0 & 0 & 0 & 0 & 0 \\ 0 & 0 & 0 & \varepsilon y & -\varepsilon x & 0 \\ 0 & 0 & 0 & 0 & 0 & 0 \\ 0 & 0 & 0 & 0 & 0 & 0 \\ x & y & 0 & 0 & 0 & 0 \end{pmatrix}$$

Le lemme 6.1.2.1. nous conduit alors à regarder les orbites

$$\varphi_{a'} = \xi(\varphi\big|_{\mathcal{D}}\ \varphi_a(\varphi\big|_{\mathcal{D}})^{-1}) - ad(x\ e_2 + y\ f_2)\big|_{\mathcal{D}}$$

c'est-à-dire, en posant $\mathcal{J} = \begin{pmatrix} 0 & 1 \\ -1 & 0 \end{pmatrix}$,

(a) $\varepsilon \mathcal{J} = \xi \varepsilon \ \alpha \mathcal{J} \ \phi^{-1}$

(b) $\underline{a}' = \xi(\varepsilon^\tau \beta \mathcal{J} + k\underline{a})\phi^{-1} + \xi\frac{2k}{r}\ \varepsilon^\tau v \phi^{-1} - {}^\tau u \mathcal{J}$

(c) $-2\varepsilon = -\frac{2k}{r}\ \xi\varepsilon$

(d) $-\mathcal{J} = -\xi\phi\ \mathcal{J}\alpha^{-1}$

(e) $0 = -\xi\ {}^\tau v \mathcal{J}\ \alpha^{-1} + \xi\frac{2r}{k}\ {}^\tau\beta\alpha^{-1} + {}^\tau u$

(f) $-2 = -\xi\ \frac{2r}{k}$

où $\underline{a} = (0, a)$ et $u = (x, y)$.

(a)(c)(d) et (f) livrent
$$k = r \qquad \xi = 1 \qquad \phi = -\mathcal{J}\alpha\mathcal{J}$$

Dès lors, on a ${}^\tau u = ({}^\tau v\ \mathcal{J} - 2{}^\tau\beta)\alpha^{-1}$ et
$$\underline{a}'\phi = k\underline{a} + (\varepsilon + 2)^\tau \beta\ \mathcal{J} + (1 + 2\varepsilon)^\tau v$$

De là, on voit que l'on peut annuler "a" et que les valeurs $\varepsilon = +1$ et $\varepsilon = -1$ livrent des paires symétriques, notées $p_\varepsilon^{(3)}$, non isomorphes.

Dans la base $\{e_1, f_1, e_2, f_2\}$ la forme symplectique prend la forme
$$\Omega_{x,y,z} = \begin{pmatrix} x\ \mathcal{J} & \begin{matrix} 0 & 0 \\ -y & z \end{matrix} \\ \begin{matrix} 0 & y \\ 0 & -z \end{matrix} & x\ \mathcal{J} \end{pmatrix}$$

où $x \in \mathbb{R}_0$, $y, z \in \mathbb{R}$.

En vertu de ce qui précède, une transformation de $\mathcal{P}$ du type
$$\varphi = \begin{pmatrix} r\,I & 0 \\ 0 & \phi \end{pmatrix} \qquad\qquad\text{où } \phi \in GL_2(\mathbb{R})$$

et $r = \det \phi$, s'étend en un automorphisme de $(\mathcal{G}, \sigma)$

Dès lors, on peut ramener la forme symplectique à
$$\Omega_{\varepsilon', 0, \eta}$$

où $\eta = 0, 1$ et $\varepsilon' = \pm 1$.

Remarquons que les triples obtenus ne sont pas exacts et que $\mathcal{Z}(\mathcal{G}) = \{0\}$. Nous noterons $t_{\varepsilon, \varepsilon', \eta}^{4\ (3)}$ le triple $(p_\varepsilon^3, \Omega_{\varepsilon', 0, \eta})$.

**La table (II) dans la preuve du lemme 6.1.1.6.**

En posant $u = \varepsilon\ \nu^2$, $\frac{z}{u} = \eta\mu^2$ et en effectuant le changement de base

$$\begin{array}{rcl} U & \leftarrow & \frac{1}{\mu\ \nu}(U + \frac{x}{2u}\ X) \\ X & \leftarrow & \frac{\eta}{\mu\nu}\ X \\ e_1 & \leftarrow & \frac{\eta}{\mu\nu^2}\ e_1 \\ e_2 & \leftarrow & \frac{\eta}{\nu}e_2 \\ f_1 & \leftarrow & \frac{1}{\nu}(f_1 - \frac{x}{2u}\ e_2) \\ f_2 & \leftarrow & \frac{1}{u\mu}\ (f_2 - y\ e_2) \end{array}$$

où $\eta = \pm 1$ et $\varepsilon = \pm 1$, on obtient la table :

$$\begin{array}{rcl} [U, \ e_2] & = & e_1 \\ [U, \ f_1] & = & -\varepsilon \ f_2 \\ [U, \ f_2] & = & e_2 \\ [X, \ f_1] & = & e_1 \\ [e_1, \ f_1] & = & \varepsilon X \\ [e_2, \ f_2] & = & X \\ [f_1, \ f_2] & = & U \end{array} \qquad \varepsilon = \pm 1$$

Comme dans le cas précédent, on remarque que

$$\mathcal{G} = \mathcal{D} \times_\varphi \mathbb{R}$$

où $\mathcal{D} = [\mathcal{G} \ \mathcal{G}] = \rangle U, X, e_2, f_2, e_1 \langle$ et dans la base $\{U, X, e_1, e_2, f_2\}$,

$$\varphi = ad(f_1)\big|_\mathcal{D} = \begin{pmatrix} 0 & 0 & 0 & 0 & 1 \\ 0 & 0 & -\varepsilon & 0 & 0 \\ 0 & -1 & 0 & 0 & 0 \\ 0 & 0 & 0 & 0 & 0 \\ \varepsilon & 0 & 0 & 0 & 0 \end{pmatrix}$$

De nouveau, en utilisant le lemme 6.1.2.1.; on voit que les valeurs $\varepsilon = +1$ et $\varepsilon = -1$ livrent des algèbres non isomorphes.
Il est intéressant de remarquer que dans ce cas-ci, toute forme symplectique $\mathcal{K}$-invariante sur $\mathcal{P}$ est exacte, c'est-à-dire

$$\Omega = \delta(a \ X_\star + b \ U_\star)\big|_{\mathcal{P} \times \mathcal{P}}$$

où $a \in \mathbb{R}_0$, $b \in \mathbb{R}$; et où l'on définit

$$\begin{array}{rcl} X_\star(x \ X + u \ U) & = & x \\ U_\star(x \ X + u \ U) & = & u \end{array} \qquad \forall x, u \in \mathbb{R}$$

Nous noterons $t_{\varepsilon,\alpha}^{4\ (4)}$ les triples ainsi construits.

## La table (III) dans la preuve du lemme 6.1.1.6.

Posons $a = \widetilde{u} \ y + u$, $\mu = \sqrt{|a|}$, $a = \varepsilon \ \mu^2$ et $b = x + \widetilde{x} \ y$.
   Deux cas se présentent

(a) $a \neq 0$. En affectant la transformation

$$\begin{array}{rcl} U & \leftarrow & \mu U \\ X & \leftarrow & \mu^2 \ X \\ e_1 & \leftarrow & \mu e_1 \\ e_2 & \leftarrow & (1 - \frac{\widetilde{u} \ y}{a}) e_2 + (\frac{\widetilde{u} \ b}{2a^2} - \frac{\widetilde{x}}{a}) \ e_1 + \frac{\widetilde{u}}{a} \ f_2 \\ f_1 & \leftarrow & \frac{1}{\mu}(f_1 - \frac{b}{2a} \ e_2) \\ f_2 & \leftarrow & f_2 - y \ e_2 + \frac{b}{2a} \ e_1 \end{array}$$

on obtient la table

$$\begin{array}{rcl} [U \ e_2] & = & e_1 \\ [U \ f_1] & = & -f_2 \\ [X \ f_1] & = & e_1 \\ [e_1 \ f_1] & = & \varepsilon X \\ [e_2 \ f_2] & = & \varepsilon X \\ [f_1 \ f_2] & = & \varepsilon U \end{array} \qquad \varepsilon = \pm 1$$

L'ensemble des valeurs propres de la forme de Killing, dans la base $\{U, X, e_1, f_1, e_2, f_2\}$ est $\{0, 4\varepsilon\}$; dès lors les valeurs $\varepsilon = +1$ et $\varepsilon = -1$ livrent des paires symétriques, $p_\varepsilon^{(5)}$, non isomorphes.

Un automorphisme du type

$$\begin{aligned} e_1 &\leftarrow r\,e_2 \\ X &\leftarrow r\,X \\ e_1 &\leftarrow r\,e_1 \\ f_1 &\leftarrow f_1 \\ U &\leftarrow U \\ f_2 &\leftarrow f_2 \end{aligned}$$

ramène toute forme symplectique sur $\mathcal{K}$-invariante sur $\mathcal{P}$ à une des formes suivantes :

$$\Omega_{\alpha,\eta} = (\delta\alpha)\big|_{\mathcal{P}\times\mathcal{P}} + \eta\, de_2 \wedge df_1$$

où $\eta = -1, 0, 1$ et $\alpha = X_\star + u\,U_\star$. Nous noterons $t_{\varepsilon,\alpha,\eta}^{4\ (5)}$ le triple $(p_\varepsilon^{(5)}, \Omega_{\alpha,\eta})$.

(b) $a = 0$. On a $u = -\widetilde{u}y$, dès lors si $b = 0$, $x = -\widetilde{x}y$ et $[f_1\ f_2]$ est proportionnel à $[e_2\ f_1]$, ce qui contredit $[\mathcal{P}\ \mathcal{P}] = \mathcal{K}$.

De plus, si $\widetilde{u}$ était nul, $u$ le serait ce qui contredit $[\mathcal{P}\ \mathcal{P}] = \mathcal{K}$. On a donc $b \neq 0 \neq \widetilde{u}$.

1. $x \neq 0$.
   La transformation
   $$\begin{aligned} U &\leftarrow \eta\,\widetilde{\mu}\nu U \\ X &\leftarrow \widetilde{u}\,b^2\,X \\ e_1 &\leftarrow \eta\widetilde{u}\,b^2\,e_1 \\ e_2 &\leftarrow \eta(x\,e_2 + \widetilde{x}\,f_2) \\ f_1 &\leftarrow \eta\,f_1 \\ f_2 &\leftarrow \widetilde{\mu}\nu(f_2 - y\,e_2) \end{aligned}$$
   où $\widetilde{\mu} = |b|$, $b = \varepsilon\widetilde{\mu}$, $|\widetilde{u}| = \nu$ et $\widetilde{u} = \varepsilon\eta\nu$; livre la table
   $$\begin{aligned} [U\ e_2] &= e_1 \\ [U\ f_1] &= -f_2 \\ [X\ f_1] &= e_1 \\ [f_1\ f_2] &= X \\ [e_2\ f_1] &= U \end{aligned}$$
   En calculant explicitement le groupe des automorphismes de $(\mathcal{G}, \sigma)$ on voit que l'on peut se ramener à la forme symplectique
   $$\Omega_\varepsilon = \varepsilon\begin{pmatrix} 0 & I \\ -I & 0 \end{pmatrix}$$
   dans la base $\{e_1, e_2, f_1, f_2\}$ où $\varepsilon = \pm 1$. Nous noterons $t_\varepsilon^{4\ (6)}$ le triple ainsi construit.

2. $x = 0$. Dans ce cas, $u \neq 0 \neq \widetilde{x}$ et comme $u = -\widetilde{u}y$, $\widetilde{u} \neq 0 \neq y$. En posant $\sqrt{|uy|} = \nu^2$, $uy = \varepsilon\nu^2$ et $c = \varepsilon\frac{\widetilde{x}}{\nu^2}$; la transformation
   $$\begin{aligned} U &\leftarrow y^{-1}U \\ X &\leftarrow \varepsilon X \\ e_1 &\leftarrow \varepsilon\nu^{-1}e_1 \\ e_2 &\leftarrow \nu^{-1}(e_2 + \tfrac{\widetilde{u}}{u}\,f_2) \\ f_1 &\leftarrow \nu^{-1}f_1 \\ f_2 &\leftarrow \varepsilon\nu^{-1}f_2 \end{aligned}$$

livre la table
$$\begin{aligned}
[U\ f_1] &= e_2 \\
[U\ f_2] &= e_1 \\
[X\ f_1] &= U \\
[f_1\ f_2] &= U \\
[e_2\ f_1] &= cX
\end{aligned}$$

En renommant $f_2 \leftrightarrow e_2$, $U \leftrightarrow -U$, $e_1 \leftrightarrow -e_1$ et $X \leftrightarrow -X$, on a

$$\begin{aligned}
[U\ e_2] &= e_1 \\
[U\ f_1] &= -f_2 \\
[X\ f_1] &= e_1 \\
[e_2\ f_1] &= U \\
[f_1\ f_2] &= cX
\end{aligned}$$

c'est-à-dire la table (III) où

$$y = 0 = a = u = \widetilde{x} \quad \text{et} \quad x = c \neq 0$$

on est donc ramené au cas $x \neq 0$.

**La table 1 dans la preuve du lemme 6.1.1.7.**

La transformation
$$\begin{aligned}
U &\leftarrow U \\
e_1 &\leftarrow d^{-1}\ e_1 \\
e_2 &\leftarrow d^{-1}\ e_2 \\
f_1 &\leftarrow f_1 - x\ e_2 \\
f_2 &\leftarrow -f_2 + y\ e_2
\end{aligned}$$

livre la table
$$\begin{aligned}
[U\ e_2] &= e_1 \\
[U\ f_1] &= f_2 \\
[e_2\ f_1] &= U
\end{aligned}$$

Dans la base $\{e_1, f_2\} \cup \{e_2, f_1\}$ toute transformation de $\mathcal{P}$ du type

$$\varphi = \begin{pmatrix} \det(A) \cdot A & X \\ 0 & A \end{pmatrix} \qquad \text{où } A \in GL_2(\mathbb{R})$$

et $X \in \text{End}(\mathbb{R}^2)$, s'étend en un automorphisme de $(\mathcal{G}, \sigma)$. Réciproquement, la restriction à $\mathcal{P}$ d'un automorphisme de $(\mathcal{G}, \sigma)$ est de cette forme.

Dans cette même base, toute forme symplectique $\mathcal{K}$-invariante s'écrit

$$\Omega = \begin{pmatrix} 0 & S \\ -S & \xi \mathcal{J} \end{pmatrix}$$

où $S = {}^\tau S$ et $\xi \in \mathbb{R}$. On voit alors que l'on peut se ramener à

$$\Omega_\varepsilon = \begin{pmatrix} 0 & & 1 & 0 \\ & 0 & 0 & \varepsilon \\ -1 & 0 & & 0 \\ 0 & -\varepsilon & & 0 \end{pmatrix}$$

où $\varepsilon = \pm 1$. Nous noterons $t_\varepsilon^{4\ (7)}$ le triple ainsi construit.

**La table 2.(a) dans la preuve du lemme 6.1.1.7.**

Deux cas se présentent :
(a) $[\mathcal{K}\,\mathcal{P}]$ est central. On se ramène alors à

$$[U\ f_2] = e_2 \qquad [e_1\ f_2] = U$$

Dès lors, par un automorphisme de $(\mathcal{G}, \sigma)$, on ramène, dans la base $\{e_1, e_2, f_1, f_2\}$ de $\mathcal{P}$, toute forme symplectique $\mathcal{K}$-invariante à

$$\Omega = \begin{pmatrix} 0 & I \\ -I & 0 \end{pmatrix}$$

Nous noterons ce triple $t^{4\,(8)}$.
(b) $[\mathcal{K}\,\mathcal{P}]$ n'est pas central.
On se ramène alors à

$$[U\ f_2] = e_2 \qquad [e_2\ f_2] = \varepsilon U$$

où $\varepsilon = \pm 1$.

Un automorphisme de $(\mathcal{G}, \sigma)$ ramène la forme symplectique, dans la base $\{e_1, e_2, f_1, f_2\}$ à

$$\Omega = \begin{pmatrix} 0 & I \\ -I & 0 \end{pmatrix}$$

**Remarque**.
On voit que dans ce cas, $t = t_\varepsilon^2 \oplus t_0^2$ où $t_0^2$ est plat et où $t_\varepsilon^2$ est le triple de dimension 2 dont la table est donnée plus haut.

**La table 2.(b) du lemme 6.1.1.7.**

On se ramène à
$$\begin{array}{rcl} [U\ f_1] & = & e_1 \\ [U\ f_2] & = & e_2 \\ [f_1\ f_2] & = & U \end{array}$$

et donc au cas $t_\varepsilon^{4\,(7)}$.

## 6.2 Classification

**Proposition 6.2.1.** *La liste des T.S.S. de dimension 2 est*

- $t_0^2 =$ *le triple plat*
- $t_\varepsilon^2 (\varepsilon = \pm 1) =$ *les triples décrits page 107*
- *Les triples associés aux E.S.S. simples:* $SU(1,1)\big/{SO(2)}$, $SU(1,1)\big/{\mathbb{R}}$, $SU(2)\big/{SO(2)}$

*Nous les noterons respectivement* $t_D$, $t_{H^1}$, $t_{S^2}$

**Proposition 6.2.2.** *La liste des T.S.S. résolubles de dimension 4 est*

- $t_0^4$: *le triple plat*
- $t_{1,\varepsilon,x}^{4\,(1)}, t_{2,\varepsilon,0,x}^{4\,(1)}, t_{2,0,\varepsilon,x}^{4\,(1)}, t_{3,\varepsilon,x}^{4\,(1)}(\varepsilon = \pm 1, x \in \mathbb{R}) =$ *les triples décrits page 87 (Proposition 5.1.2.3.)*

- $t_{1,a,x}^{4\ (2)}, t_{2,\varepsilon,\varepsilon',a,x}^{4\ (2)}(\varepsilon = \pm 1, \varepsilon' = \pm 1, x \in \mathbb{R}, a \in \mathbb{R}_0) = $ les triples décrits page 90 (Proposition 5.2.2.2.)
- $t_{\varepsilon,\varepsilon',\eta}^{4\ (3)}(\varepsilon = \pm 1, \varepsilon' = \pm 1, \eta = 0, 1) = $ les triples décrits page 103
- $t_{\varepsilon,\alpha}^{4\ (4)}(\varepsilon = \pm 1, \Omega = \delta\alpha) = $ les triples décrits page 104
- $t_{\varepsilon,\alpha,\eta}^{4\ (5)}(\varepsilon = \pm 1; \eta = -1, 0, 1) = $ les triples décrits page 105
- $t_{\varepsilon}^{4\ (6)}(\varepsilon = \pm 1) = $ les triples décrits page 105-106
- $t_{\varepsilon}^{4\ (7)}(\varepsilon = \pm 1) = $ les triples décrits page 107
- $t_{\varepsilon}^{4\ (8)} = $ les triples décrits page 107
- $t_0^2 \oplus t_\varepsilon^2$

    **Remarques.**

(i) Le nombre de structures symétriques non-isomorphes sous-jacentes aux E.S.S. simplement connexes resolubles de dimension 4 est fini (et vaut 25).

(ii) Les triples $t_\varepsilon^2 \oplus t_{\varepsilon'}^2$ n'apparaissent pas explicitement dans la liste de la proposition 6.2.2.. Ils sont, en fait, isomorphes aux triples $t_\cdot^{4\ (1)}$ pour la valeur nulle du paramètre $x$.

**Proposition 6.2.3.** *La liste des T.S.S. ni résolubles, ni semi-simples de dimension 4 est*

- $t_0^2 \oplus t_{S^2}$
- $t_0^2 \oplus t_{H^1}$
- $t_0^2 \oplus t_D$
- *Les triples associés aux E.S.S. que constituent*
    1. le fibré cotangent à la sphère $S^2 = SU(2)\big/SO(2)$
    2. le fibré cotangent à l'hyperboloïde à une nappe $H^1 = SU(1,1)\big/\mathbb{R}$
    3. le fibré cotangent au disque $D = SU(1,1)\big/SO(2)$

**Proposition 6.2.4.** *La liste des T.S.S. simples de dimension 4 est :*

| $\mathcal{G}$ | $\mathcal{K}$ |
|---|---|
| $su(3)$ | $su(2) \oplus so(2)$ |
| $su(1,2)$ | $su(2) \oplus so(2)$ |
| $su(1,2)$ | $su(1,1) \oplus so(2)$ |
| $sl(3,\mathbb{R})$ | $sl(2,\mathbb{R}) \oplus \mathbb{R}$ |
| $sl(2,\mathbb{C})$ | $\mathbb{C}$ |

# Références